\documentclass[twoside,12pt,a4paper]{report}
\usepackage{hyperref}
\usepackage{graphicx}

\newif\iftth

\input{mydef}

\usepackage{makeidx}
\makeindex
  \makeatletter
  \newcommand{\eqhyperref}[2][]{\hyperref[#1]{#2}}

  \newcommand{\extref}[3]{\hyperref[#1!!#2]{#3}}
  \newcommand{\intro}[2][\relax]{\ifx#1\relax\index{#2}\else\index{#1}\fi
    \emph{#2}}
  \usepackage[extnum]{myart2k}

\newcommand{\introind}[2]{\intro{#1}\index{#2}}

\PassOptionsToPackage{pdfpagescrop={53 400 300 700}%
}{hyperref}
\providecommand{\matr}[4]{
\ensuremath{  \left(\!\begin{array}{cc}
      #1 & #2\\
      #3 & #4
  \end{array}\!\right)}
}
\providecommand{\mvec}[1]{\mathbf{#1}}
\providecommand{\tr}{\object{tr}}

\ifundefined{eoe}
   \DeclareMathSymbol{\eoe}      {\mathord}{AMSa}{"06}
\fi

     \theorembodyfont{\upshape}
\newtheorem{examplea}{\addtocounter{thm}{1}Example}
\newtheorem{exampleb}{Example}

\providecommand{\oper}[1]{\mathcal{#1}}
\newcommand{\vecbf}[1]{\mathbf{#1}}
\newcommand{\n}[1]{\mathsf{#1}}
\newcommand{\TSpace}[2]{\ensuremath{ { \widetilde{\mathbb{#1}}^{#2}} }}

\input{cyr.fd}

\usepackage{makeidx}
\makeindex

  \providecommand{\hyperlink}[2][\relax]{#2}

  \providecommand{\eqhyperref}[2][]{#2}
  \providecommand{\hypertarget}[2][]{#2}
  \providecommand{\pdfbookmark}[3][0]{}
  
\begin{document}
\title[]{Spaces of Analytical Functions and Wavelets\\Lecture Notes}

\author[Vladimir V. Kisil]%
{\href{http://maths.leeds.ac.uk/~kisilv/}{Vladimir V. Kisil}}

\maketitle
\begin{abstract} 
  This is (raw) lecture notes of the course read on 6th European intensive
  course on Complex Analysis (Coimbra, Portugal) in 2000.

  Our purpose is to describe a general framework for generalizations of
  the complex analysis. As a consequence a classification scheme for
  different generalizations is obtained.  

  The framework is based on
  wavelets (coherent states) in Banach spaces generated by
  ``admissible'' group representations. Reduced wavelet transform
  allows naturally describe in abstract term main objects of an
  analytical function theory: the Cauchy integral formula, the Hardy
  and Bergman spaces, the Cauchy-Riemann equation, and the Taylor expansion.

  Among considered examples are classical analytical function theories 
  (one complex variables, several complex variables, Clifford
  analysis, Segal-Bargmann space) as well as new function theories
  which were developed within our framework (function theory of
  hyperbolic type, Clifford version of Segal-Bargmann space).

  We also briefly discuss applications to the operator theory
  (functional calculus) and quantum mechanics.
\AMSMSC{30G30}{42C40, 43A85, 46H30, 47A13, 81R30, 81R60}

\end{abstract}

\noindent \textbf{Address}: \\
~\qquad Department of Pure Mathematics,\\
~\qquad University of Leeds,\\
~\qquad Leeds LS2\,9JT,\\
~\qquad UK\\

\noindent\textbf{Email}:\\
~\qquad \href{mailto:kisilv@maths.leeds.ac.uk}%
{\texttt{kisilv@maths.leeds.ac.uk}}\\

\noindent\textbf{URL}: \\
~\qquad \href{http://maths.leeds.ac.uk/~kisilv/}%
{\texttt{http://maths.leeds.ac.uk/\~{}kisilv/}}\\

\tableofcontents

\chapter[Generalizations of Complex Analysis]{Different 
Generalizations of Complex Analysis}

\section{Introduction}
The classic heritage of complex analysis is contested between several
complex variables theory and hypercomplex analysis. The first one
was founded long ago by Cauchy and Weierstrass themselves and sometime
thought to be the only crown-prince. The hypercomplex analysis is not
a single theory but a family of related constructions discovered quite
recently~\cite{BraDelSom82,DelSomSou92,GuerSpross90} 
(and rediscovered up to now)
under hypercomplex framework.

Such a variety of theories puts the question on their classification.
One could dream about a Mendeleev-like periodic table for hypercomplex
analysis, which clearly explains properties of different theories,
relationship between them and indicates how many blank cells are
waiting for us. Moreover, because hypercomplex analysis is the
recognized background for \introind{classic
  mechanics}{mechanics!quantum} and \introind{quantum 
  physics}{mechanics!quantum} theories like the
Maxwell and Dirac equations, such a
table could play the role of \emph{the Mendeleev table for elementary
particles and fields}. We will return to this metaphor
and find it is not very superficial.

To make a step in the desired direction we should specify the notion
of \intro{function theory} and define the concept of \intro{essential
difference}. Probably many people agree that 
\begin{defn}
  The \introind{core of complex analysis}{complex!analysis} consists
  of
  \begin{enumerate}
  \item\label{it:first} The Cauchy-Riemann equation and complex
    derivative $\frac{\partial }{\partial z}$;
  \item The Cauchy theorem;
  \item The Cauchy integral formula;
  \item The Plemeli-Sokhotski formula;
  \item\label{it:last} The Taylor and Laurent series.
  \end{enumerate}
\end{defn}
Any development of several complex variables theory or hypercomplex
analysis is beginning from analogies to these notions and results. 
Thus we adopt the following
\begin{defn}
A \intro{function theory} is a collection of notions and results, which
includes at least analogies of~\ref{it:first}--\ref{it:last}.
\end{defn}
Of course the definition is more philosophical than mathematical. For
example, the understanding of an \intro{analogy} and especially the
\introind{right analogy}{analogy!right} usually generates many
disputes.

Again as a first approximation we propose the following
\begin{defn}
Two function theories is said to be \introind{similar}{function
theory!similar} if there is a correspondence between their objects
such that analogies of~\ref{it:first}--\ref{it:last} in one theory
follow from their counterparts in another theory. Two function
theories are \intro{essentially different}{function
  theory!(essentially) different} if they are not similar.
\end{defn}
Unspecified ``correspondence'' should probably be a linear map 
and we will look for its meaning soon. It is
clear that the \intro{similarity} is an equivalence relation
and we are looking for quotient sets with respect to it.

The layout is following. 
In Subsection~\ref{ss:factorization} the 
classic scheme of hypercomplex analysis is discussed and a possible 
variety of function theories appears. But we will see in 
Subsection~\ref{ss:reduce} that not all of them are very different.
Connection between group representations and (hyper)complex analysis 
is presented in Section~\ref{se:towards}. 
It could be a base for classification of essentially different 
theories.

\section{Factorizations of the Laplacian}\label{ss:factorization}
In the next Section we repeat shortly the
scheme of development of Clifford analysis as it could be found
in~\cite{BraDelSom82,DelSomSou92}. We examine different options
arising on this way and demonstrate that some differences are only
apparent not essential.

We would like to see how the contents of~\ref{it:first}--\ref{it:last}
could be realized in a function theory. We are interested in function 
theories defined in $\Space{R}{d}$. The Cauchy theorem and
integral formula clearly indicates that the behavior of functions
inside a domain should be governed by their values on the boundary.
Such a property is particularly possessed by solutions to the
\intro{second order elliptic differential operator} $P$
\begin{displaymath}
P(x,\partial_x)=\sum_{i,j=1}^{d} a_{ij}(x)\partial_i \partial_j +
\sum_{i=1}^{d} b_i(x) \partial_i +c(x)
\end{displaymath}
with some special properties. Of course, the principal example is the
Laplacian
\begin{equation}\label{eq:Laplace}
\Delta=\sum_{i=1}^{d} \frac{\partial ^2}{\partial x_i^2}.
\end{equation}

\begin{punct}
  \begin{enumerate}
  \item \label{it:equation}\emph{Choice of different operators} (for
    example, the Laplacian or the Helmholtz operator) is the first
    option which brings the variety in the family of hypercomplex
    analysis.
  \end{enumerate}
  
  The next step is called \intro{linearization}. Namely we are looking
  for two (possibly coinciding) \intro{first order differential
  operators} $D$ and $D'$ such that
  \begin{displaymath}
    DD'=P(x,\partial_x).
  \end{displaymath}
  The Dirac motivation to do that is to ``look for an equation linear in
  in time derivative $\frac{\partial }{\partial t}$, because the
  Schr\"odinger equation is''. From the function theory point of view
  the Cauchy-Riemann operator should be linear also. But the most
  important gain of the step is an introduction of the Clifford algebra.
  For example, to factorize the Laplacian~\eqref{eq:Laplace} we put
  \begin{equation}\label{eq:dirac}
    D=\sum_{i=1}^d e_i \partial_i
  \end{equation}
  where $e_i$ are the \intro{Clifford algebra} generators:
  \begin{equation}\label{eq:anti-comm}
    e_i e_j + e_j e_i = 2\delta _{ij}, \qquad 1\leq i,j\leq d.
  \end{equation}
  \begin{enumerate}\addtocounter{enumi}{1}
  \item \label{it:algebra}\emph{Different linearizations of a second
      order operator} multiply the spectrum of theories.
  \end{enumerate}
\end{punct}
Mathematicians and physicists are looking up to now new factorization
even for the Laplacian\nocite{Keller93,ShaVas94a}. The essential
uniqueness of such factorization was already felt by Dirac himself but
it was never put as a theorem. So the idea of the \emph{genuine}
factorization becomes the philosophers' stone of our times.

After one made a choice~\ref{it:equation} and~\ref{it:algebra} the
following turns to be a routine. The equation
\begin{displaymath}
D' f(x)=0,
\end{displaymath}
plays the role of the Cauchy-Riemann equation.
Having a fundamental solution $F(x)$
to the operator $P(x,\partial _x)$ the Cauchy integral kernel defined
by
\begin{displaymath}
E(x)=D' F(x)
\end{displaymath}
with the property $DE(x)=\delta (x)$. Then the Stocks theorem implies
the Cauchy theorem and Cauchy integral formula. A decomposition of the
Cauchy kernel of the form
\begin{displaymath}
C(x-y)=\sum_{\alpha } V_\alpha (x) W_\alpha (y),
\end{displaymath}
where $V_\alpha (x)$ are some polynomials, yields via integration over
the ball the Taylor and Laurent series\footnote{Not all such 
decompositions give interesting series. The scheme from 
Section~\ref{se:towards} gives a selection rule to distinguish them.}. 
In such a way the program-minimum~\ref{it:first}--\ref{it:last} could 
be accomplished.

Thus all possibilities to alter function theory concentrated
in~\ref{it:equation} and \ref{it:algebra}. Possible universal algebras
arising from such an approach were investigated by
\person{F.~Sommen}~\cite{Sommen95a}. In spite of the apparent wide
selection, for operator $D$ and $D'$ with constant coefficients it was
found ``nothing dramatically new''~\cite{Sommen95a}:
\begin{quotation}
Of course one can study all these algebras and prove theorems or work
out lots of examples and representations of universal algebras. But in
the constant coefficient case the most important factorization seems
to remain the relation $\Delta=\sum x_j^2$, i.e., the one leading to
the definition of the Clifford algebra.
\end{quotation}
We present an example that there is no dramatical news not only on the
level of universal algebras but also for function theory (for the 
constant coefficient case). We will return to non
constant case in Section~\ref{se:towards}.

\section{Example of Connection}\label{ss:reduce}
We give a short example of similar theories with explicit
connection between them. The full account could be found
in~\cite{Kisil95c},
another example was considered in~\cite{Ryan95}.

Due to physical application we will consider equation
\begin{equation}\label{eq:mass}
\frac{\partial f}{\partial y_0}=(\sum_{j=1}^n e_j\frac{\partial
}{\partial y_j}+M)f,
\end{equation}
where $e_j$ are generators~\eqref{eq:anti-comm} of the Clifford 
algebra and $M=M_\lambda $ is an operator of
multiplication from the {\em right\/}-hand side by the Clifford number
$\lambda $. Equation~\eqref{eq:mass} is known in \introind{quantum
mechanics}{mechanics!quantum} as the
{\em \introind{Dirac
equation}{equation!Dirac} for a particle with a non-zero rest
mass\/}~\cite[\S 20]{BerLif82}, \cite[\S 6.3]{BogShir80} and
\cite{Kravchenko95a}. We will specialize our results for the case
$M=M_\lambda$, especially for the simplest (but still important!) case
$\lambda\in\Space{R}{}$.
\begin{thm}
The function $f(y)$ is a solution to the equation
\begin{displaymath}
\frac{\partial f}{\partial y_0}=(\sum_{j=1}^n e_j\frac{\partial
}{\partial y_j}+M_1)f
\end{displaymath}
 if and only if the function
\begin{displaymath}
g(y)=e^{y_0 M_2} e^{-y_0 M_1} f(y)
\end{displaymath}
is a solution to the equation
\begin{displaymath}
\frac{\partial g}{\partial y_0}=(\sum_{j=1}^n e_j\frac{\partial
}{\partial y_j}+M_2)g,
\end{displaymath}
where $M_1$ and $M_2$ are bounded operators commuting with $e_j$.
\end{thm}
\begin{cor}\label{co:mass}
The function $f(y)$ is a solution to the equation~\eqref{eq:mass}
 if and only if the function
$e^{y_0 M}f(y)$
is a solution to the generalized Cauchy-Riemann
equation~\eqref{eq:dirac}.

In the case $M=M_\lambda$ we have $e^{y_0 M_\lambda}f(y)=f(y)e^{y_0
\lambda}$ and if $\lambda\in \Space{R}{}$ then $e^{y_0
M_\lambda}f(y)=f(y)e^{y_0 \lambda}=e^{y_0 \lambda}f(y)$.
\end{cor}

In this Subsection we construct a function theory (in the sense
of~\ref{it:first}--\ref{it:last}) for $M$-solutions of
the generalized Cauchy-Riemann operator based on Clifford analysis
and Corollary~\ref{co:mass}.

The set of solutions to~\eqref{eq:dirac} and~\eqref{eq:mass} in a nice
domain $\Omega$ will be denoted by
$\algebra{M}(\Omega)=\algebra{M}_0(\Omega)$ and
$\algebra{M}_M(\Omega)$ correspondingly. In the case $M=M_\lambda$ we
use the notation
$\algebra{M}_\lambda(\Omega)=\algebra{M}_{M_\lambda}(\Omega)$ also. We
suppose that all functions from $\algebra{M}_\lambda(\Omega)$ are
continuous in the closure of $\Omega$. Let
\begin{equation}\label{eq:cauchy-ker}
E(y-x)=
\frac{\Gamma(\frac{n+1}{2})}{2\pi^{(n+1)/2}}\,
\frac{\overline{y-x}}{\modulus{y-x}^{n+1}}
\end{equation}
be the Cauchy kernel~\cite[p.~146]{DelSomSou92}  and
\begin{displaymath}
d\sigma=\sum_{j=0}^n (-1)^j e_j dx_0 \wedge \ldots \wedge[dx_j] \wedge
\ldots \wedge dx_m.
\end{displaymath}
be the differential form of the ``oriented surface
element''~\cite[p.~144]{DelSomSou92}. Then for any
$f(x)\in\algebra{M}(\Omega)$ we have the Cauchy integral
formula~\cite[p.~147]{DelSomSou92}
\begin{displaymath}
\int_{\partial \Omega} E(y-x)\,d\sigma_y\,
f(y)=\left\{\begin{array}{cl}
f(x),& x\in\Omega\\
0,& x\not\in\bar{\Omega}
\end{array}.\right.
\end{displaymath}

\begin{thm}[Cauchy's Theorem]
Let $f(y)\in \algebra{M}_M(\Omega)$. Then
\begin{displaymath}
\int_{\partial \Omega} d\sigma_y\,e^{-y_0 M}f(y)=0.
\end{displaymath}
Particularly, for $f(y)\in \algebra{M}_\lambda(\Omega)$ we have
\begin{displaymath}
\int_{\partial \Omega} d\sigma_y\,f(y)e^{-y_0 \lambda}=0,
\end{displaymath}
and
\begin{displaymath}
\int_{\partial \Omega} d\sigma_ye^{-y_0 \lambda}\,f(y)=0,
\end{displaymath}
if $\lambda\in\Space{R}{} $.
\end{thm}

\begin{thm}[Cauchy's Integral Formula]
Let $f(y)\in \algebra{M}_M(\Omega)$. Then
\begin{equation}\label{eq:m-cauchy}
e^{x_0 M}\int_{\partial \Omega} E(y-x)\,d\sigma_y\,
e^{-y_0 M}f(y)=\left\{\begin{array}{cl}
f(x),& x\in\Omega\\
0,& x\not\in\bar{\Omega}
\end{array}.\right.
\end{equation}
Particularly, for $f(y)\in \algebra{M}_\lambda(\Omega)$ we have
\begin{displaymath}
\int_{\partial \Omega} E(y-x)\,d\sigma_y\,
f(y)e^{(x_0-y_0) \lambda}=\left\{\begin{array}{cl}
f(x),& x\in\Omega\\
0, &x\not\in\bar{\Omega}
\end{array}.\right.
\end{displaymath}
and
\begin{displaymath}
\int_{\partial \Omega} E(y-x)e^{(x_0-y_0) \lambda}\,d\sigma_y\,
f(y)=\left\{\begin{array}{cl}
f(x),& x\in\Omega\\
0,& x\not\in\bar{\Omega}
\end{array}.\right.
\end{displaymath}
if $\lambda\in\Space{R}{} $.
\end{thm}
It is hard to expect that formula~\eqref{eq:m-cauchy} may be rewritten
as
\begin{displaymath}
\int_{\partial \Omega} E'(y-x)\,d\sigma_y\,
f(y)=\left\{\begin{array}{cl}
f(x),& x\in\Omega\\
0,& x\not\in\bar{\Omega}
\end{array}\right.
\end{displaymath}
with a simple function $E'(y-x)$.

Because an application of the bounded operator $e^{y_0 M}$ does not
destroy uniform convergency of functions we obtain
(cf.~\cite[Chap.~II, \S~0.2.2, Theorem~2]{DelSomSou92})
\begin{thm}[Weierstrass' Theorem]
Let $\{f_k\}_{k\in\Space{N}{}}$ be a sequence in
$\algebra{M}_M(\Omega)$, which converges uniformly to $f$ on each
compact subset $K\in \Omega$. Then
\begin{enumerate}
\item $f\in \algebra{M}_M(\Omega)$.
\item For each multi-index
$\beta=(\beta_0,\ldots,\beta_m)\in\Space{N}{n+1}$, the sequence
$\{\partial ^\beta f_k\}_{k\in\Space{N}{}} $ converges uniformly on
each compact subset $K\in \Omega$ to $\partial ^\beta f$.
\end{enumerate}
\end{thm}

\begin{thm}[Mean Value Theorem]
Let $f\in\algebra{M}_M(\Omega)$. Then for all $x\in \Omega$ and $R>0$
such that the ball $\Space{B}{}(x,R)\in\Omega $,
\begin{displaymath}
f(x)= e^{x_0 M}
\frac{(n+1)\Gamma(\frac{n+1}{2})}{2R^{n+1}\pi^{(n+1)/2}}
\int_{\Space{B}{}(x,R)} e^{-y_0 M} f(y)\,dy.
\end{displaymath}
\end{thm}

Such a reduction of theories could be pushed even 
future~\cite{Kisil95c} up to the notion of 
hypercomplex differentiability~\cite{Malonek93}, but we will 
stop here.  

\section[Analysis and Group Representations]{Hypercomplex 
Analysis and Group
Representations --- Towards a Classification}\label{se:towards}
To construct a classification of non-equivalent objects one could
use their groups of symmetries. Classical example is Poincar\'e's proof
of bi-holomorphic non-equivalence of the unit ball and polydisk via
comparison their groups of bi-holomorphic automorphisms. To employ 
this approach we need a construction of hypercomplex analysis from its
symmetry group. The following scheme will be main theme of this
Course.

Let $G$ be a group which acts via transformation of a closed domain
$\bar{\Omega}$. Moreover, let $G: \partial \Omega\rightarrow \partial
\Omega$ and $G$ act on $\Omega$ and $\partial \Omega$ transitively.
Let us fix a point $x_0\in \Omega$ and let $H\subset G$ be a
stationary subgroup of point $x_0$. Then domain $\Omega$ is naturally
identified with the  homogeneous space $G/H$. Till the moment we do
not request anything untypical. Now let 
\begin{itemize}
\item\emph{there exist a $H$-invariant measure $d\mu$ on $\partial
    \Omega$}.
\end{itemize}
We consider the Hilbert space $\FSpace{L}{2}(\partial
\Omega, d\mu)$. Then geometrical transformations of $\partial \Omega$
give us the representation $\pi$ of $G$ in $\FSpace{L}{2}(\partial
\Omega, d\mu)$.
 Let $f_0(x)\equiv 1$ and $\FSpace{F}{2}(\partial
\Omega, d\mu)$ be the closed liner subspace of $\FSpace{L}{2}(\partial
\Omega, d\mu)$ with the properties:
\begin{enumerate}
\item\label{it:begin} $f_0\in \FSpace{F}{2}(\partial \Omega, d\mu)$;
\item $\FSpace{F}{2}(\partial \Omega, d\mu)$ is $G$-invariant;
\item\label{it:end} $\FSpace{F}{2}(\partial \Omega, d\mu)$ is
  $G$-irreducible, or $f_0$ is cyclic in $\FSpace{F}{2}(\partial
  \Omega, d\mu)$. 
\end{enumerate}
The \intro{standard wavelet transform} $W$ is defined by
\begin{displaymath}
W: \FSpace{F}{2}(\partial \Omega, d\mu) \rightarrow
\FSpace{L}{2}(G): f(x) \mapsto
\widehat{f}(g)=\scalar{f(x)}{\pi(g)f_0(x)}_{\FSpace{L}{2}(\partial
\Omega,d\mu) }
\end{displaymath}
Due to the property $[\pi(h)f_0](x)=f_0(x)$, $h\in H$ and 
identification $\Omega\sim G/H$ it could be translated to the embedding:
\begin{equation}\label{eq:cauchy}
\widetilde{W}: \FSpace{F}{2}(\partial \Omega, d\mu) \rightarrow
\FSpace{L}{2}(\Omega): f(x) \mapsto
\widehat{f}(y)=\scalar{f(x)}{\pi(g)f_0(x)}_{\FSpace{L}{2}(\partial
\Omega,d\mu) },  
\end{equation}
where $y\in\Omega$ for some $ h\in H$. The imbedding~\eqref{eq:cauchy} 
is \emph{an abstract analog of the} \intro{Cauchy integral formula}. 
Let functions $V_\alpha $ be the \intro{special functions} generated by 
the representation of $H$. Then the decomposition of 
$\widehat{f}_0(y)$ by $V_\alpha $ gives us the Taylor series.

The scheme is inspired by the following interpretation of complex 
analysis.
\begin{example}
Let the domain $\Omega$ be the unit disk $\Space{D}{}$, $\partial 
\Space{D}{}=\Space{S}{}$. We select
the group $SL(2,\Space{R}{})\sim SU(1,1)$ acting on $\Space{D}{} $
via the fractional-linear transformation:
\begin{displaymath}
\matr{a}{b}{c}{d}: z\mapsto \frac{az+b}{cz+d}.
\end{displaymath}
We fix $x_0=0$. Then its stationary group is $U(1)$ of rotations of 
$\Space{D}{}$. Then the Lebesgue measure on $\Space{S}{} $ is 
$U(1)$-invariant. We obtain $\Space{D}{}\sim SL(2,\Space{R}{} )/U(1)$. 
The subspace of $\FSpace{L}{2}(\Space{S}{},dt)$ satisfying 
to~\ref{it:begin}--\ref{it:end} is the Hardy space. The wavelets 
transform\eqref{eq:cauchy} give exactly the Cauchy formula.
The proper functions of $U(1)$ are exactly $z^n$, which provide the 
base for the Taylor series. The Riemann mapping theorem allows to 
apply the scheme to any connected, simply-connected domain.
\end{example}

The conformal group of the M\"obius transformations plays the same 
role in Clifford analysis. One usually says that the conformal group 
in $\Space{R}{n}$, $n>2$ is
not so rich as the conformal group in $\Space{R}{2}$.
Nevertheless, the conformal covariance has many applications in
Clifford analysis~\cite{Cnops94a,Ryan95b}.
Notably,
groups of conformal mappings of open unit balls $\Space{B}{n} \subset
\Space{R}{n}$ \emph{on}to itself are similar for all $n$ and
as sets can be
parametrized by the product of \Space{B}{n} itself and the group
of isometries of its boundary \Space{S}{n-1}.
\begin{thm}{\cite{Kisil95i} }\label{th:ball}
Let $a\in\Space{B}{n}$, $b\in\Gamma_n$ then the M\"obius
transformations of
the form
\begin{displaymath}
\phi_{(a,b)}=\matr{b}{0}{0}{{b}^{*-1}}\matr{1}{-
a}{{a}^*}{-1}=\matr{b}{-ba}{{b}^{*-1}a^*}
{-{b}^{*-1}},
\end{displaymath}
constitute the group $B_{n}$ of conformal mappings of the open unit
ball $\Space{B}{n}$ onto itself. $B_{n}$ acts on
$\Space{B}{n}$ transitively.
Transformations of the form $\phi_{(0,b)}$  constitute a
subgroup isomorphic to $\object[(n)]{O}$. The homogeneous space
$B_{n}/\object[(n)]{O}$ is isomorphic as a set to
$\Space{B}{n}$. Moreover:
\begin{enumerate}
\item $\phi_{(a,1)}^2=1$ identically on $\Space{B}{n}$
($\phi_{(a,1)}^{-1}=\phi_{(a,1)}$).
\item $\phi_{(a,1)}(0)=a$, $\phi_{(a,1)}(a)=0$.
\end{enumerate}
\end{thm}

Obviously, conformal mappings preserve the space of null solutions to
the \intro{Laplace operator}~\eqref{eq:Laplace} and null solutions the
\intro{Dirac operator}~\eqref{eq:dirac}. The
group $B_{n}$ is sufficient for construction of the Poisson and the 
Cauchy integral representation of harmonic functions and Szeg\"o and 
Bergman projections in Clifford analysis by the 
formula~\cite{Kisil95d}
\begin{equation}\label{eq:reproduce}
K(x,y)=c\int_G
[\pi_g f](x) \overline{[\pi_g f](y)}\,dg,
\end{equation}
where $\pi_g$ is an irreducible unitary
square integrable representation of a group $G$,
$f(x)$ is an arbitrary non-zero function, and $c$ is a constant.

The scheme gives a correspondence between \intro{function theories} and 
\intro{group representations}. The last are rather well studded and 
thus such a connection could be a foundation for a classification of 
function theories. Particularly, the \emph{constant coefficient} 
function theories in the sense of \person{F.~Sommen}\cite{Sommen95a} 
corresponds to the groups acting only on the function domains in the 
Euclidean space. Between such groups the Moebius transformations play 
the leading role. On the contrary, the \emph{variable coefficient} 
case is described by groups acting on the function space in the 
non-point sense (for example, combining action on the functions domain 
and range, see~\cite{Kisil94e}). The set of groups of the second kind 
should be more profound.

\begin{rem}
It is known that many results in real analysis~\cite{McIntosh95a} several 
variables theory~\cite{MiSha95} could be obtained or even explained via 
hypercomplex analysis. One could see roots of this phenomenon in 
relationships between groups of geometric symmetries of two theories: the 
group of hypercomplex analysis is wider.
\end{rem}

Returning to our metaphor on the Mendeleev table we would like recall
that it began as linear ordering with respect to atomic masses but
have received an explanation only via representation theory of the
rotation group.

\chapter[Wavelets and Analytic Functions]{Group 
Representations, Wavelets and Analytic Spaces of Functions}

\section{Introduction}
\label{sec:introduction1}

The purpose of this Lecture to introduce the appropriate language of
\intro{coherent states} and \intro{wavelet transform}. We suppose some
knowledge about groups and their representations. The appropriate
material is included in Appendix~\ref{sec:groups-homog-spac}
and~\ref{sec:elem-repr-theory}. We will begin from 
the standard constructions of coherent states (wavelets) in a Hilbert
space (section~\ref{sec:wavel-hilb-spac}) and then will construct an
appropriate generalization for Banach spaces
(section~\ref{sec:wavel-banach-spac}).

Wavelet transform considered here is an important example of the
interesting object called \introind{token} in \cite{Kisil01b}. Tokens
are kernels of intertwining operators between actions of two
cancellative semigroups.

\section{Wavelets in Hilbert Spaces}
\label{sec:wavel-hilb-spac}

\subsection{Wavelet Transform and Coherent States}

We agree with a reader if he/she is not satisfied by the last short
proof and would like to see a more detailed account how the core of
complex analysis could be reconstructed from representation theory of
$\SL$.  We present an abstract scheme, which also could be applied to
other analytic function theories, see last two lectures and 
\cite{CnopsKisil97a,Kisil97c}.  We
start from a dry construction followed in the next Section by classic
examples, which will justify our usage of personal names.

Let $ X $ be a topological space and let $G$ be a group that acts $G:
{X}\rightarrow {X}$ as a transformation $g: x \mapsto g \cdot x$ from
the left, i.e., $g_1 \cdot(g_2 \cdot x)=(g_1 g_2)\cdot x$. Moreover,
let $G$ act on $X$ transitively. Let there exist a measure $dx$ on $
X$ such that a representation $\pi(g): f(x) \mapsto m(g,x) f(
g^{-1}\cdot x)$ (with a function $m(g,h)$) is unitary with respect to
the scalar product $\scalar{f_1(x)}{f_2(x)}_{\FSpace{L}{2}(X ) } =
\int_{ X} f_1(x) \bar{f}_2(x) \,d(x)$, i.e.,
\begin{displaymath} 
\scalar{[\pi(g)f_1](x)}{[\pi(g)f_2](x)}_{\FSpace{L}{2}( X ) }
= \scalar{f_1(x)}{f_2(x)}_{\FSpace{L}{2}( X ) }\quad \forall f_1, 
f_2 \in \FSpace{L}{2}(X).
\end{displaymath}
We consider the Hilbert space $\FSpace{L}{2}( X )$ where
representation $\pi(g)$ acts by unitary operators.
\begin{rem}
  It is well known that the most developed part of representation
  theory consider unitary representations in Hilbert spaces. By this
  reason we restrict our attention to Hilbert spaces of analytic
  functions, the area usually done by means of the functional analysis
  technique. We also assume that our functions are complex valued and
  this is sufficient for examples explicitly considered in the present
  paper. However the presented scheme is working also for vector
  valued functions and this is the natural environment for Clifford
  analysis~\cite{BraDelSom82}, for example. One also could start from
  an abstract Hilbert space $H$ with no explicit realization as $
  \FSpace{L}{2}(X)$ given.
\end{rem}  

Let $H$ be a closed compact\footnote{While the compactness will be
  explicitly used during our abstract consideration, it is not crucial
  in fact. One could make a trick
  for non-compact $H$~\cite{Kisil98a}.} 
subgroup of $G$ and let $f_0(x)$ be such a
function that $H$ acts on it as the multiplication
\begin{equation} \label{eq:homogenious}
[\pi(h)f_0](x)=\chi(h) f_0(x) \qquad , \forall h\in H,
\end{equation}
by a function $\chi(h)$, which is a character of $H$ i.e., $f_0(0)$ is
a common eigenfunction for all operators $\pi(h)$. Equivalently
$f_0(x)$ is a common eigenfunction for operators corresponding under
$\pi$ to a basis of the Lie algebra of $H$.  Note also that $
\modulus{\chi(h)}^2=1 $ because $\pi$ is unitary.  $f_0(x)$ is called
\intro{vacuum vector} (with respect to subgroup $H$).  We introduce the
$\FSpace{F}{2}(X)$ to be the closed liner subspace of
$\FSpace{L}{2}(X)$ uniquely defined by the conditions:
\begin{enumerate}
\item\label{it:begin1} $f_0\in \FSpace{F}{2}( X )$;
\item $\FSpace{F}{2}( X )$ is $G$-invariant;
\item\label{it:end1} $\FSpace{F}{2}( X )$ is $G$-irreducible, or $f_0$
  is cyclic in $\FSpace{F}{2}(\partial  \Omega, d\mu)$.
\end{enumerate}
Thus restriction of $\pi$ on $ \FSpace{F}{2}(X) $ is an irreducible
unitary representation.

The \intro{wavelet transform}\footnote{The subject of coherent states
  or wavelets have been arising many times in many \emph{applied}
  areas and the author is not able to give a comprehensive history and
  proper credits. One could mention important
  books~\cite{Daubechies92,KlaSkag85,Perelomov86}. We give our
  references by recent paper~\cite{Kisil95a}, where applications to
  \emph{pure} mathematics were considered.} $\oper{W}$ could be
defined for square-integral representations $\pi$ by the
formula
\begin{eqnarray}
\oper{W}&:& \FSpace{F}{2}( X ) \rightarrow
\FSpace{L}{ \infty }(G) \nonumber \\
&:& f(x) \mapsto
\widetilde{f}(g)=\scalar{f(x)}{\pi(g)f_0(x)}_{\FSpace{L}{2}(
X ) } \label{eq:wavelets}
\end{eqnarray}
The principal advantage of the wavelet transform $\oper{W}$ is that it
express the representation $\pi$ in geometrical terms. Namely it
\emph{intertwins} $\pi$ and left regular representation $ \lambda $ on
$G$:
\begin{equation} \label{eq:g-inter}
[\lambda_g\oper{W} f](g')
=[\oper{W}f](g^{-1}g')
= \scalar{f}{\pi_{g^{-1}g'}f_0}
= \scalar{\pi_g f}{\pi_{g'}f_0}
=[\oper{W}\pi_g f](g'),
\end{equation}
i.e., $\lambda \oper{W} = \oper{W} \pi$. Another important feature of
$W$ is that it does not lose information, namely function $ f(x) $
could \emph{be recovered} as the linear combination of \intro{coherent
  states} $f_g(x)=[\pi_g f_0](x)$ from its wavelet transform $
\widetilde{f}(g) $:
\begin{equation} \label{eq:g-inverse}
f(x)=\int_G \widetilde{f}(g) f_g(x)\, dg=\int_G \widetilde{f}(g) 
[\pi_g f_0](x) \,dg,
\end{equation}
where $dg$ is the Haar measure on $G$ normalized such that 
\begin{displaymath}
\int_G
\modulus{ \widetilde{f}_0(g)}^2\,dg=1.   
\end{displaymath}
One also has an orthogonal
\emph{projection} $ \widetilde{\oper{P}}$ from $ \FSpace{L}{2}(G,dg)$
to image $ \FSpace{F}{2}(G,dg)$ of $ \FSpace{F}{2}(X) $ under wavelet
transform $\oper{W}$, which is just a convolution on $g$ with the
image $ \widetilde{f}_0(g)=\oper{W}(f_0(x))$ of the vacuum
vector:
\begin{equation} \label{eq:g-proj}
[ \widetilde{ \oper{P}}w](g')=\int_{G} w(g)  \widetilde{f}_0(g^{-1}g')\,dg.
\end{equation}

\subsection{Reduced Wavelets Transform}
Our main observation will be that one could be much more economical
(if subgroup $H$ is non-trivial) with a help
of~\eqref{eq:homogenious}: in this case one need to know $
\widetilde{f}(g) $ not on the whole group $G$ but only on the
homogeneous space $G/H$~\cite[\S~3]{AliAntGazMue}.

Let $ \Omega=G / H$ and $s: \Omega \rightarrow G$ be a continuous
mapping~\cite[\S~13.1]{Kirillov76}. Then any $g\in G$ has a unique
decomposition of the form $g=s(a)h$, $a\in \Omega$ and we will write
$a=s^{-1}(g)$, $h=r(g)={(s^{-1}(g))}^{-1}g$. Note that $ \Omega $ is a
left $G$-homogeneous space\footnote{$ \Omega $ with binary operation $
  (a_1,a_2) \mapsto s^{-1}(s(a_1)\cdot s(a_2))$ becomes a loop of the
  most general form~\cite{Sabinin72}. Thus theory of reduced wavelet
  transform developed in this subsection could be considered as
  \intro{wavelet transform associated with loops}. However we prefer to
  develop our theory based on groups rather on loops.}  with an action
defined in terms of $s$ as follow: $g: a \mapsto s^{-1}(g\cdot s(a))
$. Due to~\eqref{eq:homogenious} one could rewrite~\eqref{eq:wavelets}
as:
\begin{eqnarray*}
\widetilde{f}(g) & = & \scalar{f(x)}{\pi(g)f_0(x)}_{\FSpace{L}{2}( X ) }\\
& = & \scalar{f(x)}{\pi(s(a)h)f_0(x)}_{\FSpace{L}{2}( X ) }\\
& = & \scalar{f(x)}{\pi(s(a))\pi(h)f_0(x)}_{\FSpace{L}{2}( X ) }\\
& = & \scalar{f(x)}{\pi(s(a))\chi(h)f_0(x)}_{\FSpace{L}{2}( X ) }\\
& = & \bar{\chi}(h)\scalar{f(x)}{\pi(s(a))f_0(x)}_{\FSpace{L}{2}( X ) }
\end{eqnarray*}
Thus $\widetilde{f}(g)= \bar{\chi}(h)\widehat{f}(a)$ where
\begin{equation} \label{eq:berg-cauch}
\widehat{f}(a)=[\oper{C} f] (a)=\scalar{f(x)}{\pi(s(a))f_0(x)}_{\FSpace{L}{2}(
X ) }
\end{equation}
and function $\widetilde{f}(g)$ on $G$ is completely defined by
function $ \widehat{f}(a) $ on $ \Omega $.
Formula~\eqref{eq:berg-cauch} gives us an embedding $ \oper{C}:
\FSpace{F}{2}(X) \rightarrow \FSpace{L}{ \infty }( \Omega ) $, which
we will call \intro{reduced wavelet transform}. We denote by $
\FSpace{F}{2}(\Omega)$ the image of $ \oper{C}$ equipped with Hilbert
space inner product induced by $ \oper{C}$ from $ \FSpace{F}{2}(X) $.

Note a special property of $ \widetilde{f}_0(g)$ and $
\widehat{f}_0(a) $:
\begin{displaymath}
\widetilde{f}_0(h^{-1}g)
= \scalar{f_0}{\pi_{h^{-1}g}f_0}
= \scalar{\pi_{h}f_0}{\pi_{g}f_0}
= \scalar{\chi(h)f_0}{\pi_{g}f_0}
=\chi(h) \widetilde{f}_0(g).
\end{displaymath}
It follows from~\eqref{eq:g-inter} that $ \oper{C} $ intertwines $
\rho\oper{C}= \oper{C} \pi$ representation $\pi$ with the
representation
\begin{equation} \label{eq:a-inter}
[\rho_g \widehat{f}](a) = \widehat{f}(s^{-1}(g\cdot s(a))) 
\chi(r(g\cdot s(a))).
\end{equation}
While $\rho$ is not completely geometrical as $\lambda$ in
applications it is still more geometrical than original $\pi$. In many
cases $\rho$ is \intro{representation induced} by the character $\chi$.

If $ f_0(x) $ is a vacuum state with respect to $H$ then
$f_g(x)=\chi(h) f_{s(a)} (x)$ and we could
rewrite~\eqref{eq:g-inverse} as follows:
\begin{eqnarray*}
f(x) & = &  \int_G \widetilde{f}(g) f_g(x)\, dg \nonumber \\
 & = & \int_{ \Omega } \int_H \widetilde{f}(s(a)h) f_{s(a)h}(x)\, 
 dh \, da \nonumber \\
 & = & \int_{ \Omega }  \int_H \widehat{f}(a) \bar{\chi}(h)  \chi(h) 
f_{s(a)} (x)   \, dh\,  da \nonumber \\ 
& = & \int_{ \Omega } \widehat{f}(a) f_{s(a)} (x)\,da \cdot
 \int_H \modulus{ \chi(h) }^2\, dh  \nonumber \\ 
 & = & \int_{ \Omega } \widehat{f}(a)  f_{s(a)}(x) \, da  ,
\end{eqnarray*}
if the Haar measure $dh$ on $H$ is set in such a way that $\int_H
\modulus{ \chi(h) }^2 \, dh=1$ and $dg=dh\,da$. We define an integral
transformation $ \oper{F} $ according to the last formula:
\begin{equation} 
[ \oper{F} \widehat{f}](x) =  \int_{ \Omega } \widehat{f}(a) 
f_{s(a)}(x) \, da \label{eq:a-inverse} ,
\end{equation}
which has the property $ \oper{F} \oper{C} = I$ on $ \FSpace{F}{2}(X)
$ with $ \oper{C} $ defined in~\eqref{eq:berg-cauch}. One could
consider the integral transformation
\begin{equation} \label{eq:szego}
[\oper{P} f](x)=[\oper{F} \oper{C} f](x)=
\int_{ \Omega } \scalar{f(y)}{f_{s(a)}(y)}_{\FSpace{L}{2}(
X ) }  f_{s(a)}(x) \, da
\end{equation}
as defined on whole $ \FSpace{L}{2}(X)$ (not only $ \FSpace{F}{2}(X)
$).  It is known that $ \oper{P} $ \emph{is an orthogonal projection $
  \FSpace{L}{2}(X) \rightarrow \FSpace{F}{2}(X) $}.
If we formally use linearity of the scalar product $
\scalar{\cdot}{\cdot}_{ \FSpace{L}{2}(X)}$ (i.e., assume that the
Fubini theorem holds) we could obtain from~\eqref{eq:szego}
\begin{eqnarray}
[\oper{P} f](x) & = & 
\int_{ \Omega } \scalar{f(y)}{f_{s(a)}(y)}_{\FSpace{L}{2}(
X ) }  f_{s(a)}(x) \, da \nonumber \\
& = &  \scalar{f(y)}{\int_{ \Omega }f_{s(a)}(y) \bar{f}_{s(a)}(x) \, da
}_{\FSpace{L}{2}(
X )}   \nonumber \\
& = & \int_X f(y) K(y,x)\, d\mu(y) \label{eq:bergman},
\end{eqnarray}
where
\begin{displaymath}
K(y,x)=\int_{ \Omega } \bar{ f}_{s(a)}(y)  f_{s(a)}(x) \, da
\end{displaymath}
With the ``probability $ \frac{1}{2}$'' (see discussion on the Bergman
and the Szeg\"o kernels bellow) the integral~\eqref{eq:bergman} exists
in the standard sense, otherwise it is a singular integral operator
(i.e, $K(y,x)$ is a regular function or a distribution).

Sometimes a reduced form $ \widehat{\oper{P}}: \FSpace{L}{2}( \Omega )
\rightarrow \FSpace{F}{2}( \Omega )$ of the projection
$\widetilde{\oper{P}}$~\eqref{eq:g-proj} is of a separate interest. It
is an easy calculation that
\begin{eqnarray} \label{eq:s-b-proj}
[\widehat{\oper{P}} f](a')=\int_{ \Omega } f(a) \widehat{f}_0  ( 
s^{-1}( a^{-1}\cdot a')) \bar{\chi}(r(a^{-1}\cdot a')) \, da,
\end{eqnarray}
where $ a^{-1}\cdot a'$ is an informal abbreviation for $ {\left (s(a)
  \right )}^{-1}\cdot s(a') $. As we will see its explicit form could
be easily calculated in practical cases.

And only at the very end of our consideration we introduce the Taylor
series and the Cauchy-Riemann equations. One knows that they are
\emph{starting points} in the Weierstrass and the Cauchy approaches to
complex analysis correspondingly.

For any decomposition $f_a(x)=\sum_\alpha \psi_\alpha(x) V_\alpha(a)$
of the coherent states $f_a(x)$ by means of functions $V_\alpha(a)$
(where the sum could become eventually an integral) we have the
\intro{Taylor series} expansion
\begin{eqnarray} 
\widehat{f}(a) & = & \int_X f(x) \bar{f}_a(x)\, dx= \int_X f(x) \sum_\alpha 
\bar{\psi}_\alpha(x)\bar{V}_\alpha(a)\, dx  \nonumber \\
 & = &  \sum_\alpha 
\int_X f(x)\bar{\psi}_\alpha(x)\, dx \bar{V}_\alpha(a) \nonumber \\
 & = & \sum_{\alpha}^{\infty} \bar{V}_\alpha(a) f_\alpha,\label{eq:taylor}
\end{eqnarray}
where $f_\alpha=\int_X f(x)\bar{\psi}_\alpha(x)\, dx$.  However to be
useful within the presented scheme such a decomposition should be
connected with structures of $G$ and $H$. For example, if $G$ is a
semisimple Lie group and $H$ its maximal compact subgroup then indices
$\alpha$ run through the set of irreducible unitary representations of
$H$, which enter to the representation $\pi$ of $G$.

The \intro{Cauchy-Riemann equations} need more discussion. One could
observe from~\eqref{eq:g-inter} that the image of $\oper{W}$ is
invariant under action of the left but right regular representations.
Thus $ \FSpace{F}{2}(\Omega)$ is invariant under
representation~\eqref{eq:a-inter}, which is a pullback of the left
regular representation on $G$, but its right counterpart. Thus
generally there is no way to define an action of left-invariant vector
fields on $\Omega$, which are infinitesimal generators of right
translations, on $ \FSpace{L}{2}(\Omega) $. But there is an exception.
Let $ \algebra{X}_j $ be a maximal set of left-invariant vector fields
on $G$ such that
\begin{displaymath} 
\algebra{X}_j \widetilde{f}_0(g)=0.
\end{displaymath}
Because $\algebra{X}_j $ are left invariant we have $\algebra{X}_j
\widetilde{f}_g'(g)=0$ for all $g'$ and thus image of $\oper{W}$,
which the linear span of $\widetilde{f}_g'(g)$, belongs to
intersection of kernels of $\algebra{X}_j$. The same remains true if
we consider pullback $ \widehat{\algebra{X}}_j$ of $\algebra{X}_j$ to
$ \Omega$.  Note that the number of linearly independent
$\widehat{\algebra{X}}_j$ is generally less than for
${\algebra{X}}_j$. We call $ \widehat{\algebra{X}}_j$ as
\intro{Cauchy-Riemann-Dirac operators} in connection with their
property
\begin{equation} \label{eq:dirac1}
\widehat{\algebra{X}}_j \widehat{f}(g)=0 \qquad \forall \widehat{f}(g) 
\in \FSpace{F}{2}(\Omega).
\end{equation}
Explicit constructions of the Dirac type operator for a discrete
series representation could be found in
\cite{AtiyahSchmid80,KnappWallach76}.

We do not use Cauchy-Riemann-Dirac operator in our construction, but
this does not mean that it is useless. One could found at least such
its nice properties:
\begin{enumerate}
\item Being a left-invariant operator it naturally encodes an
  information about symmetry group $G$.
\item It effectively separates irreducible components of the
  representation $\pi$ of $G$ in $\FSpace{L}{2}(X )$.
\item It has a local nature in a neighborhood of a point vs.
  transformations, which act globally on the domain.
\end{enumerate}

\section{Wavelets in Banach Spaces}
\label{sec:wavel-banach-spac}

\subsection{Abstract Nonsence}\label{ss:abstract}
Let $G$ be a group and $H$ be its closed normal subgroup.  Let $X=G/H$
be the corresponding homogeneous space with an invariant measure
$d\mu$ and $s: X \rightarrow G$ be a Borel section in the principal
bundle $G \rightarrow G/H$.  Let $\pi$ be a continuous representation
of a group $G$ by invertible isometry operators $\pi_g$, $g \in G$ in
a (complex) Banach space $B$.

The following definition simulates ones from the Hilbert space
case~\cite[\S~3.1]{AliAntGazMue}.
\begin{defn} \label{de:coherent1}
  Let $G$, $H$, $X=G/H$, $s: X \rightarrow G$, $\pi: G \rightarrow
  \oper{L}(B)$ be as above. We say that $b_0 \in B$ is a \intro{vacuum
    vector} if for all $h\in H$
  \begin{equation} \label{eq:h-char}
    \pi(h) b_0 = \chi(h) b_0, \qquad \chi(h) \in \Space{C}{}.
  \end{equation}
  We will say that set of vectors $b_x=\pi(x) b_0$, $x\in X$ form a
  family of \intro{coherent states} if there exists a continuous non-zero
  linear functional $l_0 \in B^*$ such that
  \begin{enumerate}
  \item \label{it:norm} $\norm{b_0}=1$, $\norm{l_0}=1$,
    $\scalar{b_0}{l_0}\neq 0$;
  \item \label{it:h-char} $\pi(h)^* l_0=\bar{\chi}(h) l_0$, where
    $\pi(h)^*$ is the adjoint operator to $\pi(h)$;
  \item \label{it:coher-eq} The following equality holds
    \begin{equation} \label{eq:coher-eq}
      \int_X \scalar{\pi(x^{-1}) b_0}{l_0}\, \scalar{\pi(x) b_0}{l_0}\, 
      d\mu(x) = \scalar{b_0}{l_0}.
    \end{equation}
  \end{enumerate} 
  The functional $l_0$ is called the \intro{test
  functional}.  According to the strong tradition we call the set
  $(G,H,\pi,B,b_0,l_0)$ \introind{admissible}{wavelets set!admissible}
  if it satisfies to the above conditions.
\end{defn}
We note that mapping $h \rightarrow \chi(h)$ from~\eqref{eq:h-char}
defines a character of the subgroup $H$.  The following Lemma
demonstrates that condition~\eqref{eq:coher-eq} could be relaxed.
\begin{lem}\label{le:exist1}
  For the existence of a vacuum vector $b_0$ and a test functional
  $l_0$ it is sufficient that there exists a vector $b_0'$ and
  continuous linear functional $l'_0$ satisfying to \eqref{eq:h-char}
  and \textup{\ref{it:h-char}} correspondingly such that the constant
\begin{equation} \label{eq:sq-int}
c = \int_X \scalar{\pi(x^{-1})b'_0}{l'_0}\, \scalar{\pi(x) 
b'_0}{l'_0}\, d\mu(x)
\end{equation}
is non-zero and finite.
\end{lem}
\begin{proof}
  There exist a $x_0\in X$ such that $ \scalar{\pi(x_0^{-1})
    b_0'}{l'_0} \neq 0 $, otherwise one has $c=0$.  Let
  $b_0=\pi(x^{-1}) b_0' \norm{\pi(x^{-1}) b_0'}^{-1} $ and $l_0=l_0'
  \norm{l_0'}^{-1}$. For such $b_0$ and $l_0$ we have \ref{it:norm}
  already fulfilled.  To obtain~\eqref{eq:coher-eq} we change the
  measure $d\mu(x)$.  Let $c_0=\scalar{b_0}{l_0} \neq 0 $ then $d\mu'=
  \norm{\pi(x^{-1}) b_0'} \norm{l_0'} c_0 c^{-1}d\mu$ is the desired
  measure.
\end{proof}
\begin{rem}
  Conditions~\eqref{eq:coher-eq} and~\eqref{eq:sq-int} are known for
  unitary representations in Hilbert spaces as \intro{square
    integrability} (with respect to a subgroup $H$).  Thus our
  definition describes an analog of square integrable representations
  for Banach spaces.  Note that in Hilbert space case $b_0$ and $l_0$
  are often the same function, thus condition~\ref{it:h-char} is
  exactly~\eqref{eq:h-char}.  In the particular but still important
  case of trivial $H=\{e\}$ (and thus $X=G$) all our results take
  simpler forms.
\end{rem}
\begin{conv}
  In that follow we will usually write $x\in X$ and $x^{-1}$ instead
  of $s(x)\in G$ and $s(x)^{-1}$ correspondingly.  The right meaning
  of ``$x$'' could be easily found from the context (whether an
  element of $X$ or $G$ is expected there).
\end{conv}

The wavelet transform (similarly to the Hilbert space case) could be
defined as a mapping from $B$ to a space of bounded continuous
functions over $G$ via representational coefficients
\begin{displaymath}
v \mapsto \widehat{v}(g)= \scalar{\pi(g^{-1})v}{l_0}= 
\scalar{v}{\pi (g)^*l_0}.
\end{displaymath}
Due to~\ref{it:h-char} such functions have simple transformation
properties along orbits $gH$, i.e.
$\widehat{v}(gh)=\bar{\chi}(h)\widehat{v}(g)$, $g\in G$, $h\in H$.
Thus they are completely defined by their values indexed by points of
$X=G/H$.  Therefore we prefer to consider so called reduced wavelet
transform.
\begin{defn}
  The \intro{reduced wavelet transform} $\oper{W}$ from a Banach space
  $B$ to a space of function $\FSpace{F}{}(X)$ on a homogeneous space
  $X=G/H$ defined by a representation $\pi$ of $G$ on $B$, a vacuum
  vector $b_0$ and a test functional $l_0$ is given by the formula
\begin{equation} \label{eq:wave-tr}
\oper{W}: B \rightarrow \FSpace{F}{}(X): v \mapsto \widehat{v}(x)= 
[\oper{W}v] (x)=\scalar{\pi(x^{-1}) v}{l_0}=
\scalar{v}{\pi^*(x)l_0}.
\end{equation}
\end{defn}
There is a natural representation of $G$ in $\FSpace{F}{}(X)$.  For
any $g\in G$ there is a unique decomposition of the form $g=s(x)h$,
$h\in H$, $x\in X$.  We will define $r: G \rightarrow H:
r(g)=h=(s^{-1}(g))^{-1}g$ from the previous equality and write a
formal notation $x=s^{-1}(g)$.  Then there is a geometric action of
$G$ on $X \rightarrow X$ defined as follows
\begin{displaymath}
g: x \mapsto g^{-1} \cdot x = s^{-1} (g^{-1} s(x)).
\end{displaymath}
We define a representation $\lambda(g): \FSpace{F}{}(X) \rightarrow
\FSpace{F}{}(X)$ as follow
\begin{equation} \label{eq:l-rep}
[\lambda(g) f] (x) = \chi(r(g^{-1}\cdot x)) f(g^{-1}\cdot x).
\end{equation}
We recall that $\chi(h)$ is a character of $H$ defined
in~\eqref{eq:h-char} by the vacuum vector $b_0$. For the case of
trivial $H=\{e\}$ \eqref{eq:l-rep} becomes the left regular
representation $\rho_l(g)$ of $G$.
\begin{prop} \label{pr:inter1}
  The reduced wavelet transform $\oper{W}$ intertwines $\pi$ and
  the representation $\lambda$~\eqref{eq:l-rep} on $\FSpace{F}{}(X)$:
\begin{eqnordisp}
  \oper{W} \pi(g) = \lambda(g) \oper{W}.
\end{eqnordisp}
\end{prop}
\begin{proof}
  We have:
\begin{eqnarray*}{}
[\oper{W}( \pi(g) v)] (x) & = & \scalar{\pi(x^{-1}) \pi(g) v }{ l_0} \\
  & = & \scalar{\pi((g^{-1}s(x))^{-1}) v }{ l_0} \\
  & = & \scalar{\pi(r(g^{-1}\cdot x)^{-1})\pi(s(g^{-1}\cdot x)^{-1}) v }{ l_0} \\
  & = & \scalar{\pi(s(g^{-1}\cdot x)^{-1}) v }{\pi^*(r(g^{-1}\cdot x)^{-1}) 
l_0} \\
  & = & \chi(r(g^{-1}\cdot x)^{-1}) [\oper{W} v] (g^{-1}x) \\
  & = & \lambda(g) [\oper{W}v] (x).
\end{eqnarray*}
\end{proof}
\begin{cor}\label{co:pi}
  The function space $\FSpace{F}{}(X)$ is invariant under the
  representation $\lambda$ of $G$.
\end{cor}
We will see that $\FSpace{F}{}(X)$ posses many properties of the
\intro{Hardy space}. 
The duality between $l_0$ and $b_0$ generates a transform dual to $\oper{W}$.
\begin{defn}
  The \intro{inverse wavelet transform} $\oper{M}$ from $
  \FSpace{F}{}(X) $ to $B$ is given by the formula:
\begin{eqnarray}
\oper{M}:  \FSpace{F}{}(X) \rightarrow B: \widehat{v}(x) \mapsto \oper{M} 
[\widehat{v}(x)] & = & \int_X \widehat{v}(x) b_x\,d\mu(x) \nonumber\\
 & = & \int_X \widehat{v}(x) \pi(x)\,d\mu(x) b_0. \label{eq:m-tr}
\end{eqnarray}
\end{defn}
\begin{prop} \label{pr:inter2}
  The inverse wavelet transform $ \oper{M} $ intertwines the
  representation $ \lambda $ on $ \FSpace{F}{}(X) $ and $ \pi $ on
  $B$:
\begin{eqnordisp}
  \oper{M} \lambda(g) = \pi(g) \oper{M}.
\end{eqnordisp}
\end{prop}
\begin{proof}
  We have:
\begin{eqnarray*}
\oper{M} [\lambda(g)\widehat{v}(x)] 
& = & \oper{M} [ \chi(r(g^{-1}\cdot x)) \widehat{v}(g^{-1}\cdot x)] \\
& = & \int_X \chi(r(g^{-1}\cdot x)) \widehat{v}(g^{-1}\cdot x) 
        b_x \,d\mu(x)\\
& = & \chi(r(g^{-1}\cdot x)) \int_X \widehat{v}(x') b_{g\cdot x'}\,d\mu(x')
\\
& = & \pi_g \int_X \widehat{v}(x') b_{x'}\,d\mu(x')\\
& = & \pi_g \oper{M} [\widehat{v}(x')],
\end{eqnarray*}  
where $x'=g^{-1} \cdot x$.
\end{proof}
\begin{cor}\label{co:lambda}
  The image $\oper{M}(\FSpace{F}{}(X))\subset B$ of subspace
  $\FSpace{F}{}(X)$ under the inverse wavelet transform $\oper{M}$ is
  invariant under the representation $\pi$.
\end{cor}
The following proposition explain the usage of the name for
$\oper{M}$.
\begin{thm}
  The operator
\begin{equation} \label{eq:szego1}
\oper{P}= \oper{M} \oper{W}: B \rightarrow B
\end{equation}
is a projection of $B$ to its linear subspace for which $b_0$ is
cyclic. Particularly if $\pi$ is an irreducible representation then the
inverse wavelet transform $\oper{M}$ is a \introind{left inverse
  operator}{operator!left inverse}
on $B$ for the wavelet transform $\oper{W}$:
\begin{eqnordisp}
  \oper{M}\oper{W}=I.
\end{eqnordisp}
\end{thm}
\begin{proof}
  It follows from Propositions~\ref{pr:inter1} and~\ref{pr:inter2}
  that operator $\oper{M}\oper{W}: B \rightarrow B$ intertwines $\pi$
  with itself.  Then Corollaries~\ref{co:pi} and~\ref{co:lambda} imply
  that the image $\oper{M}\oper{W}$ is a $\pi$-invariant subspace of
  $B$ containing $b_0$.  Because $\oper{M}\oper{W}b_0=b_0$ we conclude
  that $\oper{M}\oper{W}$ is a projection.
  
  From irreducibility of $\pi$ by Schur's
  Lemma~\cite[\S~8.2]{Kirillov76} one concludes that
  $\oper{M}\oper{W}=cI$ on $B$ for a constant $c\in\Space{C}{}$.
  Particularly
\begin{displaymath}
\oper{M}\oper{W} b_0= \int_X \scalar{\pi(x^{-1})b_0}{l_0}\, \pi(x) 
b_0\,d\mu(x)=cb_0.
\end{displaymath}
From the condition~\eqref{eq:coher-eq} it follows that
$\scalar{cb_0}{l_0}=\scalar{\oper{M}\oper{W}
  b_0}{l_0}=\scalar{b_0}{l_0}$ and therefore $c=1$.
\end{proof}
We have similar
\begin{thm}
  Operator $\oper{W}\oper{M}$ is a projection of $\FSpace{L}{1}(X)$ to
  $\FSpace{F}{}(X)$. 
\end{thm}

We denote by $\oper{W}^*: \FSpace{F}{}^*(X) \rightarrow B^* $ and
$\oper{M}^*: B^* \rightarrow \FSpace{F}{}^*(X)$ the adjoint (in the
standard sense) operators to $\oper{W}$ and $\oper{M}$ respectively.
\begin{cor}
  We have the following identity:
\begin{equation} \label{eq:isom1}
\scalar{\oper{W} v }{ \oper{M}^* l}_{ \FSpace{F}{}(X) } = \scalar{v}{l}_B, 
\qquad \forall v\in B, \quad l\in B^*
\end{equation}
or equivalently
\begin{equation} \label{eq:isom2}
\int_X \scalar{\pi(x^{-1}) v}{l_0}\, \scalar{\pi(x) b_0}{l}\, d\mu(x) 
= \scalar{v}{l}.
\end{equation}
\end{cor}
\begin{proof}
  We show the equality in the first form~\eqref{eq:isom2} (but will
  apply it often in the second one):
\begin{displaymath}
\scalar{\oper{W} v }{ \oper{M}^* l}_{ \FSpace{F}{}(X) }
 = \scalar{\oper{M}\oper{W} v }{l}_B =\scalar{v}{l}_B.
\end{displaymath}
\end{proof}
\begin{cor}
  The space $ \FSpace{F}{}(X) $ has the reproducing formula
\begin{equation} \label{eq:reprod}
\widehat{v}(y)=\int_X \widehat{v}(x) \, 
   \widehat{b}_0(x^{-1}\cdot y)\,d\mu(x),
\end{equation}
where $ \widehat{b}_0(y)=[\oper{W}b_0] (y)$ is the wavelet transform
of the vacuum vector $b_0$.
\end{cor}
\begin{proof}
  Again we have a simple application of the previous formulas:
\begin{eqnarray}
\widehat{v}(y)& =& \scalar{\pi(y^{-1})v}{l_0} \nonumber \\ 
&= & \int_X \scalar{\pi(x^{-1}) \pi(y^{-1}) v}{l_0}\, \scalar{\pi(x) 
b_0}{l_0}\,d\mu(x) \label{eq:rep-tr1} \\
&= & \int_X \scalar{\pi(s(y\cdot x)^{-1}) v}{l_0}\, \scalar{\pi(x) 
b_0}{l_0}\,d\mu(x) \nonumber \\
&= & \int_X \widehat{v} (y\cdot x)\, \widehat{b}_0(x^{-1})  \,d\mu(x)
\nonumber \\
& = & \int_X \widehat{v}(x)\, \widehat{b}_0(x^{-1}y)\,d\mu(x), \nonumber
\end{eqnarray}
where transformation~\eqref{eq:rep-tr1} is due to~\eqref{eq:isom2}.
\end{proof}
\begin{rem}
  To possess a reproducing kernel---is a well-known property of spaces
  of analytic functions.  The space $\FSpace{F}{}(X)$ shares also
  another important property of analytic functions: it belongs to a
  kernel of a certain first order differential operator with Clifford
  coefficients (the Dirac operator) and a second order operator with
  scalar coefficients (the Laplace
  operator)~\cite{AtiyahSchmid80,Kisil97a,Kisil97c,KnappWallach76}.
\end{rem}
Let us now assume that there are two representations $\pi'$ and
$\pi''$ of the same group $G$ in two different spaces $B'$ and $B''$
such that two admissible sets $(G,H,\pi',B',b_0',l_0')$ and
$(G,H,\pi'',B'',b_0'',l_0'')$ could be constructed for the same normal
subgroup $H\subset G$.
\begin{prop} \label{pr:2repr}
  In the above situation if $F'(X) \subset F''(X)$ then the
  composition $\oper{T}=\oper{M}''\oper{W}'$ of the wavelet transform
  $\oper{W}'$ for $\pi'$ and the inverse wavelet transform
  $\oper{M}''$ for $\pi''$ is an intertwining operator between $\pi'$
  and $\pi''$:
\begin{eqnordisp}[]
  \oper{T}\pi'=\pi''\oper{T}.
\end{eqnordisp}
$\oper{T}$ is defined as follows
\begin{eqnordisp}[eq:T-intertw]
  \oper{T}: b \mapsto \int_X \scalar{\pi'(x^{-1})b}{l'_0}\,
  \pi''(x)b_0''\,d\mu(x).
\end{eqnordisp}
This transformation defines a $B''$-valued linear functional (a
distribution for function spaces) on $B'$.
\end{prop}
The Proposition has an obvious proof.  This simple result is a base
for an alternative approach to functional calculus of
operators~\cite{Kisil95i,Kisil97a,Kisil02c}.  Note also that
formulas~\eqref{eq:wave-tr} and \eqref{eq:m-tr} are particular cases
of~\eqref{eq:T-intertw} because $\oper{W}$ and $\oper{M}$ intertwine
$\pi$ and $\lambda$.

\subsection{Singular Vacuum Vectors} \label{ss:singular}
In many important cases the above general scheme could not be carried
out because the representation $\pi$ of $G$ is not square-integrable
or even not square-integrable modulo a subgroup $H$.  Thereafter the
vacuum vector $b_0$ could not be selected within the original space
$B$ which the representation $\pi$ acts on.  The simplest mathematical
example is the Fourier transform (see\cite{Kisil98a}). In
physics this is the well-known problem of \intro{absence of vacuum
  state} in the constructive algebraic \introind{quantum field
  theory}{field theory!quantum}
\cite{Segal90,Segal94,Segal96a}.  The absence of the vacuum within the
linear space of system's states is another illustration to the old
thesis \emph{Natura abhorret vacuum}\footnote{Nature is horrified by
  (any) vacuum (Lat.).} or even more specifically \emph{Natura
  abhorret vectorem vacui}\footnote{Nature is horrified by a carrier
  of nothingness (Lat.).  This illustrates how far a humane beings
  deviated from Nature.}.

We will present a modification of our construction which works in such
a situation.  For a singular vacuum vector the algebraic structure of
group representations could not describe the situation alone and
requires an essential assistance from analytical structures.
\begin{defn} \label{de:coherent2}
  Let $G$, $H$, $X=G/H$, $s: X \rightarrow G$, $\pi: G \rightarrow
  \oper{L}(B)$ be as in Definition~\ref{de:coherent1}. We assume that
  there exist a topological linear space $\widehat{B}\supset B$ such
  that
\begin{enumerate}
\item \label{it:bexist} $B$ is dense in $\widehat{B}$ (in topology of
  $\widehat{B}$) and representation $\pi$ could be uniquely extended
  to the continuous representation $\widehat{\pi}$ on $\widehat{B}$.
\item There exists $b_0 \in \widehat{B}$ be such that for all $h\in H$
\begin{equation} \label{eq:h-char2}
\widehat{\pi}(h) b_0 = \chi(h) b_0, \qquad \chi(h) \in \Space{C}{}.
\end{equation}
\item There exists a continuous non-zero linear functional $l_0 \in
  B^*$ such that \label{it:h-char2a} $\pi(h)^* l_0=\bar{\chi}(h) l_0$,
  where $\pi(h)^*$ is the adjoint operator to $\pi(h)$;
\item \label{it:b2b} The composition $\oper{M}\oper{W}: B \rightarrow
  \widehat{B}$ of the wavelet transform \eqref{eq:wave-tr} and the
  inverse wavelet transform \eqref{eq:m-tr} maps $B$ to $B$.
\item \label{it:coher-eq2} For a vector $p_0\in B$ the following
  equality holds
\begin{equation} \label{eq:coher-eq2}
\scalar{\int_X \scalar{\pi(x^{-1}) p_0}{l_0}\, \pi(x) b_0\, 
d\mu(x) }{l_0}= \scalar{p_0}{l_0},
\end{equation}
where the integral converges in the weak topology of $\widehat{B}$.
\end{enumerate}
As before we call the set of vectors $b_x=\pi(x) b_0$, $x\in X$ by
\intro{coherent states}; the vector $b_0$---a \intro{vacuum vector}; the
functional $l_0$ is called the \intro{test functional} and finally
$p_0$ is the \intro{probe vector}.
\end{defn}
This Definition is more complicated than
Definition~\ref{de:coherent1}.  The equation~\eqref{eq:coher-eq2} is a
substitution for~\eqref{eq:coher-eq} if the linear functional $l_0$ is
not continuous in the topology of $\widehat{B}$.
 The function theory in \Space{R}{1,1} constructed in
the next lecture provides a more exotic example of a singular vacuum
vector.

We shall show that \ref{it:coher-eq2} could be satisfied by an
adjustment of other components.
\begin{lem}\label{le:exist2}
  For the existence of a vacuum vector $b_0$, a test functional $l_0$,
  and a probe vector $p_0$ it is sufficient that there exists a vector
  $b_0'$ and continuous linear functional $l'_0$ satisfying to
  \textup{\ref{it:bexist}--\ref{it:b2b}} and a vector $p'_0\in B$ such
  that the constant
\begin{eqnordisp}[]
  c=\scalar{\int_X \scalar{\pi(x^{-1}) p_0}{l_0}\, \pi(x) b_0\,
    d\mu(x) }{l_0}
\end{eqnordisp}
is non-zero and finite.
\end{lem}
The proof follows the path for Lemma~\ref{le:exist1}.  The following
Proposition summarizes results which could be obtained in this case.
\begin{prop} Let the wavelet transform $\oper{W}$ \eqref{eq:wave-tr},  
  its inverse $\oper{M}$ \eqref{eq:m-tr}, the representation
  $\lambda(g)$~\eqref{eq:l-rep}, and functional space
  $\FSpace{F}{}(X)$ be adjusted accordingly to
  Definition~\textup{\ref{de:coherent2}}. Then
\begin{enumerate}
\item $\oper{W}$ intertwines $\pi(g)$ and $\lambda(g)$ and the image
  of $\FSpace{F}{}(X)=\oper{W}(B)$ is invariant under $\lambda(g)$.
\item $\oper{M}$ intertwines $\lambda(g)$ and $\widehat{\pi}(g)$ and
  the image of $\oper{M}(\FSpace{F}{}(B))=\oper{M}\oper{W}(B) \subset
  B$ is invariant under $\pi(g)$.
\item If $\oper{M}(\FSpace{F}{}(X))=B$ (particularly if $\pi(g)$ is
  irreducible) then $\oper{M}\oper{W}=I$ otherwise $\oper{M}\oper{W}$
  is a projection $B \rightarrow \oper{M}(\FSpace{F}{}(X))$. In both
  cases $\oper{M}\oper{W}$ is an operator defined by integral
\begin{eqnordisp}[eq:mw2]
  b \mapsto \int_X \scalar{\pi(x^{-1}) b}{l_0}\, \pi(x) b_0\, d\mu(x),
\end{eqnordisp}
\item Space $\FSpace{F}{}(X)$ has a reproducing formula
\begin{eqnordisp}[eq:reprod2]
  \widehat{v}(y)=\scalar{\int_X \widehat{v}(x) \, \pi(x^{-1}y)b_0\,
    dx}{l_0}
\end{eqnordisp}
which could be rewritten as a singular convolution
\begin{eqnordisp}[]
  \widehat{v}(y)=\int_X \widehat{v}(x) \, \widehat{b}(x^{-1}y)\, dx
\end{eqnordisp}
with a distribution $b(y)=\scalar{\pi(y^{-1})b_0}{l_0}$ defined by
\eqref{eq:reprod2}.
\end{enumerate}
\end{prop}
The proof is algebraic and completely similar to
Subsection~\ref{ss:abstract}.

\chapter[Hyperbolic Function Theory]{Analytical 
Function Theory of Hyperbolic Type}

\section{Introduction}
\epigraph{You should complete your own original research in order to 
learn when it was done before.}{}{}
Connections between complex analysis (one variable, several complex 
variables, Clifford analysis) and its symmetry groups are known from its 
earliest days. They are an obligatory part of the textbook on the 
subject~\cite{Cnops94a}, \cite{DelSomSou92}, \cite[\S~1.4, 
\S~5.4]{GilbertMurray91}, \cite{Krantz82}, \cite[Chap.~2]{Rudin80} and 
play an essential role in many research 
papers~\cite{PeetreQian94,Ryan95d,Ryan98a} just to mention only few.
However ideas about fundamental role of symmetries in function
theories outlined in~\cite{GilbertKunzeTomas86,Koranyi72a} were not
incorporated in a working toolkit of researchers yet.

It was proposed in the first lecture to distinguish essentially different
function theories by corresponding group of symmetries.  Such a
classification is needed because not all seemingly different function
theories are essentially different, see the first lecture 
and~\cite{Kisil95c}.  But it is also
important that the group approach gives a constructive way to develop
essentially different function theories (see the last two
lectures~\cite{Kisil96c,Kisil97a,Kisil02c}), as well 
as outlines an alternative ground for functional calculi of
operators~\cite{Kisil95i}.  In the mentioned papers all given examples
consider only well-known function theories.  While rearranging of known
results is not completely useless there was an appeal to produce a new
function theory based on the described scheme.

The theorem proved in~\cite{Kisil96d} underlines the similarity between
structure of the group of M\"obius transformations in spaces $ \Space{R}{n}
$ and $ \Space{R}{pq} $. This generates a hope that there exists a
non empty function theory in $ \Space{R}{pq} $. We construct such a theory
in the present lecture for the case of $ \Space{R}{1,1} $. Other new function
theories based on the same scheme will be described
in the next lecture~\cite{CnopsKisil97a}.

The format of the lecture is as follows. In Section~\ref{se:preliminaries} we 
introduce basic notations and definitions. We construct two function 
theories---the standard complex analysis and a function theory in 
$\Space{R}{1,1}$---in Section~\ref{se:theories}. Our consideration is based
on two different series of representation of $\SL$: discrete and principal.
We deduce in their terms the Cauchy integral formula, the Hardy spaces, the
Cauchy-Riemann equation, the Taylor expansion and their counterparts for 
$\Space{R}{1,1}$. Finally we collect in
Appendices~\ref{sec:groups-homog-spac} and~\ref{sec:elem-repr-theory} several 
facts, which we would like (however can not) to assume well known. It may 
be a good idea to look through the 
Appendixes~~\ref{sec:miscellanea} between the reading of
Sections~\ref{se:preliminaries} and~\ref{se:theories}. 
Finally Appendix~\ref{pt:o-problems} states few among many open problems.
Our examples will be rather lengthy thus their (not always obvious) ends 
will be indicated by the symbol $\eoe$.

\section{Preliminaries} \label{se:preliminaries}
Let $\Space{R}{pq}$ be a real $n$-dimensional vector space, where
$n=p+q$ with a fixed frame $e_1$, $e_2$, \ldots, $e_p$, $e_{p+1}$,
\ldots, $e_n$ and with the nondegenerate bilinear
form $B(\cdot,\cdot)$ of the signature $(p,q)$, which is diagonal in the
frame $e_i$, i.e.:
\begin{displaymath}
B(e_i,e_j)=\epsilon_i \delta_{ij}, \textrm{ where }
\epsilon_i=\left\{\begin{array}{ll}
1, & i=1,\ldots,p\\
-1, & i=p+1,\ldots,n
\end{array}\right.
\end{displaymath}
and $\delta _{ij}$ is the Kronecker delta. In particular the usual
Euclidean space $\Space{R}{n}$ is $\Space{R}{0n}$. Let $\Cliff[p]{q}$ be
the \intro{real Clifford algebra} generated by $1$, $e_j$, $1\leq j\leq n$
and the relations 
\begin{displaymath}
e_i e_j + e_j e_i =-2B(e_i,e_j).
\end{displaymath}
We put $e_0=1$ also. Then there is the natural embedding
$\algebra{i}: \Space{R}{pq}\rightarrow
\Cliff[p]{q}$. We identify $\Space{R}{pq}$ with its image under
$\algebra{i}$ and
call its elements \emph{vectors}. There are two linear
anti{-}automorphisms $*$ (reversion) and $-$
(main anti{-}automorphisms) 
and automorphism $'$
of $\Cliff[p]{q}$ defined on its
basis $A_\nu=e_{j_1}e_{j_2}\cdots e_{j_r}$, $1\leq j_1 <\cdots<j_r\leq
n$ by the rule:
\begin{eqnarray*}
(A_\nu)^*= (-1)^{\frac{r(r-1)}{2}} A_\nu, \qquad
\bar{A}_\nu= (-1)^{\frac{r(r+1)}{2}} A_\nu,\qquad
A_\nu'= (-1)^{r} A_\nu.
\end{eqnarray*}
In particular, for vectors, $\bar{\vecbf{x}}=\vecbf{x}'=-\vecbf{x}$ and
$\vecbf{x}^*=\vecbf{x}$.

It is easy to see that $\vecbf{x}\vecbf{y}=\vecbf{y}\vecbf{x}=1$
for any $\vecbf{x}\in\Space{R}{pq}$ such that
$B(\vecbf{x},\vecbf{x})\neq 0$
and $\vecbf{y}={\bar{\vecbf{x}}}\,{\norm{\vecbf{x}}^{-2}}$, which is
the \intro{Kelvin inverse} of $\vecbf{x}$.
Finite products of invertible vectors are invertible in $\Cliff[p]{q}$
and form the \intro{Clifford group} $\Gamma(p,q)$. Elements
$a\in\Gamma(p,q)$ such that
$a\bar{a}=\pm 1$ form the $\object[(p,q)]{Pin}$ group---the double cover
of the group of orthogonal rotations $\object[(p,q)]{O}$. We also
consider~\cite[\S~5.2]{Cnops94a} $T(p,q)$ to be the set of all
products of vectors in $\Space{R}{pq} $.

Let $(a, b, c, d)$ be a quadruple from $T(p,q)$ with
the properties:
\begin{enumerate}
\item $(ad^*-bc^*)\in \Space{R}{}\setminus {0}$;
\item $a^*b$, $c^*d$, $ac^* $, $ bd^*$ are vectors.
\end{enumerate}
Then~\cite[Theorem~5.2.3]{Cnops94a}
$2\times 2$-matrixes $\matr{a}{b}{c}{d}$ form the
group $\Gamma(p+1,q+1)$ under the
usual matrix multiplication. It has a representation
$\pi_{\Space{R}{pq} }$ \comment{we denote its restriction to any
subgroup by the same notation.} by transformations of
$\overline{\Space{R}{pq} }$ given by:
\begin{equation}\label{eq:sp-rep}
\pi_{\Space{R}{pq}}\matr{a}{b}{c}{d} :
\vecbf{x} \mapsto (a\vecbf{x}+b)(c\vecbf{x}+d)^{-1},
\end{equation}
which form the \intro{M\"obius} (or
the \intro{conformal}) group of $\overline{\Space{R}{pq}}$. 
Here $\overline{\Space{R}{pq}}$ the compactification of $\Space{R}{pq} $
by the ``necessary number 
of points'' (which form the light cone) at infinity 
(see~\cite[\S~5.1]{Cnops94a}).
The analogy with fractional-linear transformations of the complex line
\Space{C}{} is useful, as well as representations of shifts
$\vecbf{x}\mapsto \vecbf{x}+y$, orthogonal rotations
$\vecbf{x}\mapsto k(a)\vecbf{x}$, dilations
$\vecbf{x}\mapsto \lambda \vecbf{x}$, and the Kelvin inverse
$\vecbf{x}\mapsto \vecbf{x}^{-1}$ by the
matrixes $\matr{1}{y}{0}{1}$, 
$\matr{a}{0}{0}{{a}^{*-1}}$,
\matr{\lambda^{1/2}}{0}{0}{\lambda^{-1/2}}, \matr{0}{-1}{1}{0}
respectively. We also use the agreement of~\cite{Cnops94a} that $ 
\frac{a}{b} $ always denotes $ab^{-1}$ for $a$, $b\in \Cliff[p]{q}$.

\section[Two Function Theories from $\SL$]{Two Function 
Theories Associated with Representations of $\SL$} 
\label{se:theories}
\subsection{Unit Disks in $ \Space{R}{0,2} $ and $ \Space{R}{1,1} $}
The main example is provided by group 
$G=\SL$ (books~\cite{HoweTan92,Lang85,MTaylor86} are our
standard references about $ \SL $ and its
representations) consisting of $ 2\times 2$ matrices $
\matr{a}{b}{c}{d}$ with real entries and determinant $ad-bc=1$. 

The Lie algebra $ \algebra{sl}(2,\Space{R}) $ of $\SL$ consists of all 
$2\times 2$ real matrices of trace zero. One can introduce a basis
\begin{displaymath}
A= \frac{1}{2} \matr{-1}{0}{0}{1},\quad B= \frac{1}{2} \matr{0}{1}{1}{0}, 
\quad Z=\matr{0}{1}{-1}{0}.
\end{displaymath}
The commutator relations are
\begin{displaymath}
[Z,A]=2B, \qquad [Z,B]=-2A, \qquad [A,B]=- \frac{1}{2} Z.
\end{displaymath}

We will construct two series of examples. One is connected with discreet 
series representation and produces the core of standard complex analysis. 
The second will be its mirror in principal series representations and 
create parallel function theory. $ \SL $ has also other type 
representation, which can be of particular interest in other 
circumstances. However the discreet series and principal ones stay 
separately from others (in particular by being the support of the
Plancherel measure~\cite[\S~VIII.4]{Lang85},~\cite[Chap.~8,
(4.16)]{MTaylor86}) and are in a good resemblance each other.

\begin{examplea}
Via identities
\begin{displaymath}
\alpha= \frac{1}{2}(a+d-ic+ib), \qquad \beta= \frac{1}{2}(c+b -ia+id)
\end{displaymath}
we have isomorphism of $ \SL $ with group $SU(1,1)$ of $
2\times 2$ matrices with complex entries of the form $ \matr{ \alpha}{
\beta}{ \bar{\beta}}{\bar{ \alpha}}$ such that $ \modulus{ \alpha }^2 - 
\modulus{ \beta }^2=1 $. We will use the last form for  $ \SL$ for complex 
analysis in unit disk $\Space{D}{}$.

$ \SL $ has the only non-trivial compact closed subgroup
$K$, namely the group of matrices of the form
$h_{\psi}=\matr{e^{i\psi}}{0}{0}{e^{-i\psi}}$.
Now any $ g\in \SL $ has a unique decomposition of the
form 
\begin{eqnarray}
\matr{\alpha}{\beta}{\bar{\beta}}{\bar{\alpha}} 
& =& \modulus{\alpha} 
\matr{1}{\beta\bar{\alpha}^{-1}}{\bar{\beta}\alpha^{-1}}{1} \matr{ 
\frac{{\alpha}}{ \modulus{\alpha} } }{0}{0}{\frac{\bar{\alpha}}{ 
\modulus{\alpha} } }
\nonumber \\
& =& \frac{1}{\sqrt{1- \modulus{a}^2 }}
\matr{1}{a}{\bar{a}}{1} 
\matr{e^{i\psi}}{0}{0}{e^{-i\psi}} 
\end{eqnarray}
where $\psi=\Im \ln \alpha$, $a=\beta\bar{\alpha}^{-1}$, and $
\modulus{ a } < 1$ because $ \modulus{ \alpha }^2 - \modulus{ \beta
}^2=1 $. Thus we can identify $\SL / H$ with the unit
disk $ \Space{D}{} $ and define mapping $s: \Space{D}{} \rightarrow
\SL $ as follows
\begin{equation} \label{eq:def-s-a}
s: a \mapsto \frac{1}{\sqrt{1- \modulus{a}^2 }}
\matr{1}{a}{\bar{a}}{1}.
\end{equation}
Mapping $ r: G \rightarrow H$ associated to $s$  is
\begin{equation} \label{eq:def-r-a}
r: \matr{\alpha}{\beta}{\bar{\beta}}{\bar{\alpha}} \mapsto \matr{ 
\frac{{\alpha}}{
\modulus{\alpha} }}{0}{0}{\frac{\bar{\alpha}}{ \modulus{\alpha} }}
\end{equation}

The invariant measure $ d\mu(a) $ on $ \Space{D}{} $ coming from
decomposition $dg=d\mu(a)\, dk$, where $dg$ and $dk$ are Haar measures
on $G$ and $K$ respectively, is equal to 
\begin{equation} \label{eq:def-m-a}
d\mu(a)= \frac{da}{{(1- \modulus{a}^2)}^2}
\end{equation}
with $da$---the standard Lebesgue measure on $ \Space{D}{} $.

The formula $g: a \mapsto g\cdot a=s^{-1}(g^{-1} * s(a)) $ associates with
a matrix $ g^{-1}=\matr{ \alpha}{ \beta}{ \bar{\beta}}{ \bar{\alpha}}$ the
fraction-linear transformation of $ \Space{D}{} $ of the form
\begin{equation} \label{eq:fr-lin-a}
g: z \mapsto g\cdot z=
\frac{\alpha z + \beta}{\bar{ \beta}z +\bar{\alpha}},  
\qquad g^{-1}=\matr{\alpha}{\beta}{\bar{\beta}}{\bar{\alpha}},
\end{equation}
which also can be considered as a transformation of $\dot{ \Space{C}{} }$
(the one-point compactification of $ \Space{C}{} $). $\eoe$
\end{examplea}

\begin{exampleb}
We will describe a version of previous formulas corresponding to 
geometry of unit disk in $ \Space{R}{1,1}$. For generators $e_1$ and $e_2$ 
of $ \Space{R}{1,1} $ (here $e_1^2=-e^2_2=-1$) we see that matrices 
$\matr{a}{be_2}{ce_2}{d}$ again give a realization of $\SL$. Making 
composition with the Caley transform
\begin{displaymath}
T=\frac{1}{ {2} } \matr{1}{e_2}{e_2}{-1} \matr{1}{e_1}{e_1}{1} =
\frac{1}{ {2} } \matr{1+e_2e_1}{e_1+e_2}{e_2-e_1}{e_2e_1-1}
\end{displaymath} 
and its inverse 
\begin{displaymath}
T^{-1}=\frac{1}{ {2} } \matr{1}{-e_1}{-e_1}{1} \matr{1}{e_2}{e_2}{-1}=
\frac{1}{2} \matr{1-e_1e_2}{e_2+e_1}{e_2-e_1}{-1-e_1e_2}
\end{displaymath}
(see analogous calculation in~\cite[\S~IX.1]{Lang85}) we obtain another
realization of $\SL$:
\begin{equation} \label{eq:caley}
\frac{1}{4} \matr{1-e_1e_2}{e_2+e_1}{e_2-e_1}{-1-e_1e_2}
\matr{a}{be_2}{ce_2}{d} 
\matr{1+e_2e_1}{e_1+e_2}{e_2-e_1}{e_2e_1-1}
=\matr{\n{a}}{\n{b}}{\n{b}'}{\n{a}'},
\end{equation}                     
where                              
\begin{equation}                   
\n{a}= \frac{1}{2}(a(1-e_1e_2)+d(1+e_1e_2)), \qquad 
\n{b}= \frac{1}{2} (b(e_1-e_2)+c(e_1+e_2)).
\end{equation}
It is easy to check that the condition $ad-bc=1$ implies the following
value of the pseudodeterminant of the matrix $ \n{a}(\n{a}')^*
-\n{b}(\n{b}')^*= \n{a}\bar{\n{a}}-\n{b} \bar{\n{b}}=1$.  We also observe
that $\n{a}$ is an even Clifford number and $\n{b}$ is a vector thus
$\n{a}'=\n{a}$, $\n{b}'=-\n{b}$.

Now we consider the decomposition
\begin{equation}
\matr{\n{a}}{\n{b}}{-\n{b}}{\n{a}} = 
\modulus{\n{a}} \matr{1}{\n{b}\n{a}^{-1}}{-\n{b}\n{a}^{-1}}{1}   \matr{ 
\frac{\n{a}}{ \modulus{\n{a}}} }{0}{0}{\frac{\n{a}}{ \modulus{\n{a}}}}.
\end{equation}
It is seen directly, or alternatively follows from general characterization 
of $\Gamma(p+1,q+1)$~\cite[Theorem~5.2.3(b)]{Cnops94a}, that
$\n{b}\n{a}^{-1}\in \Space{R}{1,1}$.  Note that now we cannot derive from $
\n{a}\bar{\n{a}}-\n{b} \bar{\n{b}}=1$ that $\n{b}\n{a}^{-1}
\overline{\n{b}\n{a}^{-1}}= -{(\n{b}\n{a}^{-1})}^2 <1 $ because $
\n{a}\bar{\n{a}}$ can be positive or negative (but we are sure that
${(\n{b}\n{a}^{-1})}^2\neq -1$).  For this reason we cannot define the unit
disk in $ \Space{R}{1,1} $ by the condition $ \modulus{\n{u}}<1 $ in a way
consistent with its M\"obius transformations.  This topic will be discussed
elsewhere with more illustrations~\cite{CnopsKisil96a}.

\begin{figure}[tbh]
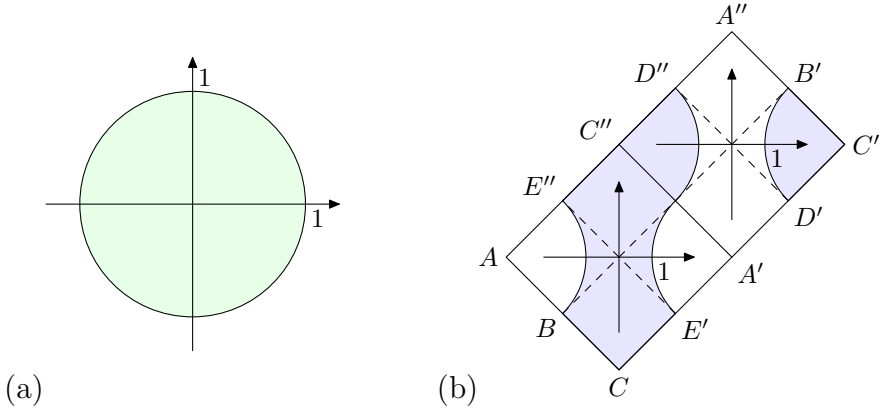

  \begin{center}
    \href{http://maths.leeds.ac.uk/~kisilv/r11.gif}{
      (a)\includegraphics{nccg2.1} \hspace{1cm}
      (b)\includegraphics{nccg2.2}}
    \caption[Two unit disks]
    {
      Two unit disks in elliptic (a) and hyperbolic (b)
      metrics. In (b) squares $ACA'C''$ and $A'C'A''C''$ represent
      two copies of $\Space{R}{2}$, their boundaries are the image
      of the light cone at infinity. These cones should be glued in
      a way to merge points with the same letters (regardless number
      of dashes). The hyperbolic unit disk $\TSpace{D}{}$ (shaded
      area) is bounded by four branches of hyperbola. Dashed lines
      are light cones at origins.}
    \label{pi:disk}
  \end{center}
\end{figure}

We are taking two copies $\mathbb{R}^{1,1}_+ $ and $ \mathbb{R}^{1,1}_- $
of $ \Space{R}{1,1} $ glued over their light cones at infinity in such a
way that the construction is invariant under natural action of the M\"obius
transformation.  This aggregate denoted by $ \TSpace{R}{1,1}$ is a two-fold
cover of $ \Space{R}{1,1} $. Topologically $\TSpace{R}{1,1}$ is equivalent
to the Klein bottle. Similar conformally invariant two-fold cover of the 
Minkowski space-time was constructed in~\cite[\S~III.4]{Segal76} in 
connection with the red shift problem in extragalactic astronomy.

We define \intro{(conformal) unit disk} in $ \TSpace{R}{1,1} $ as follows:
\begin{equation}
\TSpace{D}{}=\{\n{u} \such \n{u}^2<-1,\ \n{u}\in \mathbb{R}^{1,1}_+ 
\}
\cup \{\n{u} \such \n{u}^2>-1,\ \n{u}\in \mathbb{R}^{1,1}_- \}.
\end{equation}
It can be shown that $ \TSpace{D}{}$ is conformally invariant and 
has a boundary $\TSpace{T}{}$---the two glued copies of unit circles in $ 
\mathbb{R}^{1,1}_+$ and $\mathbb{R}^{1,1}_-$. 

We call $\TSpace{T}{}$ 
the \intro{(conformal) unit circle} in $\Space{R}{1,1}$. $\TSpace{T}{}$ 
consists of four parts---branches of hyperbola---with subgroup $A\in\SL$ 
acting simply transitively on each of them. Thus we will regard 
$\TSpace{T}{}$ as $ \Space{R}{} \cup \Space{R}{} \cup\Space{R}{} 
\cup\Space{R}{}$ with an exponential mapping $\exp: t \mapsto (+ \mathrm{ 
or } -) e_1^{+ \mathrm{ or } -}$, $e_1^\pm \in \mathbb{R}^{1,1}_\pm$, where
each of four possible sign combinations is realized on a particular copy of
$ \Space{R}{}$. More generally we define a set of concentric circles for 
$-1\leq\lambda<0$:
\begin{equation} \label{eq:circle-l}
\TSpace{T}{\lambda}=\{\n{u} \such \n{u}^2=-\lambda^2,\ \n{u}\in 
\mathbb{R}^{1,1}_+
\}
\cup \{\n{u} \such \n{u}^2=-\lambda^{-2},\ \n{u}\in \mathbb{R}^{1,1}_- \}.
\end{equation}

Figure~\ref{pi:disk} illustrates geometry of the conformal unit disk in 
$\TSpace{R}{1,1}$ as well as the ``left'' half plane conformally equivalent
to it.

Matrices of the form 
\begin{displaymath}
\matr{\n{a}}{0}{0}{\n{a}'}=\matr{e^{e_1e_2\tau}}{0}{0}{e^{e_1e_2\tau}}, 
\qquad \n{a}={e^{e_1e_2\tau}}=\cosh\tau+e_1e_2\sinh\tau, \quad \tau\in 
\Space{R}{}
\end{displaymath}
comprise a subgroup of $\SL$ which we denote by $A$.  This subgroup is an
image of the subgroup $A$ in the Iwasawa decomposition
$\SL=ANK$~\cite[\S~III.1]{Lang85} under the
transformation~\eqref{eq:caley}.

We define an embedding $s$ of $\TSpace{D}{}$ for our realization of 
$\SL$ by the formula:
\begin{equation} \label{eq:def-s-b}
s: \n{u} \mapsto \frac{1}{ \sqrt[]{1+\n{u}^2} } \matr{1}{\n{u}}{-\n{u}}{1}.
\end{equation}
The formula $g: \n{u} \mapsto s^{-1}(g\cdot s(\n{u}))$ associated with a 
matrix $g^{-1}=\matr{\n{a}}{\n{b}}{-\n{b}}{\n{a}}$ gives the 
fraction-linear transformation $ \TSpace{D}{} \rightarrow \TSpace{D}{}$ of 
the form:
\begin{equation} \label{eq:fr-lin-b}
g: \n{u} \mapsto g\cdot \n{u}= \frac{\n{a}\n{u}+\n{b}}{-\n{b}\n{u}+\n{a}}, 
\qquad g^{-1}=\matr{\n{a}}{\n{b}}{-\n{b}}{\n{a}}
\end{equation}
The mapping $ r: G \rightarrow H$ associated to $s$ defined 
in~\eqref{eq:def-s-b} is
\begin{equation} \label{eq:def-r-b}
r: \matr{\n{a}}{\n{b}}{-\n{b}}{\n{a}} \mapsto \matr{ \frac{\n{a}}{ 
\modulus{\n{a}} }}{0}{0}{\frac{\n{a}}{ \modulus{\n{a}} }}
\end{equation}

And finally the invariant measure on $ \TSpace{D}{} $ 
\begin{equation} \label{eq:def-m-b}
d\mu(\n{u})= \frac{d\n{u}}{{(1+\n{u}^2)}^2} = \frac{du_1 
du_2}{{(1-u_1^2+u_2^2)}^2}.
\end{equation}
follows from the elegant consideration in~\cite[\S~6.1]{Cnops94a}. $\eoe$
\end{exampleb}
We hope the reader notes the explicit similarity between these two 
examples.
Following examples will explore it further.

\subsection[Reduced Wavelet Transform]{Reduced Wavelet 
Transform---the Cauchy Integral Formula}
\begin{examplea}
We continue to consider the case of $G=\SL$ and $H=K$. The compact group
$K\sim \Space{T}{}$ has a discrete set of characters $ \chi_m(h_\phi)=
e^{-i m\phi} $, $m\in \Space{Z}{} $. We drop the trivial character $\chi_0$
and remark that characters $\chi_m$ and $ \chi_{-m} $ give similar
holomorphic and \emph{anti}holomorphic series of representations.  Thus we
will consider only characters $\chi_m$ with $m=1,2,3,\ldots$.

There is a difference in behavior of characters $\chi_1$ and $ \chi_m $ for 
$m=2,3,\ldots$ and we will consider them separately.

First we describe $\chi_1$. Let us take $X=\Space{T}{}$---the unit circle 
equipped with the standard Lebesgue measure $ d\phi $ normalized in such a
way that 
\begin{equation} \label{eq:n-measure}
\int_{ \Space{T}{} } \modulus{f_0(\phi)}^2\, d\phi=1 \textrm{ with }
f_0(\phi)\equiv 1. 
\end{equation}
From ~\eqref{eq:def-s-a} and~\eqref{eq:def-r-a} one can find that
\begin{displaymath}
r(g^{-1}*s(e^{i\phi}))= \frac{ \bar{\beta}e^{i\phi}+\bar{\alpha}}
{ \modulus{ \bar{\beta}e^{i\phi}+\bar{\alpha}} }, 
\qquad
g^{-1}=\matr{\alpha}{\beta}{ \bar{\beta} }{ \bar{\alpha} }.
\end{displaymath}
Then the action of $G$ on $\Space{T}{}$ defined by~\eqref{eq:fr-lin-a}, 
the equality  $ {d(g\cdot \phi)}/{d\phi}= \modulus{\bar{\beta} e^{i\phi} + 
\bar{\alpha} }^{-2} $  and the character $ \chi_1 $ give the following
formula: 
\begin{equation} \label{eq:g-transform}
[\pi_1(g) f](e^{i\phi})= \frac{1}{ \bar{\beta} e^{i\phi} + \bar{ \alpha }} 
f \left( \frac{  { \alpha }e^{i\phi}+{\beta}}{\bar{\beta} e^{i\phi} 
+ \bar{ \alpha }} \right).
\end{equation}
This is a unitary representation---the \intro{mock discrete 
series} of $ \SL$~\cite[\S~8.4]{MTaylor86}. It is easily seen that $K$ acts 
in a trivial way by multiplication by $\chi(e^{i\phi})$.
The function $ f_0(e^{i\phi})\equiv 1 $ mentioned in~\eqref{eq:n-measure}
transforms as follows
\begin{equation} \label{eq:t-vacuum}
[\pi_1(g) f_0](e^{i\phi})= \frac{1}{ \bar{\beta} e^{i\phi} + \bar{ \alpha 
}}
\end{equation}
and in particular has an obvious property
$[\pi_1(h_\psi)f_0](\phi)=e^{i\psi} f_0(\phi)$, i.e.  it is a \intro{vacuum
vector} with respect to the subgroup $H$.  The smallest linear subspace $
\FSpace{F}{2}(X) \in \FSpace{L}{2}(X) $ spanned by~\eqref{eq:t-vacuum}
consists of boundary values of analytic functions in the unit disk and is
the \intro{Hardy space}.  Now the reduced wavelet
transform~\eqref{eq:berg-cauch} takes the form
\begin{eqnarray}
\widehat{f}(a)=[\oper{W} f] (a) & =  &
\scalar{f(x)}{\pi_1(s(a))f_0(x)}_{\FSpace{L}{2}(
X ) } \nonumber \\
& = & \int_{\Space{T}{}} f(e^{i\phi}) \frac{ \sqrt{ 1-\modulus{a}^2 }
}{ \overline{ \bar{a}e^{i\phi} + 1} }\,d\phi \nonumber \\
& = & \frac{\sqrt{ 1-\modulus{a}^2 }}{i}  \int_{\Space{T}{}} 
\frac{f(e^{i\phi})}{{a}+e^{i\phi}} ie^{i\phi}\,d\phi 
\nonumber \\
& = & \frac{\sqrt{ 1-\modulus{a}^2 }}{i}
 \int_{\Space{T}{}}  \frac{f(z)}{{a}+z}\,dz, \label{eq:cauchy1}
\end{eqnarray}
where $z=e^{i\phi}$. Of course~\eqref{eq:cauchy1} is the \intro{Cauchy
integral formula} up to factor ${2\pi } {\sqrt{ 1-\modulus{a}^2 }} $.
Thus we will write $f(a)={\left({2\pi \sqrt{ 1-\modulus{a}^2 
}}\right)}^{-1} \widehat{f}(-a) $ for analytic extension of $f(\phi)$ to
the unit disk.  The factor $2\pi $ is due to our
normalization~\eqref{eq:n-measure} and $\sqrt{ 1-\modulus{a}^2 }$ is
connected with the invariant measure on $ \Space{D}{} $.

Let us now consider characters $\chi_m$ ($m=2,3,\ldots$). 
These characters together with action~\eqref{eq:fr-lin-a} of $G$ give following
representations:
\begin{equation} \label{eq:rep-dis}
[\pi_m(g) f](w)=f\left(  \frac{\alpha w + \beta}{\bar{\beta} w + 
\bar{\alpha}} \right) {(\bar{\beta} w + \bar{\alpha} )}^{-m}.
\end{equation}
For any integer $m\geq 2$ one can select a measure 
\begin{displaymath}
d\mu_m(w)=4^{1-m}{(1- \modulus{w}^2 )}^{m-2} dw,
\end{displaymath}
where $dw$ is the standard Lebesgue measure on $ \Space{D}{} $, such
that~\eqref{eq:rep-dis} become unitary
representations~\cite[\S~IX.3]{Lang85},~\cite[\S~8.4]{MTaylor86}. 
These are discrete series.
 
If we again select $f_0(w)\equiv 1$ then
\begin{displaymath}
[\pi_m(g) f_0](w)= {(\bar{\beta} w + \bar{\alpha} )}^{-m}.
\end{displaymath}
In particular $[\pi_m(h_\phi) f_0](w)= e^{im\phi} f_0(w)$ so this  
again is a vacuum vector with respect to $K$. The irreducible subspace $ 
\FSpace{F}{2}( \Space{D}{} )$ generated by $f_0(w)$ consists of analytic
functions and  is the $m$-th Bergman space (actually \person{Bergman}
considered only $m=2$). Now the transformation~\eqref{eq:berg-cauch}
takes the form                                                                 \begin{eqnarray*}
  \widehat{f}(a) & = & \scalar{f(w)}{[\pi_m(s(a)) f_0](w)}\\
  & = & {\left(\sqrt{1- \modulus{a}^2 }\right)}^m \int_{ \Space{D}{} }
  \frac{f(w)}{{(a\bar{w}+1)}^m} \frac{dw}{{(1- \modulus{w}^2 )}^{2-m}},
\end{eqnarray*}
which for $m=2$ is the classical Bergman formula up to factor 
${\left(\sqrt{1- \modulus{a}^2 }\right)}^m$. Note that calculations in 
standard approaches are ``rather lengthy and must be done in
stages''~\cite[\S~1.4]{Krantz82}. $\eoe$
\end{examplea}
\begin{exampleb}
Now we consider the same group $G=\SL$ but pick up another subgroup 
$H=A$. Let $e_{12}:=e_1e_2$. It follows from~\eqref{eq:p-prop} that 
the mapping from the subgroup $A\sim \Space{R}{}$ to
even numbers\footnote{See Appendix~\ref{pt:bivec-fun} for a definition of 
functions of
even Clifford numbers.} $\chi_\sigma: a \mapsto a^{e_{12}\sigma}= 
{(\exp (e_1e_2 \sigma\ln{a}))}= (a \n{p}_1 + a^{-1} \n{p}_2) ^\sigma$ 
parametrized by $\sigma \in \Space{R}{}$ is a character (in a somewhat
generalized sense). 
It represents an isometric rotation of $\TSpace{T}{}$ by the angle $a$.

Under the present conditions we have from~\eqref{eq:def-s-b} 
and~\eqref{eq:def-r-b}:
\begin{displaymath}
r(g^{-1}*s(\n{u}))=\matr{ \frac{-\n{b}\n{u}+\n{a}}{ 
\modulus{-\n{b}\n{u}+\n{a}} } }{0}{0}{ \frac{-\n{b}\n{u}+\n{a}}{ 
\modulus{-\n{b}\n{u}+\n{a}} } }, \qquad 
g^{-1}=\matr{\n{a}}{\n{b}}{-\n{b}}{\n{a}}.
\end{displaymath}
If we again introduce the exponential coordinates $t$ on $\TSpace{T}{}$ 
coming from the subgroup $A$ (i.e., $\n{u}= e_1 e^{e_1e_2t}\cosh t e_1 -
\sinh t e_2= (x+ \frac{1}{x}) e_1 - (x- \frac{1}{x}) e_2$, $x=e^t$) then
the measure $dt$ on $\TSpace{T}{}$ will satisfy the transformation
condition
\begin{displaymath}
\frac{d(g\cdot t)}{dt}= \frac{1}{(be^{-t}+a)(ce^t+d)}= 
\frac{1}{(-\n{b}\n{u}+\n{a})(\n{u}\n{b}-\n{a})},  
\end{displaymath}
where
\begin{displaymath}
g^{-1}=\matr{a}{b}{c}{d}=\matr{\n{a}}{\n{b}}{-\n{b}}{\n{a}}.
\end{displaymath}

Furthermore we can construct a representation on the 
functions defined on $ \TSpace{T}{} $ by the formula:
\begin{equation} \label{eq:ind-b}
[\pi_\sigma(g) f](\n{v})= 
\frac{(- \n{v} \n{b} + \bar{\n{a}})^{\sigma}}
{(- \n{b} \n{v} + {\n{a}})^{1+\sigma}}
f\left( \frac{\n{a}\n{v}+\n{b}}{-\n{b}\n{v} +\n{a}} \right), \quad 
g^{-1}=\matr{\n{a}}{\n{b}}{-\n{b}}{\n{a}}.
\end{equation}
It is induced by the character $\chi_{\sigma}$ due to formula
$-\n{b} \n{v}+\n{a}= (cx+d)\n{p}_1+(bx^{-1}+a)\n{p}_2$, where $x=e^t$ and 
it is a cousin of the principal series representation (see~\cite[\S~VI.6,
Theorem~8]{Lang85},~\cite[\S~8.2, Theorem~2.2]{MTaylor86} 
and Appendix~\ref{pt:principal}).
The subspaces of vector valued and even number valued functions are invariant 
under~\eqref{eq:ind-b} and the representation is unitary with respect to 
the following inner product (about Clifford valued inner product
see~\cite[\S~3]{Cnops94a}):
\begin{displaymath}
\scalar{f_1}{f_2}_{ \TSpace{T}{} } = \int_{ \TSpace{T}{} }
\bar{f}_2(t) f_1(t)\, dt.
\end{displaymath}
We will denote by $\FSpace{L}{2}(\TSpace{T}{})$ the space of
$\Space{R}{1,1}$-even Clifford number valued functions on $\TSpace{T}{}$
equipped with the above inner product.

We select function $f_0(\n{u})\equiv 1$ neglecting the fact that it does 
not belong to $ \FSpace{L}{2}( \TSpace{T}{} )$. Its transformations
\begin{equation} \label{eq:coher-b}
f_g(\n{v})=[\pi_\sigma(g) f_0](\n{v})= \modulus{1+\n{u}^2}^{1/2}
\frac{(- \n{v} \n{b} + \bar{\n{a}})^{\sigma}}{(- \n{b} \n{v} + 
{\n{a}})^{1+\sigma}}
\end{equation}
and in particular the identity $[\pi_\sigma(g)f_0](\n{v})=
\bar{\n{a}}^\sigma {\n{a}}^{-1-\sigma} f_0(\n{v})=\n{a}^{-1-2\sigma}
f_0(\n{v})$ for $g^{-1}=\matr{\n{a}}{0}{0}{\n{a}}$ demonstrates that it is
a vacuum vector.  Thus we define the reduced wavelet transform accordingly
to~\eqref{eq:def-s-b} and~\eqref{eq:berg-cauch} by the formula\footnote{This 
formula is not well defined in the Hilbert spaces setting.  Fortunately it
is possible (see~\ref{sec:wavel-banach-spac} and ~\cite{Kisil98a}) 
to define a theory of wavelets in Banach spaces
in a way very similar to the Hilbert space case.  So we will ignore this 
difference
in this lecture.}:
\begin{eqnarray}
[\oper{W}_\sigma f] (\n{u})                        
& = & \int_{ \TSpace{T}{} } \modulus{1+\n{u}^2}^{1/2} \overline{ \left(
\frac{(-e_1e^{e_{12} t} \n{u} + \n{1})^\sigma}
{(-\n{u}e_1 e^{e_{12} t} + \n{1})^{1+\sigma}}
\right) } f(t) \,dt \nonumber \\
& = &  \modulus{1+\n{u}^2}^{1/2} \int_{ \TSpace{T}{} }  
\frac{(-\n{u} e_1 e^{e_{12} t} + \n{1})^\sigma}
{(-e^{-e_{12} t} e_1 \n{u} + \n{1})^{1+\sigma}}
f(t) \,dt \label{eq:sio1} \\
& = &  \modulus{1+\n{u}^2}^{1/2} \int_{ \TSpace{T}{} }  
\frac{(-\n{u} e_1 e^{e_{12} t} +\n{1}  )^\sigma }
{e^{-e_{12} t(1+\sigma)}(- e_1 \n{u} + e^{e_{12} t})^{1+\sigma}}
f(t) \,dt \nonumber \\
& = &  \modulus{1+\n{u}^2}^{1/2} \int_{ \TSpace{T}{} }  
\frac {(-\n{u} e_1 e^{e_{12} t} +\n{1} )^\sigma} 
{(- e_1 \n{u} + e^{e_{12} t})^{1+\sigma}} e^{e_{12} t(1+\sigma)}
f(t) \,dt \nonumber \\
& = &  \modulus{1+\n{u}^2}^{1/2} e_{12}\int_{ \TSpace{T}{} }  
\frac{(-\n{u} e_1 e^{e_{12} t} +\n{1} )^\sigma }
{(- e_1 \n{u} + e^{e_{12} t})^{1+\sigma}}e^{e_{12} t\sigma}
( e_{12} e^{e_{12} t} \,dt )\,
f(t)  \nonumber \\
& = &  \modulus{1+\n{u}^2}^{1/2} e_{12}\int_{ \TSpace{T}{} }  
\frac{(-\n{u} e_1 \n{z}  + \n{1} )^\sigma \n{z}^{\sigma} }
{(- e_1 \n{u} + \n{z})^{1+\sigma}} 
 \,d\n{z} \, f(\n{z})  \label{eq:cauchy-b}
\end{eqnarray}
where $\n{z}=e^{e_{12}t}$ and $d\n{z}=e_{12}e^{e_{12}t}\,dt$ are the new 
monogenic variable and its differential respectively. The 
integral~\eqref{eq:cauchy-b} is a singular one, its four singular 
points are intersections of the light cone with the origin in $\n{u}$ with 
the unit circle $\TSpace{T}{}$. See Appendix~\ref{pt:sio} about the 
meaning of this singular integral operator.

The explicit similarity between~\eqref{eq:cauchy1} and~\eqref{eq:cauchy-b}
allows us to consider transformation $\oper{W}_\sigma$~\eqref{eq:cauchy-b} as
an analog of the Cauchy integral formula and the linear space $
\FSpace{H}{\sigma}( \TSpace{T}{} ) $~\eqref{eq:hardy-b} generated by the
coherent states $f_\n{u}(\n{z})$~\eqref{eq:coher-b} as the correspondence
of the Hardy space.  Due to ``indiscrete'' (i.e. they are not square 
integrable) nature of principal series representations there are no
counterparts for the Bergman projection and Bergman space.  $\eoe$ 
\end{exampleb}

\subsection[The Dirac and Laplace Operators]{The 
Dirac (Cauchy-Riemann) and Laplace Operators}
Consideration of Lie groups is hardly possible without consideration of 
their Lie algebras, which are naturally represented by left and right
invariant vectors fields on groups. On a homogeneous space $\Omega=G/H$ we
have also defined a left action of $G$ and can be interested in left
invariant vector fields (first order differential operators). Due to the 
irreducibility of $ \FSpace{F}{2}( \Omega )$ under left action of $G$ every
such vector field $D$ restricted to $ \FSpace{F}{2}( \Omega )$ is a scalar
multiplier of identity $D|_{\FSpace{F}{2}( \Omega )}=cI$. We are
in particular interested in the case $c=0$.
\begin{defn} \cite{AtiyahSchmid80,KnappWallach76}
A $G$-invariant first order differential operator 
\begin{displaymath}
D_\tau: \FSpace{C}{\infty}(\Omega, \mathcal{S} \otimes V_\tau ) 
\rightarrow
\FSpace{C}{\infty}(\Omega, \mathcal{S} \otimes V_\tau )
\end{displaymath}
such that $\oper{W}(\FSpace{F}{2}(X))\subset \object{ker} D_\tau$ is called 
\intro{(Cauchy-Riemann-)Dirac operator} on $\Omega=G/H$
associated with an irreducible representation $ \tau $ of $H$ in a space 
$V_\tau$ and a spinor bundle $\mathcal{S}$.  
\end{defn}
The Dirac operator is explicitly defined by the 
formula~\cite[(3.1)]{KnappWallach76}:
\begin{equation} \label{eq:dirac-def}
D_\tau= \sum_{j=1}^n \rho(Y_j) \otimes c(Y_j) \otimes 1,
\end{equation}
where $Y_j$ is an orthonormal basis of
$\algebra{p}=\algebra{h}^\perp$---the orthogonal completion of the Lie
algebra $\algebra{h}$ of the subgroup $H$ in the Lie algebra $\algebra{g}$
of $G$; $\rho(Y_j)$ is the infinitesimal generator of the right action of
$G$ on $\Omega$; $c(Y_j)$ is Clifford multiplication by $Y_i \in
\algebra{p}$ on the Clifford module $\mathcal{S}$.  We also define 
an invariant Laplacian by the formula
\begin{equation} \label{eq:lap-def}
\Delta_\tau= \sum_{j=1}^n \rho(Y_j)^2 \otimes \epsilon_j \otimes 1,
\end{equation}
where $\epsilon_j = c(Y_j)^2$ is $+1$ or $-1$.
\begin{prop}
Let all commutators of vectors of $ \algebra{h}^\perp $ belong to $
\algebra{h}$, i.e.
$[\algebra{h}^\perp,\algebra{h}^\perp]\subset\algebra{h}$.  Let also $f_0$
be an eigenfunction for all vectors of $ \algebra{h} $ with eigenvalue $0$
and let also $\oper{W}f_0$ be a null solution to the Dirac operator $D$.
Then $\Delta f(x)=0$ for all $f(x)\in \FSpace{F}{2}(\Omega)$.
\end{prop}
\begin{proof}
Because $\Delta$ is a linear operator and $ \FSpace{F}{2}(\Omega)$ is 
generated by $\pi_0(s(a))\oper{W} f_0$ it is enough to check that $
\Delta \pi_0(s(a))\oper{W} f_0=0 $. Because $\Delta$  and $\pi_0$ commute 
it is enough to check that $\Delta \oper{W} f_0=0$. Now we observe that
\begin{displaymath}
\Delta = D^2 - \sum_{i,j} \rho([Y_i,Y_j]) \otimes c(Y_i)c(Y_j) \otimes 1.
\end{displaymath}
Thus the desired assertion is follows from two identities: $D 
\oper{W}f_0=0$ and
$\rho([Y_i,Y_j])\oper{W}f_0=0$, $[Y_i,Y_j]\in H$.
\end{proof}
\begin{examplea}
Let $G=\SL$ and $H$ be its one-dimensional compact subgroup generated by 
an element $Z \in \algebra{sl}(2, \Space{R}{})$. Then 
$\algebra{h}^\perp$ is spanned by two vectors $Y_1=A$ and $Y_2=B$. In such 
a situation we can use $ \Space{C}{} $ instead of the Clifford algebra $ 
\Cliff[0]{2} $. Then formula~\eqref{eq:dirac-def} takes a simple 
form $D=r(A+iB)$. Infinitesimal action of this operator in the upper-half 
plane follows from calculation in~\cite[VI.5(8), IX.5(3)]{Lang85}, it is 
$[D_{ \Space{H}{} } f] (z)= -2i y \frac{ \partial f(z)}{ \partial \bar{z} }
$, $z=x+iy$. Making the Caley transform we can find its action in the
unit disk $D_{ \Space{D}{} } $: again the Cauchy-Riemann operator $ \frac{ 
\partial }{ \partial \bar{z} } $ is its principal component. 
We calculate  $D_{ \Space{H}{} }$ explicitly now to stress the 
similarity with $\Space{R}{1,1}$ case. 

For the upper half plane $\Space{H}{}$ we have following formulas:
\begin{eqnarray*}
s&:&\Space{H}{} \rightarrow \SL: z=x+iy \mapsto 
g=\matr{y^{1/2}}{xy^{-1/2}}{0}{y^{-1/2}}; \nonumber \\
s^{-1}&:& \SL \rightarrow \Space{H}{}: \matr{a}{b}{c}{d} \mapsto z= 
\frac{ai+b}{ci+d};  \nonumber \\
\rho(g)&:& \Space{H}{} \rightarrow\Space{H}{} :z \mapsto s^{-1}( s(z) * g) 
\\
&& \qquad \qquad \qquad  =s^{-1}\matr{ay^{-1/2}+cxy^{-1/2}}{ 
by^{1/2}+dxy^{-1/2}}{cy^{-1/2}}{dy^{-1/2}}\\
&& \qquad \qquad \qquad =\frac{(yb+xd)+i(ay+cx)}{ci+d} 
\nonumber
\end{eqnarray*}
Thus the right action of $\SL$ on $\Space{H}{}$ is given by the formula
\begin{displaymath}
\rho(g)z=\frac{(yb+xd)+i(ay+cx)}{ci+d}= x+y \frac{bd+ac}{c^2+d^2}+iy 
\frac{1}{c^2+d^2}.
\end{displaymath}
For $A$ and $B$ in $\algebra{sl}(2, \Space{R}{} )$ we have:
\begin{displaymath}
\rho(e^{At}) z = x+iy e^{2t}, \qquad \rho(e^{Bt}) z = x 
+y\frac{e^{2t}-e^{-2t}}{e^{2t}+e^{-2t}}+iy\frac{4}{e^{2t}+e^{-
2t}}.
\end{displaymath}
Thus
\begin{eqnarray*}
[\rho(A) f](z) & = & \frac{ \partial f (\rho(e^{At}) z )}{ \partial t }
|_{t=0} = 2y \partial_2 f(z), \\ {}
[\rho(B) f](z) & = & \frac{ \partial f (\rho(e^{Bt}) z )}{ \partial t } 
|_{t=0} = 2y \partial_1 f(z),
\end{eqnarray*}
where $\partial_1$ and $\partial_2$ are derivatives of $f(z)$ with respect 
to real and imaginary party of $z$ respectively. Thus we get 
\begin{displaymath}
D_{ \Space{H}{} }= i\rho(A) + \rho( B) = 2y i\partial_2 + 2y \partial_1=
2y \frac{ \partial }{ \partial \bar{z} }
\end{displaymath}
as was expected. $\eoe$
\end{examplea}
\begin{exampleb}
In $\Space{R}{1,1}$ the element $B\in \algebra{sl}$ generates the subgroup 
$H$ and its orthogonal completion is spanned by $B$ and $Z$. Thus the
associated Dirac operator has the form $D=e_1 \rho(B) + e_2 \rho(Z)$. 
We need infinitesimal generators of the right action $\rho$ on the ``left'' 
half plane $\TSpace{H}{}$. Again we have a set of formulas similar to the 
classic case:
\begin{eqnarray*}
s&:&\TSpace{H}{} \rightarrow \SL: \n{z}=e_1y+e_2x \mapsto 
g=\matr{y^{1/2}}{xy^{-1/2}}{0}{y^{-1/2}}; \nonumber \\
s^{-1}&:& \SL \rightarrow \TSpace{H}{}: \matr{a}{b}{c}{d} \mapsto \n{z}= 
\frac{ae_1+be_2}{ce_2e_1+d};  \nonumber \\
\rho(g)&:& \TSpace{H}{} \rightarrow \TSpace{H}{} :\n{z} \mapsto s^{-1}( 
s(\n{z}) * g) \\
&&\qquad \qquad \qquad =s^{-1}\matr{ay^{-1/2}+cxy^{-1/2}}{ 
by^{1/2}+dxy^{-1/2}}{cy^{-1/2}}{dy^{-1/2}}\\
&&\qquad \qquad \qquad = \frac{(yb+xd)e_2+(ay+cx)e_1}{ce_2e_1+d}
\nonumber
\end{eqnarray*}
Thus the right action of $\SL$ on $\Space{H}{}$ is given by the formula
\begin{displaymath}
\rho(g)\n{z}=\frac{(yb+xd)e_2+(ay+cx)e_1}{ce_2e_1+d}= 
e_1 y \frac{-1}{c^2-d^2} + e_2 x+e_2y \frac{ac-bd}{c^2-d^2}.
\end{displaymath}
For $A$ and $Z$ in $\algebra{sl}(2, \Space{R}{} )$ we have:
\begin{eqnarray*}
\rho(e^{At}) \n{z}& = & e_1 y e^{2t}+ e_2 x, \\
\rho(e^{Zt}) \n{z}& = &e_1y\frac{-1}{\sin^2 t - \cos^2 t}
+e_2 y\frac{-2\sin t \cos t}{\sin^2 t -\cos^2 t}+ e_2 x\\
& =& e_1y\frac{1}{ \cos 2 t}+ e_2 y \tan 2 t + e_2 x .
\end{eqnarray*}
Thus
\begin{eqnarray*}
[\rho(A) f](\n{z}) & = & \frac{ \partial f (\rho(e^{At}) \n{z} )}{ \partial 
t }
|_{t=0} = 2y \partial_2 f(\n{z}), \\ {}
[\rho(Z) f](\n{z}) & = & \frac{ \partial f (\rho(e^{Zt}) \n{z} )}{ \partial 
t }
|_{t=0} = 2y \partial_1 f(\n{z}),
\end{eqnarray*}
where $\partial_1$ and $\partial_2$ are derivatives of $f(\n{z})$ with 
respect of $e_1$ and $e_2$ components of $\n{z}$ respectively. Thus we get 
\begin{displaymath}
D_{ \TSpace{H}{} }= e_1\rho(Z) + e_2\rho( A) = 2y (e_1\partial_1 + e_2
\partial_2).
\end{displaymath}
In this case the Dirac operator is not elliptic and as a consequence we
have in particular a singular Cauchy integral formula~\eqref{eq:cauchy-b}.
Another manifestation of the same property is that primitives in the
``Taylor expansion'' do not belong to $\FSpace{F}{2}(\TSpace{T}{})$ itself
(see Example~\ref{ex:taylor-b}).  It is known that solutions of a
hyperbolic system (unlike the elliptic one) essentially depend on the
behavior of the boundary value data.  Thus function theory in
$\Space{R}{1,1}$ is not defined only by the hyperbolic Dirac equation alone
but also by an appropriate boundary condition.  
$\eoe$
\end{exampleb}

\subsection{The Taylor expansion}
For any decomposition $f_a(x)=\sum_\alpha  \psi_\alpha(x) V_\alpha(a)$ 
of the coherent states $f_a(x)$ by means of functions $V_\alpha(a)$
(where the sum can become eventually an integral) we have the
\intro{Taylor expansion} 
\begin{eqnarray} 
\widehat{f}(a) & = & \int_X f(x) \bar{f}_a(x)\, dx= \int_X f(x) \sum_\alpha 
\bar{\psi}_\alpha(x)\bar{V}_\alpha(a)\, dx  \nonumber \\
 & = &  \sum_\alpha 
\int_X f(x)\bar{\psi}_\alpha(x)\, dx \bar{V}_\alpha(a) \nonumber \\
 & = & \sum_{\alpha}^{\infty} \bar{V}_\alpha(a) f_\alpha,\label{eq:taylor1}
\end{eqnarray}
where $f_\alpha=\int_X f(x)\bar{\psi}_\alpha(x)\, dx$.
However to be useful within the presented scheme such a decomposition 
should be connected with the structures of $G$, $H$, and the representation 
$\pi_0$. We will use a decomposition of $f_a(x)$ by the eigenfunctions of 
the operators $\pi_0(h)$, $h\in \algebra{h}$.
\begin{defn}
 Let $\FSpace{F}{2}=\int_{A} \FSpace{H}{\alpha}\,d\alpha$ be a spectral 
decomposition with respect to the operators $\pi_0(h)$, $h\in \algebra{h}$.
Then the decomposition
\begin{equation} \label{eq:spec-c}
 f_a(x)= \int_{A} V_\alpha(a) f_\alpha(x)\, d\alpha,
\end{equation}
where $f_\alpha(x)\in \FSpace{H}{\alpha}$ and $V_\alpha(a): 
\FSpace{H}{\alpha} \rightarrow \FSpace{H}{\alpha}$ is called the Taylor 
decomposition of the Cauchy kernel $f_a(x)$.
\end{defn}
Note that the Dirac operator $D$ is defined in the terms of left invariant 
shifts and therefor commutes with all $\pi_0(h)$. Thus it also has a 
spectral decomposition over spectral subspaces of $\pi_0(h)$:
\begin{equation} \label{eq:spec-d}
 D= \int_{A} D_\delta \, d\delta.
\end{equation}
We have obvious property
\begin{prop} \label{pr:cauchy-dirac}
If spectral measures $d\alpha$ and $d\delta$ 
from~\eqref{eq:spec-c} and~\eqref{eq:spec-d} have disjoint supports then 
the image of the Cauchy integral belongs to the kernel of the Dirac 
operator.
\end{prop}
For discrete series representation functions $f_\alpha(x)$ can be 
found in $\FSpace{F}{2}$ (as in Example~\ref{ex:taylor-a}), for the 
principal series representation this is not the case. To overcome confusion
one can think about the Fourier transform on the real line. It can be 
regarded as a continuous decomposition of a function $f(x)\in 
\FSpace{L}{2}(\Space{R}{})$ over a set of harmonics $e^{i\xi x}$ neither of
those belongs to $\FSpace{L}{2}(\Space{R}{})$. This has a lot of common 
with the Example~\ref{ex:taylor-b}.
\begin{examplea} \label{ex:taylor-a}
Let $G=\SL$ and $H=K$ be its maximal compact subgroup and $\pi_1$ be 
described by~\eqref{eq:g-transform}.  $H$ acts on $\Space{T}{}$ by
rotations.  It is one dimensional and eigenfunctions of its generator $Z$
are parametrized by integers (due to compactness of $K$).  Moreover, on the
irreducible Hardy space these are positive integers $n=1,2,3\ldots$ and
corresponding eigenfunctions are $f_n(\phi)=e^{i(n-1)\phi}$. Negative 
integers span the space of anti-holomorphic function and the splitting 
reflects the existence of analytic structure given by the Cauchy-Riemann 
equation. The decomposition of coherent states $f_a(\phi)$ by means of this
functions is well known:
\begin{displaymath}
f_a(\phi)= \frac{ \sqrt[]{1- \modulus{a}^2 }}{ \bar{a}e^{i\phi}-1}= 
\sum_{n=1}^\infty \sqrt[]{1- \modulus{a}^2 }\bar{a}^{n-1} e^{i(n-1)\phi}=
\sum_{n=1}^\infty V_n(a)f_n(\phi),
\end{displaymath}
where $V_n(a)=\sqrt[]{1- \modulus{a}^2 }\bar{a}^{n-1} $. This is the 
classical Taylor expansion up to multipliers coming from the invariant 
measure. $\eoe$
\end{examplea}
\begin{exampleb} \label{ex:taylor-b}
Let $G=\SL$, $H=A$, and  $\pi_\sigma$ be described by~\eqref{eq:ind-b}.
Subgroup $H$ acts on $\TSpace{T}{}$ by hyperbolic rotations: 
\begin{displaymath}
\tau: \n{z}=e_1 e^{e_{12}t} \rightarrow e^{2e_{12}\tau}\n{z}=e_1
e^{e_{12}(2\tau+t)}, \qquad t, \tau \in \TSpace{T}{}.
\end{displaymath}
Then for every $p\in \Space{R}{}$ the function 
$f_p(\n{z})=(\n{z})^p=e^{e_{12}pt}$ where $\n{z}=e^{e_{12}t}$ is an
eigenfunction in the representation~\eqref{eq:ind-b} for a generator $a$ 
of the subgroup $H=A$ with the eigenvalue $2(p-\sigma)-1$.  Again, due to
the analytical structure reflected in the Dirac operator, we only need
negative values of $p$ to decompose the Cauchy integral kernel.
\begin{prop} \label{pr:taylor-b}
For $\sigma=0$ the Cauchy integral kernel~\eqref{eq:cauchy-b} has the 
following decomposition:
\begin{equation}
\frac{1}{-e_1\n{u}+\n{z}}=
\int_0^\infty \frac{(e_1\n{u})^{[p]}-1}{e_1\n{u}-1}  
\cdot t\n{z}^{-p}\,dp,
\end{equation}
where $\n{u}=u_1 e_1 +u_2 e_2$, $\n{z}=e^{e_{12}t}$,  and $[p]$
is the integer part of $p$ (i.e. $k=[p]\in \Space{Z}{} $, $k\leq p < k+1$).
\end{prop}
\begin{proof}
Let 
\begin{displaymath}
f(t)=\int_0^\infty F(p) e^{-tp}\, dp
\end{displaymath}
be the Laplace transform. We use the formula~\cite[Laplace Transform Table,
p.~479, (66)]{CRCMTables} 
\begin{equation} \label{eq:lap-table}
\frac{1}{t(e^{kt}-a)}= \int_0^\infty \frac{a^{[p/k]}-1}{a-1} e^{-tp}\, dp
\end{equation}
with the particular value of the parameter $k=1$.
Then using $\n{p}_{1,2}$ defined in~\eqref{eq:p-def} we have
\begin{eqnarray}
\lefteqn{
\int_0^\infty \frac{(e_1\n{u})^{[p]}-1}{e_1\n{u}-1}  
\cdot t\n{z}^{-p}\,dp =}\  && \nonumber \\
& = & t\int_0^\infty \left(\frac{(-u_1-u_2)^{[p]}-1}{(-u_1-u_2)-1} \n{p}_2 
+  \frac{(-u_1+u_2)^{[p]}-1}{(-u_1+u_2)-1} \n{p}_1 
\right)  (e^{tp}\n{p}_2+ e^{-tp}\n{p}_1)\,dp \nonumber \\
& = & t\int_0^\infty  \frac{(-u_1-u_2)^{[p]}-1}{(-u_1-u_2)-1}e^{tp}\,dp\, 
\n{p}_2 +  
t\int_0^\infty \frac{(-u_1+u_2)^{[p]}-1}{(-u_1+u_2)-1}e^{-tp}\,dp\, \n{p}_1 
\nonumber \\
& = & \frac{t}{t(e^{-t} +u_1+u_2)} \n{p}_2 + 
\frac{t}{t(e^{t} +u_1-u_2)} \n{p}_1 \label{eq:lap-appl} \\
& = & \frac{1}{(e^{-t} +u_1+u_2) \n{p}_2 + 
(e^{t} +u_1-u_2) \n{p}_1 }\nonumber \\
&=& \frac{1}{-e_1\n{u}+\n{z}} \nonumber,
\end{eqnarray}
where we obtain~\eqref{eq:lap-appl} by an application 
of~\eqref{eq:lap-table}.
\end{proof}
Thereafter for a function $f(\n{z})\in \FSpace{F}{2}(\TSpace{T}{})$ we 
have the following Taylor expansion of its wavelet transform:
\begin{displaymath}
[\oper{W}_0 f](u)= \int_0^\infty \frac{(e_1\n{u})^{[p]}-1}{e_1\n{u}-1} 
f_p \,dp,
\end{displaymath}
where
\begin{displaymath}
f_p = \int_{\TSpace{T}{}} t \n{z}^{-p}\,d\n{z} f(\n{z}).
\end{displaymath}
The last integral is similar to the Mellin 
transform~\cite[\S~III.3]{Lang85}, \cite[Chap.~8, (3.12)]{MTaylor86}, 
which naturally arises in study of the principal series representations 
of $\SL$. 

I was pointed by Dr.~J.~Cnops that for the Cauchy kernel $(-e_1 
\n{u}+\n{z})$ there is still a decomposition of the form $(-e_1 
\n{u}+\n{z})=\sum_{j=0}^\infty (e_1 \n{u})^j \n{z}^{-j-1}$. In this 
connection one may note that representations $\pi_1$~\eqref{eq:g-transform}
and $\pi_\sigma$~\eqref{eq:ind-b} for $\sigma=0$ are unitary equivalent. 
(this is a meeting point between discrete and principal series). Thus 
a function theory in $\Space{R}{1,1}$ with the value $\sigma=0$ could 
carry many properties known from the complex analysis.
$\eoe$
\end{exampleb}

\section{Open problems} \label{pt:o-problems}
This lecture raises more questions than gives answers. Nevertheless it is 
useful to state some open problems explicitly.
\begin{enumerate}
\item Demonstrate that Cauchy formula~\eqref{eq:cauchy-b} is an isometry 
between $\FSpace{F}{2}(\TSpace{T}{})$ and $\FSpace{H}{\sigma}(\TSpace{D}{})$
with suitable norms chosen. This almost follows (up to some constant 
factor) from its property to intertwine two irreducible representations of 
$\SL$.
\item Formula~\eqref{eq:sio1} contains Szeg\"o type kernel, which is domain 
dependent. Integral formula~\eqref{eq:cauchy-b} formulated in terms of 
analytic kernel. Demonstrate using Stocks theorem that~\eqref{eq:cauchy-b} is 
true for other suitable chosen domains.
\item The image of Szeg\"o (or Cauchy) type formulas belong to the kernel 
of Dirac type operator only if they connected by additional condition 
(see Proposition~\ref{pr:cauchy-dirac}). Descriptive
condition for the discrete series can be found 
in~\cite[Theorem~6.1]{KnappWallach76}. Formulate a similar condition for 
principal series representations.
\end{enumerate}

\chapter[Segal-Barmann Spaces]{Segal-Barmann 
Spaces and Nilpotent Groups}

\newif\ifcoimbra
\coimbratrue
\def\ip#1#2{#1\cdot#2}
\def\pa#1{\partial_{#1}}
\def\rnuln{\R_{0n}}
\def\Rnuln{\R^{0n}}
\def\L2{{\FSpace{L}{2}(\Space{R}{n})}}
\def\BB{{\cal B}}
\def\FF{{\cal F}}
\def\hlie{\mathfrak{h}}
\def\CLie{\Space C{}\mathfrak{h}^n}
\def\CH{\Space C{}\Heisen n}\def\Lie{\mathfrak{h}^n}
\providecommand{\SL}{SL(2,\Space{R}{})}
\def\nota#1{\ \newline\hspace*{-2cm} !!!!!\hspace*{1cm}{\it #1}\newline}
\def\quork#1{{\ifmmode{{\cal C}\kern-0.18em\ell(#1)}\else
${{\cal C}\kern-0.18em\ell(#1)}$\fi}}

\def\ort{_{\perp}}
\newcommand{\uu}{{\bf u}}
\newcommand{\vv}{{\bf v}}
\renewcommand{\a}{{\bf a}}
\renewcommand{\b}{{\bf b}}
\renewcommand{\c}{{\bf c}}
\newcommand{\p}{{\bf p}}
\newcommand{\q}{{\bf q}}
\newcommand{\w}{{\bf w}}
\newcommand{\myx}{{\bf x}}
\newcommand{\y}{{\bf y}}
\newcommand{\z}{{\bf z}}
\newcommand{\vu}{\nu}
\newcommand{\SB}{\FSpace{F}{2}( \Space{C}{n} )}
\newcommand{\mult}[1]{\mathsf{#1}}

\providecommand{\hyperref}[2][\relex]{#2}

\section{Introduction}

This lecture is based on the paper~\cite{CnopsKisil97a}.

It is well known, by the celebrated Stone-von Neumann theorem,
that all models for the canonical quantisation~\cite{Mackey63} are
isomorphic and provide us with equivalent representations of 
the \hyperref[item:Heisenberg-group]{Heisenberg
group}~\cite[Chap.~1]{MTaylor86}. Nevertheless it is worthwhile to look
for some models which can act as alternatives for the Schr\"odinger
representation. In particular, the Segal-Bargmann
representation~\cite{Bargmann61,Segal60} serves to 
\begin{itemize}
\item  give a geometric 
  representation of the dynamics of the
  \extref{sp-funct}{sec:sov-quantum-problem}{harmonic oscillators};
\item present a nice model
  for the creation and annihilation operators, which is important for
  quantum field theory;
\item allow applying \extref{sp-funct}{sec:complex-contour-beta}{tools
    of analytic function theory}.
\end{itemize}
The huge abilities of the Segal-Bargmann (or Fock~\cite{Fock32}) model are
not yet completely employed, see for example new ideas in a recent
preprint~\cite{NazaikSternin96}.

We look for similar connections between nilpotent Lie groups and spaces of 
monogenic~\cite{BraDelSom82,DelSomSou92} Clifford valued functions. 
Particularly we are interested in a
third possible representation of the Heisenberg group, acting on monogenic
functions on \Space{R}{n}.  There are several reasons why such a model
can be of interest.  First of all the theory of monogenic functions is
(at least) as interesting as several complex variable theory, so the
monogenic model should share many pleasant features with the Segal-Bargmann
model.  Moreover, monogenic functions take their value in a Clifford
algebra, which is a natural environment in which to represent internal
degrees of freedom of elementary particles such as spin.  Thus from the
very beginning it has a structure which in the Segal-Bargmann model has to
be added, usually by means of the second quantization
procedure~\cite{Dirac67}.  So a monogenic representation can be even more
relevant to quantum field theory than the Segal-Bargmann one (see
Remark~\ref{re:quant-field}).

From the different aspects of the Segal-Bargmann space $ \SB $ we select the
one giving a unitary representation of the
\hyperref[item:Heisenberg-group]{Heisenberg group} 
$\Space{H}{n}$.  The
representation is unitary equivalent to the Schr\"odinger representation on
$ \FSpace{L}{2}( \Space{R}{n} ) $ and the Segal-Bargmann transform is
precisely the intertwining operator between these two representations 
(see subsection~\ref{ss:segal}).

\ifcoimbra
Monogenic functions can be introduced in this
scheme in two ways, as either $ \FSpace{L}{2}( \Space{R}{n} ) $ or $\SB$
can be substituted by a space of monogenic functions. 

In the first case one defines a new unitary irreducible representation of
the Heisenberg group on a space of monogenic functions
and constructs an analogue of the Segal-Bargmann transform as the intertwining
operator of the new representation and the Segal-Bargmann one. We examine this
possibility in section~\ref{se:heisenberg}. In a certain sense the
representation of the Heisenberg group constructed here lies between
Schr\"odinger and Segal-Bargmann ones, taking properties both of them.

In the second case we first select a substitute for the
Heisenberg group, so we can replace the Segal-Bargmann space by a
space of monogenic functions. The space $\Space{C}{n}$
underlying $ \FSpace{F}{2}( \Space{C}{n} ) $ is 
intimately connected with the structure of the Heisenberg group $ 
\Space{H}{n}$ in the sense that
$\Space{C}{n} $ is the quotient of $ \Space{H}{n} $ with respect to its 
centre. In order to define a space of monogenic functions, say on $
\Space{R}{n+1} $, we have to construct a group playing a similar
r\^ole with respect to this space. We describe an option in
section~\ref{se:new}. 

Finally we give the basics of coherent states from square integrable group
representations and an interpretation of the classic Segal--Bargmann space 
in terms of these in Appendix~\ref{sec:elem-repr-theory}.
\fi

This lecture is closely related to~\cite{Kisil97a}, where connections between
analytic function theories and group representations were described. 
Representations of another group ($\SL$) in spaces of monogenic functions 
can be found in~\cite{Kisil97c}. We hope that the present lecture make only
few first steps towards an interesting function theory and other steps will
be done elsewhere.

\ifcoimbra
\else
\section{Wavelets or Coherent States}
\label{sec:wavelets-or-coherent}

In our approach we will need some basic facts on \intro{wavelets} (or
\intro{coherent states}) and associated \intro{wavelet transform}.

Let $G$ be a group which acts via transformation of a closed domain
$\bar{\Omega}$. Moreover, let $G: \partial \Omega\rightarrow \partial
\Omega$ and $G$ act on $\Omega$ and $\partial \Omega$ transitively.
Let us fix a point $x_0\in \Omega$ and let $H\subset G$ be a
stationary subgroup of point $x_0$. Then domain $\Omega$ is naturally
identified with the  homogeneous space $G/H$. Till the moment we do
not request anything untypical. Now let 
\begin{itemize}
\item\emph{there exist a $H$-invariant measure $d\mu$ on $\partial
    \Omega$}.
\end{itemize}
We consider the Hilbert space $\FSpace{L}{2}(\partial
\Omega, d\mu)$. Then geometrical transformations of $\partial \Omega$
give us the representation $\pi$ of $G$ in $\FSpace{L}{2}(\partial
\Omega, d\mu)$.
 Let $f_0(x)\equiv 1$ and $\FSpace{F}{2}(\partial
\Omega, d\mu)$ be the closed liner subspace of $\FSpace{L}{2}(\partial
\Omega, d\mu)$ with the properties:
\begin{enumerate}
\item\label{it:begin} $f_0\in \FSpace{F}{2}(\partial \Omega, d\mu)$;
\item $\FSpace{F}{2}(\partial \Omega, d\mu)$ is $G$-invariant;
\item\label{it:end} $\FSpace{F}{2}(\partial \Omega, d\mu)$ is
  $G$-irreducible, or $f_0$ is cyclic in $\FSpace{F}{2}(\partial
  \Omega, d\mu)$. 
\end{enumerate}
The \introind{standard wavelet transform}{wavelet transform} 
$W$ is defined by
\begin{displaymath}
W: \FSpace{F}{2}(\partial \Omega, d\mu) \rightarrow
\FSpace{L}{2}(G): f(x) \mapsto
\widehat{f}(g)=\scalar{f(x)}{\pi(g)f_0(x)}_{\FSpace{L}{2}(\partial
\Omega,d\mu) }
\end{displaymath}
Due to the property $[\pi(h)f_0](x)=f_0(x)$, $h\in H$ and 
identification $\Omega\sim G/H$ it could be translated to the embedding:
\begin{equation}\label{eq:cauchy}
\widetilde{W}: \FSpace{F}{2}(\partial \Omega, d\mu) \rightarrow
\FSpace{L}{2}(\Omega): f(x) \mapsto
\widehat{f}(y)=\scalar{f(x)}{\pi(g)f_0(x)}_{\FSpace{L}{2}(\partial
\Omega,d\mu) },  
\end{equation}

We define the \introind{inverse wavelet transform}{wavelet
  transform!inverse} $ \oper{M} $ according to the formula:
\begin{equation} 
[ \oper{M} \widehat{f}](x) =  \int_{ \Omega } \widehat{f}(a) 
f_{s(a)}(x) \, da \label{eq:a-inverse} ,
\end{equation}

The following proposition explain the usage of the name for
$\oper{M}$.
\begin{thm}
  The operator
\begin{equation} \label{eq:szego1}
\oper{P}= \oper{M} \oper{W}: B \rightarrow B
\end{equation}
is a projection of $B$ to its linear subspace for which $b_0$ is
cyclic. Particularly if $\pi$ is an irreducible representation then the
inverse wavelet transform $\oper{M}$ is a \introind{left inverse
  operator}{operator!left inverse}
on $B$ for the wavelet transform $\oper{W}$:
\begin{eqnordisp}
  \oper{M}\oper{W}=I.
\end{eqnordisp}
\end{thm}

\fi

\section[The Heisenberg Group]{The Heisenberg Group 
and Spaces of Analytic Functions} 
\label{se:heisenberg}
\subsection[The Schrodinger Representation]{The 
Schr\"odinger Representation of the heisenberg Group}
We recall here some basic facts on the Heisenberg group $\Space{H}{n}$
and its Schr\"odinger representation, see~\cite[Chap.~1]{Folland89}
and~\cite[Chap.~1]{MTaylor86} for details.

The Lie algebra of the Heisenberg group is generated by the $2n+1$ elements
$p_1$, \ldots, $p_n$, $q_1$, \ldots, 
$q_n$, $e$, with the well-known Heisenberg commutator relations:
\begin{equation}\label{eq:a-heisenberg}
[p_i,q_j]=\delta_{ij} e.
\end{equation}
All other commutators vanish. In the standard quantum mechanical
interpretation the operators are momentum and coordinate 
operators~\cite[\S~1.1]{Folland89}.

It is common practice to switch between real and complex Lie algebras.
Complexify $\Lie$ to obtain the complex algebra $\CLie$, and take
four complex numbers  $a$, $b$, $c$ and $d$ such that $ad-bc\not=0$. The
{\it real\/} $2n+1$-dimensional subspace spanned by
$$
A_k=ap_k+bq_k\quad\quad B_k=cp_k+dq_k
$$
and the commutator $[A_k,B_k]=(ad-bc)e$, where $e=[p_k,q_k]$ is of course
isomorphic to $\Lie$, and exponentiating will give a group isomorphic
to the Heisenberg group.

An example of this procedure is obtained from the construction of the
so-called creation and annihilation operators of Bose particles in
the $k$-th state, $a^+_k$ and  $a^{-}_k$ (see~\cite[\S~1.1]{Folland89}).
These are defined by:
\begin{equation}\label{eq:creation}
a^{\pm}_k=\frac{q_k \mp \imath p_k}{\sqrt[]{2}},
\end{equation}
giving the commutators $[a^+_i,a^-_j]=(-\imath)\delta_{ij}e$. Putting
$-\imath e=\ell$, the real algebra spanned by $a^\pm_k$ and $\ell$ is an
alternative realization of $\Lie$, $\Lie_a$.

An element $g$ of the Heisenberg group $\Heisen{n}$ (for any positive
integer $n$, cf.~\eqref{eq:heisenberg-def}) can be represented 
as $g=(t,\z)$ with $t\in\Space{R}{}$,
$\z=(z_1,\ldots,z_n)\in \Space{C}{n}$. The group law in coordinates
$(t,\z)$ is given by
\begin{equation}\label{eq:g-heisenberg}
g*g'=(t,\z)*(t',\z')=(t+t'+\frac{1}{2}\sum_{j=1}^n\Im(\bar{z}_j
z_j'), \z+\z'),
\end{equation}
where $\Im z$ denotes the imaginary part of the complex number $z$.
Of course the Heisenberg group is non-commutative.

The relation between the Heisenberg group and its Lie algebra is given
by the exponentiation $\exp:\Lie_a\to\Heisen n$. We define the formal
vector $\a^+$ as being $(a^+_1,\ldots,a^+_n)$ and $\a^-$ as
$(a^-_1,\ldots,a^-_n)$, which allows us to use the formal inner products
\begin{eqnarray*}
\uu\cdot \a^+&=&\sum_{k=1}^nu_ka^+_k\\
\vv\cdot \a^-&=&\sum_{k=1}^nv_ka^-_k.
\end{eqnarray*}
With these we define, for real vectors $\uu$ and $\vv$, and real $s$
\begin{eqnarray}
\exp(\uu\cdot(\a^++\a^-))&=&(0,\sqrt 2\uu)\\
\exp(\vv\cdot (\a^--\a^+)&=&(0,\imath\vv)\\
\exp(s \ell)&=&(e^{-2s},0).
\end{eqnarray}

Possible Schr\"odinger representations 
\ifcoimbra
\else
(cf.~\eqref{eq:Schroding1})  
\fi
are parameterized by the
non-zero real number $\hbar$ (the Planck constant). As usual, for
considerations where the correspondence principle between classic and
quantum mechanics is irrelevant, we consider only the case 
$\hbar=1$. The Hilbert space for the Schr\"odinger representation is
$\FSpace{L}{2} (\Space{R}{n})$, where elements of the complex Lie
algebra $\CLie$ are represented by the unbounded operators
\begin{equation} \label{eq:rep-schro}
 \sigma(a^{\pm}_k)=
\frac{1}{\sqrt[]{2}}\left(x_k I \mp \frac{\partial
}{\partial x_k}\right).
\end{equation}
From which it follows, using any $j$, that
$$
\sigma(\ell)=[a^+_j,a^-_j]=-2 I.
$$
The corresponding representation $\pi$ of the Heisenberg group is
given by exponentiation of the $\sigma(a^+_k)$ and $\sigma(a^-_k)$,
but this is most readily expressed by using $p_k$ and $q_k$, and so is
generated by shifts and multiplications $s_\c: f(\myx) \mapsto f(\myx+\c)$ and
$m_\b: f(\myx)\mapsto e^{\imath\ip\myx\b}f(\myx)$, with the Weyl commutation
relation
\begin{displaymath}
s_\c m_\b = e^{\imath\ip\c\b} m_\b s_\c.
\end{displaymath}

There is an orthonormal basis of $\FSpace{L}{2}(\Space{R}{n})$ on
which the operators $\sigma(a^\pm_k)$ act in an especially simple way. It
consists of the functions:
\begin{equation}\label{eq:hbasis}
\phi_m(\y)=[2^m m! \,\, \sqrt[]{\pi}]^{-1/2} e^{-\ip\myx\myx/2}H_m(\y),
\end{equation}
where $\y=(y_1,\ldots,y_n)$, $m=(m_1,\ldots,m_n)$, and $H_m(\y)$ is the
generalized \extref{sp-funct}{sec:hermite-polynomials}{Hermite polynomial}
$$
H_m(\y)=\prod_{i=1}^nH_{m_i}(y_i).
$$
For these
\begin{displaymath}
a^{+}_k \phi_m(\y)=\sqrt{m_k+1}\,\phi_{m'}(\y), \qquad a^{-}_k \phi_m(\y)=
\sqrt{m_k}\,\phi_{m''}(\y)
\end{displaymath}
where 
\begin{eqnarray*}
m'&=&(m_1,m_2,\ldots, m_{k-1}, m_k + 1, m_{k+1},\ldots,
m_{n})\\
m''&=&(m_1,m_2,\ldots, m_{k-1}, m_k - 1, m_{k+1},\ldots,
m_{n}).
\end{eqnarray*}
This is the most straightforward way to express the creation or
annihilation of a particle in the $k$-th state.

Let us now consider the
\extref{sp-funct}{eq:hermit-generate}{generating function} of the
$\phi_m(\myx)$, 
\begin{equation}\label{eq:sum}
A(\myx,\y)=\sum_{j=0}^\infty \frac{x^j}{\sqrt[]{j!}} \phi_k(\y)= \exp(-
\frac{1}{2}(\ip\myx\myx+\ip\y\y) + \sqrt[]{2}\ip\myx\y).
\end{equation}
We state the following elementary fact in Dirac's bra-ket notation.
\begin{lem}\label{le:dirac}
Let $H$ and $H'$ be two Hilbert spaces with orthonormal bases
$\{\phi_k\}$ and $\{\phi'_k\}$ respectively. Then the sum
\begin{equation}\label{eq:operator}
U=\sum_{j=0}^\infty \ket{\phi'_j} \bra{\phi_j}
\end{equation}
defines a unitary operator $U: H \rightarrow H'$ with the following
properties:
\begin{enumerate}
\item $U \phi_k = \phi_k'$;
\item If an operator $T: H\rightarrow H$ is expressed, relative to the
basis $\phi_k$, by the matrix $(a_{ij})$ then the operator $UTU^{-1}: H'
\rightarrow H'$ is expressed relative to the basis $\phi_k'$ by the same
matrix. 
\end{enumerate}
\end{lem}

Now, if we take the function $A(\myx,\y)$ from~\eqref{eq:sum} as a kernel
for an \extref{sp-funct}{sec:introduction5}{\intro{integral transform}},
\begin{displaymath}
[Af](\y)=\int_{\Space{R}{n}} A(\y,\myx) f(\myx)\,dx
\end{displaymath}
we can consider it subject to the Lemma above. However, for this we need
to define the space $H'$ and an orthonormal basis $\{\phi_k'\}$  (we
already identified $H$ with $\FSpace{L}{2}(\Space{R}{n})$ and the
$\{\phi_k\}$ are given by~\eqref{eq:hbasis}). There is some freedom in
doing this.

For example it is possible to take the holomorphic extension $A(\z,\y)$
of $A(\myx,\y)$ with respect to the first variable.
Then
\begin{enumerate}
\item $H'$ is the Segal-Bargmann space of analytic functions over
$\Space{C}{n} $ with scalar product defined by the integral with respect to
Gaussian measure $e^{-\modulus{\z}^2}\,d\z$;
\item The Heisenberg group acts on the Segal-Bargmann space as follows:
\begin{equation} \label{eq:rep-barg}
[\beta_{(t,\z)}f](\uu)=f(\uu+\z)e^{\imath t-\ip{\bar{\z}}\uu-\modulus{\z}^2/2}.
\end{equation}
This action generates the set of coherent states 
$f_{(0,\vv)}(\uu)=e^{-\bar{\vv}\uu-\modulus{\vv}^2/2}$, $\uu$,
$\vv\in \Space{C}{n}$ from the vacuum vector $f_0(\uu)\equiv 1$;
\item The operators of creation and annihilation are $a^+_k=z_k I$,
$a^-_k=\frac{\partial }{\partial z_k}$.
\item The Segal-Bargmann space is spanned by the
orthonormal basis $\phi_k'=\frac{1}{\sqrt[]{m!}}z^n$ or by the set of
coherent states $f_{(0,\vv)}(\uu)=e^{-\bar{\vv}\uu-\modulus{\vv}^2/2}$, $\uu$,
$\vv\in \Space{C}{n}$
\item The intertwining kernel for $\sigma_{(t,\z)}$~\eqref{eq:rep-schro} and
$\beta_{(t,\z)}$~\eqref{eq:rep-barg} is
\begin{displaymath}
A(\z,\y)=e^{-(\ip\z\z+\ip\myx\myx)/2-\sqrt[]{2}\ip\z\myx}
=\sum_{k=0}^\infty \frac{\z^m}{\sqrt[]{m!}} \cdot
\frac{1}{\sqrt{2^m m!} \sqrt[4]{\pi}} e^{-\ip\myx\myx/2}H_m(\y)
\end{displaymath}
\item The Segal-Bargmann space has a reproducing kernel
\begin{displaymath}
K(\uu,\vv)= e^{\ip\uu{\bar{\vv}}}=\sum_{k=1}^\infty \phi_k(\uu)
\bar{\phi}_k(\vv)= \int
e^{\ip\uu{\bar{\z}}} e^{\ip\z{\bar{\vv}}} e^{-\modulus{\z}^2}\,d\z.
\end{displaymath}
\end{enumerate}

\subsection{The Segal-Bargmann space} \label{ss:segal}
We consider a representation of the Heisenberg group $ \Space{H}{n}
$ (see Section~\ref{se:heisenberg}) on $ \FSpace{L}{2}( \Space{R}{n} ) $ by
shift and multiplication operators~\cite[\S~1.1]{MTaylor86}:
\begin{equation} \label{eq:schrodinger}
g=(t,\z): f(\myx) \rightarrow [\pi_{(t,\z)}f](\myx)=
e^{\imath(2t-\sqrt{2}\ip\q\myx+\ip\q\p)}
f(\myx- \sqrt 2\p), \qquad \z=\p+\imath q,
\end{equation}
This is the Schr\"odinger representation with parameter $\hbar=1$. As a
subgroup $H$ we select the centre of \Space{H}{n} consisting of elements
$(t,0)$. It is non-compact but using the special form of
representation~\eqref{eq:schrodinger} we can consider the
cosets\footnote{$ \widetilde{G}$ is sometimes called the \intro{reduced
Heisenberg group}{Heisenberg group!reduced}. 
It seems that $ \widetilde{G}$ is a virtual object,
which is important in connection with a selected representation of $G$.}
$ \widetilde{G} $ and $ \widetilde{H} $ of $G$ and $H$ by the subgroup
with elements $(\pi m,0)$, $m\in \Space{Z}{}$.
Then~\eqref{eq:schrodinger} also defines a representation of $
\widetilde{G} $ and $ \widetilde{H} \sim \Gamma $. We consider the Haar
measure on $ \widetilde{G} $ such that its restriction on $
\widetilde{H} $ has total mass equal to $1$.

As ``vacuum vector'' we will select the original \intro{vacuum
vector} of quantum mechanics---the Gauss function $f_0(\myx)=e^{-\ip\myx\myx/2}$.
Its transformations are defined as follows:
\begin{eqnarray*}
w_g(\myx)=\pi_{(t,\z)} f_0(\myx) & = &
e^{\imath(2t-\sqrt{2}\ip\q\myx+\ip\q\p)}
\,e^{-{(\myx- \sqrt 2\p)}^2/2}\\
& = & e^{2\imath t-(\ip\p\p+\ip\q\q)/2}
e^{- ((\p-\imath \q)^2+\ip\myx\myx)/2+\sqrt{2}\ip{(\p-\imath \q)}\myx}
\\
& = & e^{2\imath t-\ip\z{\bar\z}/2}e^{- (\ip{\bar\z}{\bar\z}+\ip\myx\myx)/2
+\sqrt{2}\ip{\bar\z}\myx}.
\end{eqnarray*}
In particular $w_{(t,0)}(\myx)=e^{-2it}f_0(\myx)$, i.e.\ it really is
a vacuum vector with respect to $\widetilde{H}$ in the sense of our
definition. Of course $\widetilde{G} / \widetilde{H}$ is
isomorphic to  $\Space{C}{n}$. Embedding $\Space Cn$ in $G$ by the
identification of $(0,\z)$ with $\z$, the mapping $s:
\widetilde G \rightarrow \widetilde{G}$ is defined simply by
$s((t,\z))=(0,\z)=\z$; $\Omega$ then is identical with $\Space Cn$.

The Haar measure on $ \Space{H}{n} $ coincides with the standard
Lebesgue measure on $ \Space{R}{2n+1} $~\cite[\S~1.1]{MTaylor86} and so
the invariant measure on $ \Omega $ also coincides with Lebesgue
measure on $\Space{C}{n}$. Note also that the composition law sending
$\z_1$ $\z_2$ to $s((0,\z_1)(0,\z_2))$ reduces to Euclidean shifts on $
\Space{C}{n} $. We also find $s((0,\z_1)^{-1}\cdot (0,\z_2))=\z_2-\z_1$
and $r((0,\z_1)^{-1}\cdot (0,\z_2))= (\frac{1}{2} \Im\ip{\bar\z_1}{\z_2},0)$.

The reduced wavelet transform takes the form of a mapping
$\FSpace{L}{2}(\Space{R}{n} ) \rightarrow \FSpace{L}{2}( \Space{C}{n} ) 
$ and is given by the formula
\begin{eqnarray}
\widehat{\oper W}f(\z)&=&\scalar{f}{w_{(0,\z)}}\nonumber \\
     &=&\pi^{-n/4}\int_{\Space{R}{n}} f(\myx)\, e^{-\ip\z{\bar\z}/2}\,e^{-
(\ip\z\z+\ip\myx\myx)/2+\sqrt{2}\ip\z\myx}\,dx \nonumber \\
     &=&e^{-\modulus{\z}^2/2}\pi^{-n/4}\int_{\Space{R}{n}} f(\myx)\,e^{-
(\ip\z\z+\ip\myx\myx)/2+\sqrt{2}\ip\z\myx}\,dx, \label{eq:tr-bargmann}
\end{eqnarray}
where $\z=\p+\imath\q$. Then $\widehat{\oper W}f$ belongs to
$\FSpace{L}{2}( \Space{C}{n} , dg)$. This can better be expressed by
saying that the 
function $\breve{f}(\z)=e^{\modulus{\z}^2/2}\widehat{\oper W}f(\z)$ belongs
to $\FSpace{L}{2}( \Space{C}{n} , e^{- \modulus{\z}^2 }dg)$
because $\breve{f}(\z)$ is analytic in $\z$. These functions constitute the
\intro{Segal-Bargmann space}~\cite{Bargmann61,Segal60}  
$ \FSpace{F}{2}( \Space{C}{n}, e^{-
\modulus{\z}^2 }dg) $ of functions analytic in $\z$ and
square-integrable with respect the Gaussian measure $e^{-
\modulus{\z}^2}d\z$. Analyticity of $\breve{f}(\z)$ is equivalent to
the condition that $( \frac{ \partial }{ \partial\bar{\z}_j } + \frac{1}{2} \z_j
I ) \oper Wf(\z)$ equals zero.

The integral in~\eqref{eq:tr-bargmann} is the well-known
Segal-Bargmann transform~\cite{Bargmann61,Segal60}. Its inverse is 
given by a realization of~\eqref{eq:a-inverse}:
\begin{eqnarray}
f(\myx) & = & \int_{ \Space{C}{n} } \widehat{\oper Wf(\z)} w_{(0,\z)}(\myx)\,d\z
\nonumber\\
& = & \int_{ \Space{C}{n} } \breve{f}(\z)  e^{-
(\bar{\z}^2+\ip\myx\myx)/2+\sqrt{2}\bar{\z}x}\, e^{- \modulus{\z}^2}\, d\z.
\end{eqnarray}
This gives~\eqref{eq:a-inverse} the name of Segal-Bargmann 
inverse. The corresponding operator $\oper{P}$~\eqref{eq:szego1} is the
identity operator $ \FSpace{L}{2}(\Space{R}{n}) \rightarrow 
\FSpace{L}{2}(\Space{R}{n}) $ and~\eqref{eq:szego1} gives an integral
presentation of the Dirac delta. 

Meanwhile the orthoprojection  
$ \FSpace{L}{2}( \Space{C}{n},  e^{- \modulus{\z}^2 }dg)  \rightarrow
\FSpace{F}{2}( \Space{C}{n},  e^{- \modulus{\z}^2 }dg) $ is of interest
and is a principal ingredient in Berezin
quantisation~\cite{Berezin88,Coburn94a}. We can easy find its kernel. 
Indeed, $ \widehat{\oper
W}f_0(\z)=e^{-\modulus{\z}^2}$, and the kernel is
\begin{eqnarray*}
K(\z,\w) & = & \widehat{\oper W}f_0(\z^{-1}\cdot \w)
\bar{\chi}(r(\z^{-1}\cdot \w))\\
& = & \widehat{\oper W}f_0(\w-\z)\exp(\imath\Im(\ip{\bar\z}\w) \\
& = & \exp(\frac{1}{2}(- \modulus{\w-\z}^2 +\ip\w{\bar\z}-\ip\z{\bar\w}))\\
& = &\exp(\frac{1}{2}(- \modulus{\z}^2- \modulus{\w}^2) +\ip\w{\bar\z}).
\end{eqnarray*}
To obtain the reproducing kernel for functions
$\breve{f}(\z)=e^{\modulus{\z}^2} \widehat{\oper W}f(\z) $ in the
Segal-Bargmann space we multiply $K(\z,\w)$ by $e^{(-\modulus{\z}^2+
\modulus{\w}^2)/2}$ which gives the standard reproducing kernel, $\exp(-
\modulus{\z}^2 +\ip\w{\bar\z})$ \cite[(1.10)]{Bargmann61}.

The Segal-Bargmann space is an interesting and important object, but
there are also other options. In particular we can consider an
alternative representation of the Heisenberg group, this time acting on
monogenic functions, an action we introduce in the next subparagraph.

\subsection[Spaces of Monogenic 
Functions]{Representation of $\Space{H}{n}$ in Spaces of Monogenic 
Functions}\label{ss:cl-al}
We consider the real Clifford algebra \Cliff{n}, i.e. the algebra
generated by $e_0=1$, ${e_j}$, ${1\leq j \leq n}$, using the identities:
\begin{displaymath}
e_i e_j + e_j e_i = -2 \delta_{ij}, \qquad 1\leq i,j\leq n.
\end{displaymath}
For a function $f$ with values in $\Cliff n$, the action of the Dirac
operator of $\Space{R}{n+1}$ is defined by (here $x=x_0+\myx$ is the $n+1$
dimensional variable)
$$
Df(x)=\sum_{i=0}^n\partial_if(x).
$$
A function $f$ satisfying $Df=0$ in a certain domain is called monogenic
there; later on we shall use the term `monogenic' for solutions of more
general Dirac operators. Obviously the notion of monogenicity is closely
related to the one of holomorphy on the complex plane. As a matter of
fact $D^2=-\Delta$, and monogenic functions are solutions of the
Laplacian. The Clifford algebra is not commutative, and so it is
necessary to introduce a symmetrized product. For $k$ elements $a_i$,
$1\leq i\leq k$ of the algebra it is defined by
$$
a_1\times a_2\times\ldots\times a_k=
{1\over k!}\sum_\sigma a_{\sigma(1)}a_{\sigma(2)}\ldots a_{\sigma(n)},
$$
where the sum is taken over all possible permutations of $k$ elements.
If the same element appears several times, we use an exponent notation,
e.g.\ $a^2\times b^3=a\times a\times b\times b\times b$.

Let now $V_k$ be the symmetric power monomial defined by the expression
\begin{equation}
V_k(\myx)=\frac{1}{\sqrt[]{k!}}(e_1 x_0-
e_0 x_1)^{k_1} \times (e_2 x_0- e_0 x_2)^{k_2} \times \cdots \times
(e_n x_0- e_0 x_n)^{k_n}.
\end{equation}
It can be proved that these monomials are all monogenic (see e.g.\
\cite{Malonek93}), and even that they constitute a basis for the space of
monogenic polynomials (as a module over $\Cliff n$). In general the
symmetrized product is not associative, and manipulating it can become
quite formal. However, if we restrict the monomials defined above to the
hyperplane $x_0=0$, we obtain
$$
V_k(x)=\frac{1}{\sqrt[]{k!}}x_1^{k_1}x_2^{k_2}\ldots x_n^{k_n},
$$
and so we have the multiplicative property
$$
\sqrt{k!k'!\over(k+k)!}V_kV_{k'}=V_{k+k'},\quad x_0=0.
$$
Another important function is the monogenic exponential function which
is defined by
$$
E(\uu,x)=\exp(\uu\cdot\myx)\left(\cos(\norm\uu x_0)-
\frac\uu{\norm\uu}\sin(\uu x_0)\right).
$$
It is not hard to check~\cite[\S~14]{BraDelSom82} that this function is
monogenic, and of course its restriction to the hyperplane $x_0=0$ is
simply the exponential function, $E(\uu,\myx)=\exp(\uu\cdot\myx)$.

We can therefore extend the symmetric product by the so-called
Cauchy-Kovalevskaya product~\cite[\S~14]{BraDelSom82}: If $f$ and $g$ are
monogenic in $\Space R{n+1}$,
then $f\times g$ is the monogenic function equal to $fg$ on $\Space R n$.
Introducing the monogenic functions $\myx_{i}=e_ix_0-e_0x_i$ we can then
write 
$$
V_k(x)=\frac{1}{\sqrt{k!}}x_{1}^{k_1}
\times x_{2}^{k_2}\times\ldots\times x_{n}^{k_n}.
$$

It is fairly easy to check the $V_k$ form an orthonormal set with respect
to the following inner product (see~\cite[\S~3.1]{Cnops94a} on Clifford 
valued inner products):
\begin{equation} \label{eq:inner-m2}
\scalar{V_k}{V_{k'}}=\int_{\Space{R}{n+1}} \bar{V}_k( x) V_{k'}
( x)\, e^{- \modulus{x}^2 } \, dx.
\end{equation}
Let $ \FSpace{M}{2}$ be closure of the linear span of $\{V_k\}$, using
complex coefficients.

The creation and annihilation operators $a^+_k$ and $a^-_k$ can be
represented by symmetric multiplication  (see~\cite{Malonek93}) with the
monogenic variable $\myx_j$, which will be written $\myx_k I_\times
$, and by the (classical) partial derivative $\frac{\partial }{\partial
\myx_j}=\frac{\partial }{\partial {x}_j}$ with respect to $\myx_j$,
which appear in the definition of hypercomplex differentiability. On basis
elements they act as follows:
\begin{eqnarray*}
\myx_jI_\times
V_{(k_1,\ldots,k_j,\ldots,k_n)}&=&
\sqrt{k_j+1}V_{(k_1,\ldots,k_j+1,\ldots,k_n)},\\
\frac{\partial }{\partial \myx_j}
V_{(k_1,\ldots,k_j,\ldots,k_n)}&=&\sqrt{k_j} 
V_{(k_1,\ldots,k_j-1,\ldots,k_n)},
\end{eqnarray*}

It can be checked that this really is a representation of $a^\pm_k$, and
that $a^+_k$ and $a^-_k$ are each other's adjoint. We use the equalities
$a^-_j=\frac{1}{ \sqrt[]{2} }(a^+_j + a^-_j) $ and $a^+_j=
\frac{\imath}{ \sqrt 2}(a^-_j - a^+_j)$, and the commutation relations
$[a^+_i,a^-_j]=e\delta_{ij}$ to obtain a representation of the
Heisenberg group. Thus an element $(t,\z)$, $\z=\uu+i\vv$ of the
Heisenberg group can be written as
\begin{eqnarray*}
(t,\z)&=&
\left(t+\frac{\uu\cdot\uu-\vv\cdot \vv}4,0\right)
\left(0,\frac{(1+\imath)(\uu+\vv)}2\right)
\left(0,\frac{(1-\imath)(\uu-\vv)}2\right)\\
&=&
\exp\left(\left(t+\frac{\uu^2-\vv^2}4\right)e\right)
\exp\left(\frac{(\uu+\vv)q}{\sqrt 2}\right)
\exp\left(\frac{(\uu-\vv)\imath p}{\sqrt 2}\right).
\end{eqnarray*}
It is therefore represented by the operator
\begin{eqnarray}
\pi_{(t,\z)}&=&
\exp\left(-\left(t+\frac{\uu\cdot\uu-\vv\cdot\vv}4\right)\right)\nonumber\\
&&\exp\left(\frac{((\uu+\vv)\cdot\myx) I_\times }{\sqrt 2}\right)
\exp\left(\frac{(\uu-\vv)\cdot(\partial_\myx)}{\sqrt 2}\right),
\label{eq:rep-m2}
\end{eqnarray}
where obviously for a monogenic function $f$ we have
\begin{eqnarray*}
\exp\left(\frac{(\uu-\vv)\imath p}{\sqrt 2}\right)f(x)&=&
 f\left(x+\frac{\uu-\vv)}{\sqrt 2}\right)\\
\exp\left(\frac{((\uu+\vv)\cdot\myx) I_\times }{\sqrt 2}\right)f(x)&=&
E\left(\frac{\uu+\vv}{\sqrt 2},\cdot\right)\times f(x)
\end{eqnarray*}
Therefore it is easy to calculate the image of
the constant function $f_0( \myx)= V_0(\myx) \equiv 1$, and we obtain
the set of functions
\begin{eqnarray}
f_{(t,\z)}(\myx)&=&\pi_{(t,\z)} f_0(\myx) \nonumber\\
&=&
\exp\left(-\left(t+\frac{\uu\cdot\uu-\vv\cdot\vv}4\right)\right)
E\left(\frac{\uu+\vv}{\sqrt 2},\cdot\right)\times f_0(x)\nonumber\\
&&\exp\left(-\left(t+\frac{\uu\cdot\uu-\vv\cdot\vv}4\right)\right)
E\left(\frac{\uu+\vv}{\sqrt 2},x\right).
\label{eq:m-coher}
\end{eqnarray}
In the language of quantum physics $f_0(\myx)$ is the \intro{vacuum
vector} and functions $f_{(t,\z)}(\myx)$ are \intro{coherent states} (or
\intro{wavelets}) for the representation of $ \Space{H}{n}$ we described.
We can summarize the properties of the representation:
\begin{enumerate}
\item  All functions in $\FSpace{M}{2} $ are complex-vector
valued, monogenic in $\Space{R}{n+1}$, and square integrable with
respect to the measure $ e^{ - \modulus{x}^2 }dx$.
\item The representation of the Heisenberg group is given by 
\eqref{eq:rep-m2}. This representation generates a set of coherent
states $f_{(0,z)}(\myx)$~\eqref{eq:m-coher} as shifts of the vacuum
vector $f_0(\myx) \equiv 1$.
\item The creation and annihilation operators $a^+_k$ and $a^-_k$ are 
represented by  symmetric (Cauchy-Kovalevskaya) multiplication by $\myx_j$
and by derivation of monogenic functions. They are adjoint with respect to 
the inner product~\eqref{eq:inner-m2}.
\item $ \FSpace{M}{2} $ is generated as a closed linear space by
the orthonormal basis $V_k(\myx)=\frac{1}{\sqrt[]{k!}}(e_1 x_0-
e_0 x_1)^{k_1} \times (e_2 x_0- e_0 x_2)^{k_2} \times \cdots \times
(e_n x_0- e_0 x_n)^{k_n}$, and also by the set of coherent states
$f_{(t,\z)}(\myx)$ of~\eqref{eq:m-coher}.
\item The kernel of the operator intertwining the model constructed here
and the Segal-Bargmann one is given by 
\begin{displaymath}
B(\z,x)\sum_{j=0}^\infty V_j(x)
\frac{\z^j}{\sqrt[]{j!}}
=\exp(\sum_{k=1}^n \myx_k \bar{z}_k),
\end{displaymath}
which is the holomorphic extension in $\z=\uu+\imath\vv$ of
$E(\uu,x)$. The transformation pair is given by
\begin{eqnarray*}
\BB f(x)&=&\int_{\Space Cn}B(\z,x)f(\z)
\exp\left({-|\z|^2\over 2}\right)\, d\z\\
\BB^{-1}\phi(\z)&=&\int_{\Space R{n+1}}\overline{B(\z,x)}\phi(x)
\exp\left({-|x|^2\over 2}\right)\, dx
\end{eqnarray*}
\item The space $\FSpace{M}{2}$ has a reproducing kernel
\begin{displaymath}
K(x,y)=\sum_{k=0}^\infty V_k(x)
\bar{V}_k(y)=\int_{\Space{C}{n}}B(\z,x)\overline{B(\z,y)}
\, e^{-\modulus{z}^2}dz.
\end{displaymath}
Notice that $\overline{K(x,y)}$ is monogenic in $y$; it is the monogenic
extension of $\overline{E(\y,x)}$.
\end{enumerate}

One can see that some properties of $ \FSpace{M}{2}$ are closer to those
of the Segal-Bargmann space than to those of the space $\FSpace{L}{2}(
\Space{R}{n} )$ it replaces. 
It should be noted that the representation of the Heisenberg
group we obtained here is new and quite unexpected.
\begin{rem} \label{re:quant-field}
We construct $\FSpace{M}{2}$ as a space of complex-vector valued functions.
We can also consider an extended space $\FSpace{\widetilde{M}}{2}$ being
generated by the orthonormal basis $V_k(\myx)$ or coherent states
$f_{(0,z)}(\myx)$ with Clifford valued coefficients multiplied from the
right hand side. Such a space will share many properties of $\Space{M}{2}$
and have an additional structure: there is a natural representation 
$s: f(\myx) \mapsto s^* f(s\myx s^*) s$ of $\object{Spin}(n)$ group in
$\FSpace{\widetilde{M}}{2}$. Thus this space provides us with 
a representation of two main symmetries in quantum field theory: the 
Heisenberg group of quantized coordinate and momentum (external degrees of
freedom) and $\object{Spin}(n)$ group of quantified inner degrees of
freedom. Another composition of the Heisenberg group and Clifford algebras 
can be found in~\cite{Kisil93c}.
\end{rem}

\section[Another Nilpotent Lie Group]{Another Nilpotent 
Lie Group and Its Representation} 
\label{se:new}
\subsection[Complex Vectors]{Clifford 
Algebra and Complex Vectors}
\label{ss:cl-alv}

Starting from the real Clifford algebra \Cliff{n}, we consider complex
$n$-vector valued functions defined on the real line \Space{R}{1} with
values in \Space{C}{n}.  Moreover we will look at the $j$-th component of
\Space{C}{n} as being spanned by the elements $1$ and $e_j$ of the Clifford
algebra.  For two vectors $\uu=(u_1,\ldots,u_n)$ and
$\vv=(v_1,\ldots,v_n)$ we introduce the Clifford vector valued product
(see~\cite[\S~3.1]{Cnops94a} on Clifford valued inner products):
\begin{equation}\label{eq:v-pt}
\uu\cdot \vv=\sum_{j=1}^n \bar{u}_j v_j=\sum_{j=1}^n (u'_j-u''_je_j)
(v'_j+v''_je_j),
\end{equation}
where $u_j=u'_j+u''_j e_j$ and $v_j=v'_j+v''_j e_j$.
Of course, $u\cdot u$ coincides with $\norm{u}^2=\sum_1^n (u_j'^2 +
u_j''^2)$, the standard norm in \Space{C}{n}. So we can introduce the space
$\FSpace{R}{2}(\Space{R}{1})$ of \Space{C}{n}-valued functions on the
real line with the product
\begin{equation}\label{eq:f-pt}
\scalar{f}{f'}=\int_{\Space{R}{1}} f(x) \cdot f'(x) \, dx.
\end{equation}
Again $\scalar{f}{f}^{1/2}$ gives us the standard norm in the Hilbert space
of $\FSpace{L}{2} $ integrable \Space{C}{n} valued functions.

\subsection{A nilpotent Lie group}
We introduce a nilpotent Lie group, \Space{G}{n}. As a $C^\infty$-manifold
it coincides with \Space{R}{2n+1}. Its Lie algebra has generators $P$,
$Q_j$, $T_j$, $1\leq j\leq n$. The non-trivial commutators between them
are
\begin{equation}
[P,Q_j]=T_j;
\end{equation}
all others vanish.
Particularly \Space{G}{n} is a step two nilpotent Lie group and the
$T_j$ span its centre. It is easy to see that $\Space{G}{1} $ is just
the Heisenberg group $\Space{H}{1} $.

We denote a point $g$ of \Space{G}{n} by $2n+1$-tuple of reals
$(t_1,\ldots,t_n;p;q_1,\ldots,q_n)$. These are the exponential coordinates
corresponding to the basis of the Lie algebra $T_1$, \ldots, $T_n$, $P$,
$Q_1$, \ldots, $Q_n$. The group law is given in exponential coordinates by
the formula
\begin{eqnarray}\label{eq:G-mult}
\lefteqn{ (t_1,\ldots,t_n;p;q_1,\ldots,q_n)*
(t'_1,\ldots,t'_n;p';q'_1,\ldots,q'_n)=}\nonumber\\
&=&(t_1+t_1'+\frac{1}{2}(p'q_1-pq'_1),\ldots,t_n+t_n'+
\frac{1}{2}(p'q_n-pq'_n);\nonumber\\
&&\qquad p+p';q_1+q_1',\ldots,q_n+q_n').
\end{eqnarray}

We consider the homogeneous space $\Omega=\Space{G}{n}/\Space{Z}{}$. Here
$\Space{Z}{}$ is the centre of $\Space{G}{n}$; its Lie algebra is
spanned by $T_j$, $1\leq j \leq n$. It is easy to see that
$\Omega\sim\Space{R}{n+1}$. We define the mapping $s:\Omega\rightarrow
\Space{G}{n} $ by the rule
\begin{equation}
s: (a_0,a_1,\ldots,a_n) \mapsto (0,\ldots,0;
a_0;a_1,\ldots,a_n).
\end{equation}
It is the ``inverse'' of the natural projection $s^{-1}:\Space{G}{n}
\rightarrow \Omega=\Space{G}{n}/\Space{Z}{}$.

It easy to see that the mapping $\Omega\times\Omega \rightarrow \Omega$
defined by the rule $s^{-1}(s(a)*s(a'))$ is just Euclidean
(coordinate-wise) addition $a+a'$.

To introduce the Dirac operator we will need the following set of 
left-invariant differential operators, which generate right shifts on the 
group:
\begin{eqnarray}
T_j &=& \frac{ \partial }{ \partial t_j }, \label{eq:l-first} \\
P&=& \frac{ \partial }{ \partial p } +\frac 12 \sum_1^n q_j \frac{ \partial }{
\partial t_j }, \\
Q_j &=& -\frac{ \partial }{ \partial q_j } + \frac 12 p \frac{ \partial }{
\partial t_j }. \label{eq:l-last}
\end{eqnarray}
The corresponding set of right invariant vector fields generating left 
shifts is
\begin{eqnarray}
T^*_j&=& \frac{ \partial }{ \partial t_j }, \label{eq:r-first} \\
P^*&=& \frac{ \partial }{ \partial p } - \frac 12 \sum_1^n q_j \frac{ \partial }{
\partial t_j }, \\
Q^*_j&=&- \frac{ \partial }{ \partial q_j } - \frac 12 p \frac{ \partial }{
\partial t_j }. \label{eq:r-last}
\end{eqnarray}
A general property is that any left invariant operator commutes with any 
right invariant one.

\subsection{A representation of \Space{G}{n} }
We introduce a representation $\rho$ of \Space{G}{n} in the space
$\FSpace{R}{2}(\Space{R}{})$ by the formula:
\begin{equation}\label{eq:G-rep}
[\rho_g f](x)=(e^{e_1(2t_1+q_1(\sqrt{2}x-p))}f_1(x-\sqrt{2}p), \ldots,
e^{e_n(2t_n+q_n(\sqrt{2}x-p))}f_n(x-\sqrt{2}p)),
\end{equation}
where $f(x)=(f_1(x),\ldots,f_n(x))$ and the meaning of
$\FSpace{R}{2}(\Space{R}{})$ was discussed in Subsection~\ref{ss:cl-alv}.
We note that the generators $e_j$ of Clifford algebras do not interact
with each other under the representation just defined.
One can check directly that~\eqref{eq:G-rep} defines a representation
of \Space{G}{n}. Indeed:
\begin{eqnarray}
[\rho_g\rho_{g'}f](x)
&=&\rho_g(e^{e_1(2t'_1+q'_1(\sqrt{2}x-p'))}f_1(x-\sqrt{2}p'), \ldots,
\nonumber \\
&&\qquad \qquad e^{e_n(2t'_n+q'_n(\sqrt{2}x-p'))}f_n(x-\sqrt{2}p')) 
\nonumber \\
&=& (e^{e_1(2t_1+q_1(\sqrt{2}x-p))}
e^{e_1(2t'_1+q'_1(\sqrt{2}(x-\sqrt{2}p)-p'))}
f_1(x-\sqrt{2}p-\sqrt{2}p'), \nonumber \\
&& \ldots, \nonumber \\
&& e^{e_n(2t_n+q_n(\sqrt{2}x-p))}
e^{e_n(2t'_n+q'_n(\sqrt{2}(x-\sqrt{2}p)-p'))}
f_n(x-\sqrt{2}p-\sqrt{2}p')) \nonumber \\
&=& (e^{e_1(2(t_1+t'_1+\frac{1}{2}(p'q_1-pq'_1))+
(q_1+q'_1)(\sqrt{2}x-(p+p')))}
f_1(x-\sqrt{2}(p+p')), \nonumber\\
&& \ldots, \nonumber \\
&& (e^{e_n(2(t_n+t'_n+\frac{1}{2}(p'q_n-pq'_n))+
(q_n+q'_n)(\sqrt{2}x-(p+p')))}
f_n(x-\sqrt{2}(p+p')) \nonumber \\
&=& [\rho_{gg'}f](x),
\end{eqnarray}
where $gg'$ is defined by~\eqref{eq:G-mult}.

$\rho_g$ has the important property that it preserves
the product~\eqref{eq:f-pt}. Indeed:
\begin{eqnarray}
\scalar{\rho_gf}{\rho_gf'}&=& \int_{\Space{R}{}} [\rho_g f](x) \cdot
[\rho_gf'](x)  \,dx \nonumber \\
&=& \int_{\Space{R}{}}\sum_{j=1}^n
 \bar{f}_j(x-\sqrt[]{2}p)e^{-e_j(2t_j+q_j(\sqrt[]{2}x-p))} \nonumber \\
&& \qquad\quad\quad e^{e_j(2t_j+q_j(\sqrt[]{2}x-p))} f'_j(x-\sqrt[]{2}p)
 \,dx \nonumber \\
&=& \int_{\Space{R}{}}\sum_{j=1}^n
\bar{f}_j(x-\sqrt[]{2}p)
f'_j(x-\sqrt[]{2}p)  \,dx \nonumber \\
&=& \int_{\Space{R}{}}\sum_{j=1}^n
 \bar{f}_j(x)  f'_j(x)  \,dx \nonumber \\
&=&\scalar{f}{f'}. \nonumber
\end{eqnarray}
Thus $\rho_g$ is \intro{unitary} with respect to the Clifford valued inner
product~\eqref{eq:f-pt}. Notice this notion is stronger than unitarity for the scalar
valued inner product, as the latter is the trace of the Clifford valued
one. A proof of unitarity could also consist of proving the action of the
Lie algebra is skew-symmetric, i.e.\ that for an element $b$ of the Lie
algebra and $f$ arbitrary
$$
\scalar{d\rho_bf}{f}=\scalar{f}{-d\rho_bf}.
$$
Here $d\rho_b$ is derived representation of $d$ for an element $b\in 
\algebra{g}_n$ of the Lie algebra of $\Space{G}{n}$. In the next subsection
we will need the explicit form of it. For the selected basis of 
$\algebra{g}_n$ we have:
\begin{eqnarray} 
[d\rho(T_j) f] (x) & = &  (0, 0,\ldots, 0,2e_1 f_j(x),0, \ldots ,0,0);
\nonumber \\
{} [d\rho (P) f] (x) & = & (- \sqrt[]{2} \frac{\partial}{\partial x}
f_1(x), \ldots, -\sqrt[]{2} \frac{\partial}{\partial x} f_j(x), \ldots , -
\sqrt[]{2} \frac{\partial}{\partial x}
f_n(x)); \nonumber \\
{} [d\rho (Q_j) f] (x) & = & ( 0,0, \ldots,0, \sqrt[]{2} e_j x f_j(x), 0,
\ldots , 0,0).
\nonumber
\end{eqnarray}
Particularly $d\rho(Q_j) d\rho( Q_k)=0$ for all $j \neq k$. This does not 
follow from the structure of $\Space{G}{}$ but is a feature of the 
described representation.
\begin{rem}
The group $\Space{G}{n}$ is called as ``a generalized Heisenberg group'' 
in~\cite{Kumahara97} where its induced representations are considered.
\end{rem}

\subsection{The wavelet transform for \Space{G}{n} }
In $\FSpace{R}{2}(\Space{R}{})$ we have the \Space{C}{n}-valued function
\begin{equation}
f_0(x)=(e^{-x^2/2},\ldots,e^{-x^2/2}),
\end{equation}
which which is the \intro{vacuum vector} in this case.  It is a zero
eigenvector for the operator
\begin{equation}
\label{eq:g-annihil}
a^-=d\rho(P)- \sum_{j=1}^n e_j d\rho(Q_j),
\end{equation}
which is \emph{the only annihilating operator} in this model. But we still
have $n$ creation operators:
\begin{equation}
a_k^+ = d\rho(P)- \sum_{j=1}^n (1-2\delta_{jk})e_j d\rho(Q_j)= a^{-} + 
2e_kd\rho(Q_k).
\end{equation}
While $a^-$ and $a_k^+$ look a little bit exotic for 
$\Space{G}{1}=\Space{H}{1}$ they are exactly the standard annihilation and 
creation operators.
Another feature of the representation is that the $a_k^+$ do not commute
with each other and have a non-trivial commutator with $a^-$:
\begin{displaymath}
[a_j^+, a_k^+]= 2e_kd\rho(T_k) -2e_j d\rho(T_j), \qquad 
[a_j^+,a^-]=-2e_jd\rho(T_j)
\end{displaymath}

We need the transforms of $f_0(x)$ under the action~\eqref{eq:G-rep}, i.e.
the \intro{coherent states} $f_g(x)=[\rho_gf_0](x)$ in this model:
\begin{eqnarray}
f_g(x)&=&(\ldots,
e^{e_j(2t_j+q_j(\sqrt[]{2}x-p))} e^{-(x-\sqrt[]{2}p)^2/2},\ldots)
\nonumber \\
&=&(\ldots,
e^{2e_jr_j-(p^2+q_j^2)/2} e^{-((p-e_j q_j)^2+x^2)/2+
\sqrt[]{2}(p-e_jq_j)x},\ldots) \nonumber \\
&=&(\ldots,
e^{2e_jt_j-z_j\bar{z}_j/2} e^{-(\bar{z}_j^2+x^2)/2+
\sqrt[]{2}\bar{z}_jx},\ldots) \nonumber \\
\end{eqnarray}
where $z_j=p+e_j q_j$, $\bar{z}_j=p-e_j q_j$.

Having defined coherent states we can introduce the \intro{wavelet
transform} $\oper{W}: \FSpace{R}{2}(\Space{R}{}) \rightarrow 
\FSpace{L}{\infty}(\Space{G}{n})$ by the standard formula:
\begin{equation} \label{eq:g-wt}
 \oper{W} f (g) = \scalar f{f_g}.
\end{equation}

Calculations completely analogous to those of the complex case allow us to
find the images $\oper Wf_{(t',\a)}(t,\z)$ of coherent states
$f_{(t',\a)}(x)$ under~\eqref{eq:g-wt} as follows:
\begin{eqnarray}
\oper Wf_{(t',\a)}(t,\z)&=&\scalar{f_{(t',\a)}}{f_{(t,\z)}} \nonumber\\
&=&\int_{\Space{R}{}} \sum_{j=1}^n
\exp\left(-2e_jt_j-\frac{z_j\bar{z}_j}2-\frac{{z}_j^2+x^2}2+
\sqrt[]{2}z_jx\right)\nonumber\\
&&\quad
\exp\left(+2e_jt'_j-\frac{a_j\bar{a}_j}2-\frac{\bar{a}_j^2+x^2}2
+\sqrt[]{2}\bar{a}_jx\right)\,dx \nonumber \\
&=&
\sum_{j=1}^n
\exp\left(-2e_jt_j-\frac{z_j\bar{z}_j}2+2e_jt'_j-\frac{a_j\bar{a}_j}2
+\bar{a}_j z_j\right)\nonumber\\
&&\qquad\times \int_{\Space{R}{}} \exp({-x^2+
2x\frac{z_j+\bar{a}_j}{\sqrt[]{2}}-\frac{(z_j+\bar{a}_j)^2}{2}})\,dx
\nonumber \\
&=&
\sum_{j=1}^n
\exp\left(-2e_j(t_j-t'_j)-\frac{z_j\bar{z}_j+a_j\bar{a}_j}2
+\bar{a}_j z_j\right)\nonumber\\
&&\qquad\times \int_{\Space{R}{}}
\exp(-(x-\frac{z_j+\bar{a}_j}{\sqrt{2}})^2)\,dx \nonumber \\
&=&
\sum_{j=1}^n
\exp\left(-2e_j(t_j-t'_j)-\frac{z_j\bar{z}_j+a_j\bar{a}_j}2
+\bar{a}_j z_j\right)
\end{eqnarray}
Here $a_j=a_0+e_j a_j$, $\bar{a}_j=a_0-e_ja_j$; $z_j$, $\bar{z}_j$ were
defined above. 

In this case all $\oper Wf_{(t',\a)}(t,\z)$ are \intro{monogenic}
functions with respect to the following Dirac operator:
\begin{equation} \label{eq:dirac-f}
\frac{\partial}{\partial p}-\sum_{j=1}^n e_j \frac{\partial}{\partial q_j}
+\frac 12\sum_{j=1}^n  (e_j p+ q_j)\frac{\partial}{\partial t_j},
\end{equation}
with $z_j$ related to $p$ and $q_j$ as above. This can be checked by the
direct calculation or follows from the observation: the 
Dirac operator~\eqref{eq:dirac-f} is the image of the annihilation operator 
$a^-$~\eqref{eq:g-annihil} under the wavelet transform~\eqref{eq:g-wt}. 
The situation is completely analogous to the Segal-Bargmann case, where
holomorphy is defined by the operators $\frac{\partial}{\partial
\bar{z}_k}$, which are the images of the annihilation operators $a^-_k$.
Actually, it is the Dirac operator associated with the unique left invariant
metric on $\Space Gn/\Space Z{}$ for which $P$ together with the $Q_k$
forms an orthonormal basis in the origin, and therefore everywhere.

The operator~\eqref{eq:dirac-f} is a realization of a generic Dirac
operator constructed for a nilpotent Lie group, see~\cite{ConMosc82}.
Indeed the operator~\eqref{eq:dirac-f} is defined by the formula $D= P + 
\sum_1^n e_j Q_j$, where $P$ and $Q_j$  are the left invariant vector fields 
in \eqref{eq:l-first}--\eqref{eq:l-last}. So the operator~\eqref{eq:dirac-f}
is left invariant and one has only to check the monogenicity of $ 
\oper Wf_{(0,0)}(t,\z) $---all other functions
$\oper Wf_{(t',\a)}(t,\z) $ are its left shifts.

Of course all linear combinations of the $\oper Wf_{(t',\a)}(t,\z)$ are
also monogenic.
So if we define two function spaces, $\FSpace{R}{2}$ and
$\FSpace{M}{2}$, as being the closure of the linear span of all $f_g(x)$ and
$\oper Wf_{(t',\a)}(t,\z)$ respectively, then
\begin{enumerate}
\item $\FSpace{M}{2}$ is a space of monogenic function on \Space{G}{n} in
the sense above.
\item \Space{G}{n} has representations both in $\FSpace{R}{2}$ and in
$\FSpace{M}{2}$. On the second space the group acts via left regular
representation.
\item These representation are intertwining by the integral
transformation with the kernel $T(t',\a,x)=f_{(t',\a)}(x)$.
\item The space $\FSpace{M}{2}$ has a reproducing kernel
$K(t',\a,t,z)=\oper Wf_{(t',\a)}(t,\z)$.
\end{enumerate}
The standard wavelet transform can be processed as expected. 

For the reduced wavelet transform associated with the mapping 
$s:\Omega \rightarrow \Space{G}{n}$ in particular we have
$$
\widehat{\oper W}f_\a(\z)=\oper Wf_{(0,\a)}(0,\z)=\sum_{j=1}^n\exp\bar{a}_j z_j.
$$

However the reduced wavelet transform cannot be constructed from a
single vacuum vector.  We need exactly $n$ linearly independent vacuum
vectors and the corresponding multiresolution wavelet analysis (wavelet
transform with several independent vacuum vectors) which is outlined
in~\cite{BratJorg97a} (see also M.G.~Krein's works~\cite{Krein48a} on
``directing functionals").  Indeed we have $n$ different vacuum vectors
$(\ldots, 0, e^{-x^2/2},0,\dots)$ each of which is an eigenfunction for the
action of the centre of $ \Space{G}{n}$.  All functions $\widehat{\oper
W}f_\a(\z)$ are \intro{monogenic} with respect to the Dirac operator
\begin{equation}
D=\frac{\partial}{\partial p}+ \sum_{j=1}^n e_j \frac{\partial}{\partial q_j}.
\end{equation}

\begin{rem}
  \index{quantum mechanics}
  It may seem on the first glance that the theory constructed in this
  section is only (if any) of pure mathematical interest and could not
  be related to physical reality. However the situation with an
  matching number of coordinates and momenta appears in very promising
  approach to \introind{quantum field theory}{field theory!quantum}, see~\cite{Kanatchikov01a} and
  references therein. This series of papers also come to the
  conclusion that a proper quantum picture in that setting requires
  Clifford algebras.
\end{rem}

\appendix

\chapter{Groups and Homogeneous Spaces}
\label{sec:groups-homog-spac}

The group theory and the representation theory are two enormous
and interesting subjects themselves. However they are auxiliary in our 
consideration and we are forced to restrict ourselves only to brief
and very dry overview.

Besides introduction to that areas presented
in~\cite{Miller68,Vilenkin68} we recommend additionally the
books~\cite{Kirillov76,MTaylor86}.  The representation theory
intensively uses tools of functional analysis and on the other hand
inspires its future development. We use the book~\cite{KirGvi82} for
references on functional analysis here and recommend it as a nice
reading too.

\section[Groups]{Basics of Group Theory}
\label{sec:basics-group-theory}
We start from the definition of central object which formalizes the
universal notion of symmetries.
\begin{defn} \label{de:group}
  A \introind{transformation group}{group!transformation}
  $G$ is a nonvoid set of mappings of a certain set $X$ into itself
  with the following properties:
  \begin{enumerate}
  \item if $g_1 \in G$ and $g_2 \in G$ then $g_1g_2\in G$;
  \item if $g\in G$ then $g^{-1}$ exists and belongs to $G$.
  \end{enumerate}
\end{defn}
\begin{exercise}
\label{exe:3-elem}
  List all transformation groups on a set of three elements.
\end{exercise}
\begin{exercise}
\label{exe:basic-examples}
  Verify that the following are groups in fact: 
  \begin{enumerate}
  \item Group of permutations of $n$ elements;
  \item Group of $n\times n$ matrixes with non zero determinant over a
    field $\Space{F}{}$ under matrix multiplications;
  \item Group of \extref{math255}{sec:rotations}{rotations} of the
    unit circle $\Space{T}{}$;
  \item Groups of \extref{math255}{sec:translations}{shifts} of the
    real line $\Space{R}{}$ and plane $\Space{R}{2}$;
  \item Group of linear fractional transformations of the extended
    complex plane.
  \end{enumerate}
\end{exercise}
\begin{defn}
  \label{de:abstract-group}
  An \introind{abstract group}{group!abstract} (or simply
  \intro{group}) is a nonvoid set $G$ on which there is a law of
  \intro{group multiplication} (i.e. mapping $G \times G\rightarrow
  G$) with the properties   
  \begin{enumerate}
  \item \intro{associativity}: $g_1(g_2g_3)=(g_1g_2)g_3$;
  \item the existence of \intro{identity}: $e\in G$ such that
    $eg=ge=g$ for all $g\in G$;
  \item the existence of \intro{inverse}: for every $g\in G$ there
    exists $g^{-1}\in G$ such that $g g^{-1}=g^{-1}g=e$.
  \end{enumerate}
\end{defn}
\begin{exercise}
  Check that any transformation group is an abstract group.
\end{exercise}
\begin{exercise}
  Check that the following transformation groups
  (cf. Example~\ref{exe:basic-examples}) have the same law of
  multiplication, i.e. are equivalent as abstract groups:
  \begin{enumerate}
  \item The group of isometric mapping of an equilateral triangle
    onto itself;
  \item The group of all permutations of a set of free elements;
  \item The group of invertible matrix of order $2$ with coefficients
    in the field of integers modulo $2$;
  \item The group of linear fractional transformations of the extended 
    complex plane generated by the mappings $z\mapsto z^{-1}$ and
    $z\mapsto 1-z$.
  \end{enumerate}
\end{exercise}
\begin{exercises}
  Expand the list in the above exercise. 
\end{exercises}
It is simpler to study groups with the following additional property.
\begin{defn}
  A group $G$ is \introind{commutative}{group!commutative} if for all
  $g_1$, $g_2\in G$, we have $g_1g_2=g_2g_1$.
\end{defn}
Most of interesting and important groups are
\introind{noncommutative}{group!noncommutative}, however.   
\begin{exercise}
  \begin{enumerate}
  \item Which groups among found in Exercise~\ref{exe:3-elem} are
    commutative?
  \item Which groups among listed in Exercise~\ref{exe:basic-examples} 
    are noncommutative?
\end{enumerate}
\end{exercise}
Groups could have some additional
\extref{math150}{sec:definition-limit}{analytical} structures,
e.g. they could be a topological sets with a corresponding notion of
\extref{math150}{sec:definition-limit}{limit}. We always assume that
our groups are \introind{locally compact}{group!locally
compact}~\cite[\S~2.4]{Kirillov76}.
\begin{defn}
  If for a group $G$ the group multiplication and the taking of
  inverse are \extref{math150}{de:continuity}{continuous} mappings
  then $G$ is \introind{continuous group}{group!continuous}.
\end{defn}
Even a better structure could be found among \introind{Lie
  groups}{group!Lie}~\cite[\S~6]{Kirillov76}, e.g. groups with
a \extref{math150}{de:derivative}{differentiable} law of
multiplication.  Investigating such groups we could employ the whole
arsenal of analytical tools, thereafter most of groups studied in this
notes will be Lie groups.
\begin{exercise}
  Check that the following are noncommutative Lie (and thus
  continuous) groups:
  \begin{enumerate}
  \item \label{item:ax+b-group} 
    \cite[Chap.~7]{MTaylor86} The \introind{$ax+b$
        group}{group!$ax+b$}: set of elements
    $(a,b)$, $a\in \Space[+]{R}{}$, $b\in \Space{R}{}$ with the group
    law:
    \begin{eqnordisp}
      (a, b) * (a', b') = (aa', ab'+b).
    \end{eqnordisp}
    The identity is $(1,0)$, and $(a,b)^{-1}=(a^{-1},-b/a)$.
  \item \label{item:Heisenberg-group}
    The \introind{Heisenberg group}{group!Heisenberg}~\cite{Howe80a},
      \cite[Chap.~1]{MTaylor86}: a set of
    triples of real numbers $(s,x,y)$ with the group multiplication:
    \begin{eqnordisp}[eq:heisenberg-def]
      (s,x,y)*(s',x',y')=(s+s'+\frac{1}{2}(x'y-xy'),x+x',y+y').
    \end{eqnordisp}
    The identity is $(0,0,0)$, and $(s,x,y)^{-1}=(-s,-x,-y)$.
  \item \label{item:SL-group}
    The \introind{$\SL$}{group!$\SL$}
    group~\cite{HoweTan92,Lang85}: a set of $2\times 2$ matrixes %
    \iftth\else %
    $\matr{a}{b}{c}{d}$ %
    \fi %
    with real entries $a$, $b$, $c$,
    $d\in\Space{R}{}$, the determinant $\object{det}=ad-bc$ equal to 1 
    and the group law coinciding with matrix multiplication:
    \begin{eqnordisp}
      \matr{a}{b}{c}{d}\matr{a'}{b'}{c'}{d'}=
      \matr{aa'+bc'}{ab'+bd'}{ca'+dc'}{cb'+dd'}.
    \end{eqnordisp}
    The identity is the unit matrix and
    \begin{eqnordisp}
      \matr{a}{b}{c}{d}^{-1} =\matr{d}{-b}{-c}{a}.
    \end{eqnordisp}
  \end{enumerate}
\end{exercise} 
The above three groups are behind many important
results of real and complex
analysis~\cite{Howe80a,HoweTan92,Lang85} and we meet them many
times in these notes.

\section[Homogeneous Spaces, Invariant 
Measures]{Homogeneous Spaces and Invariant Measures}
\label{sec:homog-spac-invar}

While \hyperref[de:abstract-group]{abstract group} are a suitable
language for investigation of their general properties we meet groups
in applications as \hyperref[de:group]{transformation groups} acting
on a set $X$.

Let $X$ be a set and let be defined an operation $G: X\rightarrow X$
of $G$ on $X$.  There is an
\hyperlink{equivalence-relation}{equivalence relation} on $X$, say,
$x_1\sim x_2 \Leftrightarrow \exists g\in G: gx_1=x_2$, with respect
to which $X$ is a disjoint union of distinct
\intro{orbits}~\cite[\S~I.5]{Lang69}. 
\begin{exercise}
  \label{item:SL-action-UHP} 
  Let action of \hyperref[item:SL-group]{$\SL$ group} on
  $\Space{C}{}$ by means of \introind{linear-fractional 
    transformations}{transformation!linear-fractional}:
  \begin{eqnordisp}
    \matr{a}{b}{c}{d}: z \mapsto \frac{az+b}{cz+d}.
  \end{eqnordisp}
  Show that there three orbits: the real axis $\Space{R}{}$,
  \introind{upper (lower) half plane}{lower (upper) half plane} 
  \intro{$\Space[\pm]{R}{n}$}:
  \begin{eqnordisp}
    \Space[\pm]{R}{n}=\{ x\pm iy \such x,y\in \Space{R}{},\ y>0\}.
  \end{eqnordisp}
\end{exercise}

Thus from now on, without lost
of a generality, we assume that the operation of $G$ on $X$ is
\introind{transitive}{action!transitive}, i.~e. for every $x\in X$ we
have
\begin{displaymath}
  Gx:=\relstack{\bigcup}{g\in G} g(x)=X.
\end{displaymath} 
In this case $X$ is $G$-\intro{homogeneous space}.
\begin{exercise}
  Show that for any group $G$ we could define its action on $X=G$ as
  follows: 
  \begin{enumerate}
  \item \label{item:repr-conjugation} 
    The \intro{conjugation} $g: x \mapsto g x g^{-1}$ (which is even a
    group homomorphism, but is trivial for all commutative groups).
  \item \label{item:repr-shifts}
    The \introind{left shift}{shift!left} $\lambda(g): x \mapsto g x$ and
    the \introind{right shift}{shift!right} $\rho(g): x \mapsto  x g^{-1}$.
  \end{enumerate}
\end{exercise}

If we fix a point $x\in X$ then the set of elements
$G_x=\{g\in G\such g(x)=x\}$ obviously forms the \introind{isotropy
(sub)group}{group!subgroup!isotropy} of $x$ in
$G$~\cite[\S~I.5]{Lang69}. The set $X$ is in the bijection with the
factor set $G/G_x$ for any $x\in X$.

\begin{exercise}
  Find a subgroup which correspond to the given action of $G$ on $X$:
  \begin{enumerate}
  \item Action of \hyperref[item:ax+b-group]{$ax+b$ group} on
    $\Space{R}{}$ by the formula: $(a,b): x \mapsto ax+b$.
  \item Action of \hyperref[item:SL-group]{$\SL$ group} on one of
    three orbit from Exercise~\ref{item:SL-action-UHP}.
  \end{enumerate}
\end{exercise}

To do some analysis on groups we need suitably defined basic
operation: \extref{math150}{de:derivative}{differentiation} and
\extref{math150}{sec:definite-integral}{integration}. The first
operation is naturally defined for Lie group.
If $G$ is a Lie group then the homogeneous
space $G/G_x$ is a smooth manifold (and a \intro{loop} as an
algebraic object) for every $x\in X$. Therefore the one-to-one mapping
$G/G_x \rightarrow X: g\mapsto g(x)$ induces a structure of
$C^\infty$-manifold on $S$. Thus the class $\FSpace[\infty]{C}{0}(X)$
of smooth functions with compact supports on $x$ has the evident
definition.

In order to perform an integration we need a suitable \intro{measure}.
A smooth measure $d\mu$ on $X$ is called (left) \introind{invariant
  measure}{measure!left invariant} with respect to an operation of $G$
  on $X$ if
\begin{equation}
  \int_X f(x)\, d\mu(x) = \int_X f(g(x))\, d\mu(x),\quad
  \textrm{for all } g\in G,\ 
  f(x)\in\FSpace[\infty]{C}{0}(X).\label{eq:invar-m}
\end{equation} 
\begin{exercise}
  \label{th:SL2-measure-invariant}
  Show that measure $y^{-2}dy\,dx$ on the upper half plane
  $\Space[+]{R}{2}$ is invariant under action from
  Exercise~\ref{item:SL-action-UHP}.
\end{exercise} 
\hypertarget{Haar-measure}{Left invariant measures} on
$X=G$ is called the \introind{Haar measure}{measure!Haar}. It always
exists and is uniquely defined up to a scalar
multiplier~\cite[\S~0.2]{MTaylor86}. An equivalent formulation
of~\eqref{eq:invar-m} is: \emph{$G$ operates on $\FSpace{L}{2}(X,d\mu)$
by
  unitary operators}. We will transfer the Haar measure
$d\mu$ from $G$ to $\algebra{g}$ via the exponential map $\exp:
\algebra{g}\rightarrow G$ and will call it as the \emph{invariant
measure on a Lie algebra \algebra{g}}.
\begin{exercise}
  Check that the following are Haar measures for corresponding groups:
  \begin{enumerate}
  \item \label{item:R-measure-invariant}
    The \introind{Lebesgue measure}{measure!Lebesgue} $dx$ on the
    real line $\Space{R}{}$.
  \item \label{item:SO-measure-invariant} 
    The Lebesgue measure $d\phi$ on the
    unit circle $\Space{T}{}$.
  \item $dx/x$ is a Haar measure on the multiplicative group
    $\Space[+]{R}{}$;
  \item $dx\,dy/(x^2+y^2)$ is a Haar measure on the multiplicative
    group $\Space{C}{}\setminus \{0\}$, with coordinates $z=x+iy$.
  \item \label{item:ax+b-measure-invariant}
    $a^{-2}\,da\,db$ and
    $a^{-1}\,da\,db$ are the left and right invariant measure on
    \hyperref[item:ax+b-group]{$ax+b$ group}. 
  \item \label{item:H-measure-invariant}
    The Lebesgue measure $ds\,dx\,dy$ of $\Space{R}{3}$ for the
    \hyperref[item:Heisenberg-group]{Heisenberg group} $\Space{H}{1}$.
  \end{enumerate}
\end{exercise} 
In this notes we assume \emph{all integrations on
groups performed over the Haar measures}.
\begin{exercise}
  \label{ex:compact-measure}
  Show that invariant measure on a compact group $G$ is finite and
  thus may be normalized to total measure $1$.
\end{exercise}
The above simple result has surprisingly
\hyperref[th:compact-represent]{important consequences}. 
\begin{defn}
  \label{de:convolution}
  The left \intro{convolution} $f_1*f_2$ of two functions $f_1(g)$ and
  $f_2(g)$ defined on a group $G$ is
  \begin{eqnordisp}
    f_1*f_2(g)=\int_G f_1(h)\,f_2(h^{-1}g)\,dh
  \end{eqnordisp}
\end{defn}
\begin{exercise}
  Let $k(g)\in \FSpace{L}{1}(G,d\mu)$ and operator $K$ on
  $\FSpace{L}{1}(G,d\mu)$ is the left \introind{convolution
    operator}{operator!convolution}  with $k$, .i.e. $K: f \mapsto
  k*f$. Show that $K$ commutes with all
  \hyperref[item:repr-shifts]{right shifts} on $G$.  
\end{exercise}

The following Lemma characterizes {\em linear
subspaces\/} of $\FSpace{L}{2}(G,d\mu)$ invariant under shifts in the term of {\em ideals of convolution
algebra\/} $\FSpace{L}{2}(G,d\mu)$ and is of the separate
interest.

\begin{lem} 
  \label{le:ideal}
  A closed linear subspace $H$ of $\FSpace{L}{2}(G,d\mu)$ is invariant
  under left (right) shifts if and only if $H$ is a left (right) ideal of
  the right group convolution algebra $\FSpace{L}{2}(G,d\mu)$.
  
  A closed linear subspace $H$ of $\FSpace{L}{2}(G,d\mu)$ is invariant
  under left (right) shifts if and only if $H$ is a right (left) ideal of
  the left group convolution algebra $\FSpace{L}{2}(G,d\mu)$.
\end{lem}
\begin{proof} Of course we consider only the ``right-invariance and
  right-convolution'' case. Then the other three 
  cases are analogous. Let $H$ be a closed linear subspace of
  $\FSpace{L}{2}(G,d\mu)$ invariant under right shifts and $k(g)\in H$. We
  will show the inclusion
  \begin{equation}\label{eq:rt-convol}
    [f*k]_r(h)=\int_G f(g)k(hg)\,d\mu(g)\in H,
  \end{equation}
  for any $f\in\FSpace{L}{2}(G,d\mu)$. Indeed, we can treat
  integral~\eqref{eq:rt-convol} as a limit of sums
  \begin{equation}\label{eq:rt-sum}
    \sum_{j=1}^{N} f(g_j)k(hg_j)\Delta_j.
  \end{equation}
  But the last sum is simply a linear combination of vectors $k(hg_j)\in
  H$ (by the invariance of $H$) with coefficients $f(g_j)$. Therefore
  sum~\eqref{eq:rt-sum} belongs to $H$ and this is true for
  integral~\eqref{eq:rt-convol} by the closeness of $H$.
  
  Otherwise, let $H$ be a right ideal in the group convolution algebra
  $\FSpace{L}{2}(G,d\mu)$ and let $\phi_j(g)\in\FSpace{L}{2}(G,d\mu)$ be
  an approximate unit of the algebra~\cite[\S~13.2]{Dixmier69}, i.~e. for
  any $f\in\FSpace{L}{2}(G,d\mu)$ we have
  \begin{displaymath}
    [\phi_j*f]_r(h)=\int_G \phi_j(g)f(hg)\, d\mu(g) \rightarrow f(h)\mbox{,
      when } j\rightarrow\infty.
  \end{displaymath}
  Then for $k(g)\in H$ and for any $h'\in G$ the right convolution
  \begin{displaymath}
    [\phi_j*k]_r(hh')=\int_G \phi_j(g)k(hh'g)\, d\mu(g)= \int_G
    \phi_j(h'^{-1}g')k(hg')\, d\mu(g'),\ g'=h'g,
  \end{displaymath}
  from the first expression is tensing to $k(hh')$ and from the second
  one belongs to $H$ (as a right ideal). Again the closeness of $H$
  implies $k(hh')\in H$ that proves the assertion.
\end{proof}

\chapter[Representation Theory]{Elements of the Representation Theory}
\label{sec:elem-repr-theory}

\section[Representations]{Representations of Groups}
\label{sec:repr-groups}
Objects unveil their nature in actions. Groups act on other sets by
means of
\introind{representations}{group!representation}\index{representation}.
A representation of a group $G$ is a group homomorphism of $G$ in a
transformation group of a set. It is a fundamental observation that
\emph{linear} objects are easer to study. Therefore we begin from
linear representations of groups.
\begin{defn}
  A linear continuous \intro{representation of a group}
  \index{representation!linear}\index{representation!continuous} $G$
  is a continuous function $T(g)$ on $G$ with values in the group of
  non-degenerate linear continuous transformation in a linear space
  $H$ (either finite or infinite dimensional) such that $T(g)$
  satisfies to the functional identity:
  \begin{equation}
    \label{eq:repr-defn}
    T(g_1 g_2) =T(g_1)\, T(g_2).
  \end{equation}
\end{defn}
\begin{exercise}
  Show that $T(g^{-1})=T^{-1}(g)$ and $T(e)=I$, where $I$ is the
  identity operator on $B$.
\end{exercise}
\begin{exercise}
  \label{exe:ax+b-repr}
  Show that these are linear continuous representations of
  corresponding groups:
  \begin{enumerate}
  \item Operators $T(x)$ such that $[T(x)\,f](t)=f(t+x)$ form
    a representation of $\Space{R}{}$ in $\FSpace{L}{2}(\Space{R}{})$.  
  \item Operators $T(n)$ such that $T(n) a_k=a_{k+n}$ form a
    representation of $\Space{Z}{}$ in $\FSpace{l}{2}$. 
  \item \label{item:ax+b-repr} Operators $T(a,b)$ defined by
    \begin{eqnordisp}[eq:ax+b-repr]
      [T(a,b)\, f](x)= \sqrt{a}f(ax+b), \qquad a \in \Space[+]{R}{},\ 
      b\in\Space{R}{}
    \end{eqnordisp}
    form a representation of \hyperref[item:ax+b-group]{$ax+b$ group} in
    $\FSpace{L}{2}(\Space{R}{})$.
  \item Operators $T(s,x,y)$ defined by
    \begin{eqnordisp}[eq:Schrodinger-repr]
      [T(s,x,y)\, f] (t)=e^{i(2s-\sqrt{2}yt+xy)} f(t- \sqrt{2}x)
    \end{eqnordisp}
    form \introind{Schr\"odinger
      representation}{representation!Schr\"odinger}
    \index{Heisenberg group!Schr\"odinger representation} of the
    \hyperref[item:Heisenberg-group]{Heisenberg group $\Space{H}{1}$}
    in $\FSpace{L}{2}(\Space{R}{})$.
  \item Operators $T(g)$ defined by
    \begin{eqnordisp}[eq:sl2R-repres]
      [T(g) f](t) = \frac{1}{ct+d} f\left( \frac{at+b}{ct+d}\right),
      \quad \textrm{ where } g=\matr{a}{b}{c}{d},
    \end{eqnordisp}
    form a representation of \hyperref[item:SL-group]{\SL} in
    $\FSpace{L}{2}(\Space{R}{})$. 
\end{enumerate}
\end{exercise}
In the sequel a representation \emph{always means} linear continuous
representation. $T(g)$ is an \introind{exact
  representation}{representation!exact} (or \introind{faithful
  representation}{representation!faithful} if $T(g)=I$ only for 
$g=e$. The opposite case when $T(g)=I$ for all $g\in G$ is a
\introind{trivial representation}{representation!trivial}. The space
$H$ is \intro{representation space} and in most cases will be a
\emph{Hilber space}~\cite[\S~III.5]{KirGvi82}. If dimensionality of $H$
is finite then $T$ is a \introind{finite dimensional representation}%
{representation!finite dimensional}, in the opposite case it is
\introind{infinite dimensional representation}%
{representation!infinite dimensional}.

We denote the \intro{scalar product} on $H$ by
\intro{$\scalar{\cdot}{\cdot}$}. Let $\{\mvec{e}_j\}$ be an
(finite or infinite) \intro{orthonormal basis} in $H$, i.e.
\begin{eqnordisp}
  \scalar{\mvec{e}_j}{\mvec{e}_j}=\delta_{jk},
\end{eqnordisp}
where \intro{$\delta_{jk}$} is the \intro{Kroneker delta}, and linear
span of $\{\mvec{e}_j\}$ is dense in $H$. 

\begin{defn}
  \label{de:matrix-elements}
  The \intro{matrix elements} $t_{jk}(g)$ of a representation $T$ of a
  group $G$ (with respect to a basis $\{\mvec{e}_j\}$ in $H$) are
  complex valued functions on $G$ defined by
  \begin{eqnordisp}[eq:matrix-elements]
    t_{jk}(g) = \scalar{T(g)\mvec{e}_j}{\mvec{e}_k}.
  \end{eqnordisp}
\end{defn}
\begin{exercise}
  Show that \cite[\S~1.1.3]{Vilenkin68}
  \begin{enumerate}
  \item $T(g)\,\mvec{e}_k=\sum_j t_{jk}(g)\,\mvec{e}_j$.
  \item 
    \label{item:matrix-addition}
    $t_{jk}(g_1g_2)=\sum_n t_{jn}(g_1)\,t_{nk}(g_2)$.
\end{enumerate}
\end{exercise}

It is typical mathematical questions to determine identical objects
which may have a different appearance. For representations it is
solved in the following definition.
\begin{defn}
  \label{de:repr-equivalent}
  Two representations $T_1$ and $T_2$ of the same group $G$ in spaces
  $H_1$ and $H_2$ correspondingly are \introind{equivalent
    representations}{representations!equivalent} if there exist a
  linear operator $A: H_1 \rightarrow H_2$ with the continuous inverse
  operator $A^{-1}$ such that:
  \begin{eqnordisp}
    T_2(g)= A\, T_1(g)\, A^{-1}, \qquad \forall g\in G.
  \end{eqnordisp}
\end{defn}
\begin{exercise}
  \label{exe:ax+b-repr1}
  Show that representation $T(a,b)$ of
  \hyperref[item:ax+b-group]{$ax+b$ group} in
  $\FSpace{L}{2}(\Space{R}{})$ from Exercise~\ref{item:ax+b-repr} is 
  equivalent to the representation
  \begin{eqnordisp}[eq:ax+b-repr-1]
    [T_1(a,b)\,f] (x)= \frac{e^{i\frac{b}{a}}}{\sqrt{a}}\,
    f\left(\frac{x}{a}\right).
  \end{eqnordisp}
\end{exercise}
\begin{proof}[Hint]
  Use the Fourier transform.
\end{proof} 

The \hypertarget{equivalence-relation}{\intro{relation of
equivalence}} is reflexive, symmetric, and transitive. Thus it splits
the set of all representations of a group $G$ into \intro{classes of
equivalent representations}. \emph{In the sequel we study group
representations up to their equivalence classes only}.
\begin{exercise}
  Show that equivalent representations have the same
  \hyperref[de:matrix-elements]{matrix elements} 
  in appropriate basis.
\end{exercise}

\begin{defn}
  \label{de:adjoint-repr}
  Let $T$ is a representation of a group $G$ in $H$
  The \introind{adjoint representation}{representation!adjoint}
  $T'(g)$ of $G$ in $H$ is defined by
  \begin{eqnordisp}
    T'(g)=\left( T(g^{-1})\right)^*,
  \end{eqnordisp}
  where $^*$ denotes the adjoint operator in $H$.
\end{defn}
\begin{exercise}
  Show that 
  \begin{enumerate}
  \item $T'$ is indeed a representation.
  \item $t'_{jk}(g)=\bar{t}_{kj}(g^{-1})$.
\end{enumerate}
\end{exercise}

Recall~\cite[\S~III.5.2]{KirGvi82} that a bijection $U: H \rightarrow H$
is a \introind{unitary operator}{operator!unitary} if
\begin{eqnordisp}
  \scalar{Ux}{Uy}=\scalar{x}{y}, \qquad \forall x, y \in H.
\end{eqnordisp}
\begin{exercise}
  Show that $UU^*=I$.
\end{exercise}
\begin{defn}
  \label{de:unitary-repr}
  $T$ is a \introind{unitary representation}{representation!unitary} of a 
  group $G$ in a space $H$ if $T(g)$ is a unitary operator for all
  $g\in G$. $T_1$ and $T_2$ are \introind{unitary equivalent
    representations}{representations!equivalent!unitary} if
  $T_2=UT_2U^{-1}$ for a unitary operator $U$. 
\end{defn}
\begin{exercise}
  \begin{enumerate}
  \item Show that all representations from Exercises~\ref{exe:ax+b-repr}
    are unitary.  
  \item Show that representations from Exercises~\ref{item:ax+b-repr}
    and~\ref{exe:ax+b-repr1} are unitary equivalent.
\end{enumerate}
\end{exercise}
\begin{proof}[Hint]
  Take that the Fourier transform is unitary for granted. 
\end{proof}
\begin{exercise}
  \label{ex:Lie-group-alg-repr}
  Show that if a Lie group $G$ is represented by unitary operators in $H$
  then its Lie algebra $\algebra{g}$ is represented by self-adjoint
  (possibly unbounded) operators in $H$.
\end{exercise}

The following definition have a sense for \emph{finite} dimensional
representations.  
\begin{defn}
  \label{def:chracter-of-representation}
  A \introind{character of representation}{character!of
    representation} $T$ is equal $\chi(g)= \tr (T(g))$, where
    \intro{$\tr$} is the \intro{trace}~\cite[\S~III.5.2
    (Probl.)]{KirGvi82} of operator.
\end{defn}
\begin{exercise}
  Show that 
  \begin{enumerate}
  \item Characters of a representation $T$ are constant on the
    \hyperref[item:repr-conjugation]{adjoint elements}     
    $g^{-1}hg$, for all $g\in G$.
  \item Character is an algebra homomorphism from an algebra of
    representations with Kronecker's (tensor)
    multiplication~\cite[\S~1.9]{Vilenkin68} to complex numbers.     
  \end{enumerate}
\end{exercise}
\begin{proof}[Hint]
  Use that $\tr(AB)=\tr(BA)$, $\tr(A+B)=\tr A + \tr B$, and $\tr ( A
  \otimes B)= \tr A\,\tr B$.
\end{proof}
For \emph{infinite} dimensional representation characters could be
defined either as distributions~\cite[\S~11.2]{Kirillov76} or in
infinitesimal terms of Lie algebras~\cite[\S~11.3]{Kirillov76}.

The characters of a representation should not be confused with the
following notion.
\begin{defn}
  \label{def:character-of-group}
  A \introind{character of a group}{character!of group} $G$ is a 
  one-dimensional   representation of $G$.
\end{defn}
\begin{exercise}
  \begin{enumerate}
  \item Let $\chi$ be a \hyperref[def:character-of-group]{character of a
      group} $G$. Show that a
    \hyperref[def:chracter-of-representation]{character of
      representation} $\chi$ coincides with it and thus is a character of
    $G$.
  \item \label{item:matrix-character}
    A \hyperref[de:matrix-elements]{matrix element} of a group
    character $\chi$ coincides with $\chi$. 
  \item \label{item:character-group}\index{group!of characters}
    Let $\chi_1$ and $\chi_1$ be
    \hyperref[def:character-of-group]{characters of a 
      group} $G$. Show that $\chi_1 \otimes \chi_2 =\chi_1\chi_2$ and
    $\chi'(g)=\chi_1(g^{-1})$ are
    again characters of $G$. In other words \emph{characters of a
      group form a group themselves}.

\end{enumerate}
\end{exercise}

\section{Decomposition of Representations}
\label{sec:decomposition}
The important part of any mathematical theory is classification
theorems on structural properties of objects. Very well known examples 
are:
\begin{enumerate}
\item The main theorem of arithmetics on unique representation an
  integer as a product of powers of prime numbers.
\item Jordan's normal form of a matrix.
\end{enumerate} 
The similar structural results in the representation
theory are very difficult. The easiest (but still rather difficult)
questions are on classification of \hyperref[de:unitary-repr]{unitary
representations} up to \hyperref[de:unitary-repr]{unitary
equivalence}.
\begin{defn}
  \label{th:invariant-subspaces}
  Let $T$ be a representation of $G$ in $H$. A linear subspace
  $L\subset H$ is \introind{invariant subspace}{subspace!invariant}
  for $T$ if for any $\mvec{x}\in L$ and any $g\in G$ the vector
  $T(g)\mvec{x}$ again belong to $L$.  
\end{defn}
There are always two trivial invariant subspaces: the null and entire
$H$. All other are \intro{nontrivial invariant subspaces}.
\begin{defn}
  \label{de:irreducible}
  If there are only two trivial invariant subspaces then $T$ is
  \introind{irreducible
    representation}{representation!irreducible}. In the opposite case
  we have \introind{reducible
    representation}{representation!reducible}. 
\end{defn} 
For any nontrivial invariant subspace we could define the
\intro{restriction of representation} of $T$ on it. In this way we
obtain a \intro{subrepresentation} of $T$. 
\begin{example}
  Let $T(a)$, $a\in\Space[+]{R}{}$ be defined as follows:
  $[T(a)]f(x)=f(ax)$. Then spaces of
  \extref{math150}{se:incr-funct}{even and odd functions} are
  invariant.
\end{example}
\begin{defn}
  \label{de:cyclic-vector}
  If the closure of liner span of all vectors $T(g) v$ is dense in $H$
  then $v$ is called \introind{cyclic vector}{vector!cyclic} for $T$.
\end{defn}
\begin{exercise}
  Show that for an irreducible representation any non zero vector is
  cyclic. 
\end{exercise}
The important property of unitary representation is complete
reducibility.
\begin{exercise}
  Let a unitary representation $T$ has an invariant subspace $L\subset 
  H$, then its orthogonal completion $L^\perp$ is also invariant.
\end{exercise}
\begin{thm}
  \label{th:comple-reducibility}
  \cite[\S~8.4]{Kirillov76}
  Any unitary representation $T$ of a locally compact group $G$
  could be decomposed in a (continuous) direct sum irreducible
  representations: $T=\int_X T_x \,d\mu(x)$.
\end{thm}
The necessity of continuous sums appeared in very simple examples: 
\begin{exercise}
  \label{ex:decomp-multipl}
  Let $T$ be a representation of $\Space{R}{}$ in
  $\FSpace{L}{2}(\Space{R}{})$ as follows: $[T(a)f](x)=e^{iax}f(x)$. 
  Show that
  \begin{enumerate}
  \item Any measurable set $E\subset \Space{R}{}$ define an invariant
    subspace of functions vanishing outside $E$.
  \item $T$ does not have invariant irreducible subrepresentations.
  \end{enumerate}
\end{exercise}
\begin{defn}
  The set of equivalence classes of unitary irreducible
  representations of a group $G$ is denoted be $\hat{G}$ and called
  \intro{dual object} (or \intro{dual space}) of the group $G$.
\end{defn}

\begin{defn}
  \label{de:regular-representation} A left \introind{regular
  representation}{representation!regular} \intro{$\Lambda(g)$} of a group
  $G$ is the representation by \hyperref[item:repr-shifts]{left shifts} in
  the space $\FSpace{L}{2}(G)$ of square-integrable function on $G$
  with the left \hyperlink{Haar-measure}{Haar measure}
  \begin{eqnordisp}[eq:-left-reg-repr]
    \Lambda{g}: f(h) \mapsto f(g^{-1}h).
  \end{eqnordisp}
  The \intro{main
    problem of representation theory} is to decompose a left regular
  representation $\Lambda(g)$ into irreducible components.
\end{defn}

\section[Schur's Lemma]{Invariant Operators and Schur's Lemma} 

It is a pleasant feature of an abstract theory that we obtain
important general statements from simple
observations. \hyperref[ex:compact-measure]{Finiteness of invariant
measure} on a compact group is one such example. Another example is
\hyperref[le:schur]{Schur's Lemma} presented here. 

To find different classes of representations we need to compare them
each other. This is done by \emph{intertwining operators}.
\label{sec:invar-oper-schur}
\begin{defn}
  \label{de:intertwining}
  Let $T_1$ and $T_2$ are representations of a group $G$ in a spaces
  $H_1$ and $H_2$ correspondingly. An operator $A: H_1 \rightarrow
  H_2$ is called an \introind{intertwining
    operator}{operator!intertwining} if
  \begin{eqnordisp}
    A\,T_1(g) = T_2(g)\,A, \qquad \forall g\in G.
  \end{eqnordisp}
  If $T_1=T_2=T$ then $A$ is \emph{interntwinig operator} or
  \introind{commuting operator}{operator!commuting} for $T$.
\end{defn}
\begin{exercise}
  Let $G$, $H$, $T(g)$, and $A$ be as above. Show that
  \cite[\S~1.3.1]{Vilenkin68}
  \begin{enumerate}
  \item Let $\mvec{x}\in H$
    be an eigenvector for $A$ with eigenvalue $\lambda$. Then
    $T(g)\mvec{x}$ for all $g\in G$ are eigenvectors of $A$ with the
    same eigenvalue $\lambda$.
  \item All eigenvectors of $A$ with a fixed eigenvalue $\lambda$ for
    a linear subspace invariant under all $T(g)$, $g\in G$.
  \item If an operator $A$ is commuting with
    \hyperref[de:irreducible]{irreducible representation} $T$ then
    $A=\lambda I$. 
  \end{enumerate}
\end{exercise}
\begin{proof}[Hint]
  Use the spectral decomposition of selfadjoint
  operators~\cite[\S~V.2.2]{KirGvi82}.
\end{proof}
The next result have very important applications.
\begin{lem}[Schur]
  \label{le:schur}
  \index{Schur's lemma}\index{lemma!Schur's}
  \cite[\S~8.2]{Kirillov76}
  If two representations $T_1$ and $T_2$ of a group $G$ are
  irreducible, then every \hyperref[de:intertwining]{intertwining
    operator} between them either zero or is invertible.
\end{lem}
\begin{proof}[Hint]
  Consider subspaces $\ker A\subset H_1$ and $\object{im}A\subset H_2$.
\end{proof}
\begin{exercise}
  Show that
  \begin{enumerate}
  \item Two irreducible representations either equivalent or
    disjunctive. 
  \item All operators commuting with an irreducible representation
    form a field.
  \item \label{item:comm-represent}
    Irreducible representation of commutative group are
    one-dimensional. 
  \item If $T$ is unitary irreducible representation in $H$ and
    $B(\cdot,\cdot)$ is a bounded semi linear form in $H$ invariant
    under $T$: $B(T(g)\mvec{x},T(g)\mvec{y})=B(\mvec{x},\mvec{y})$
    then $B(\cdot,\cdot)=\lambda\scalar{\cdot}{\cdot}$.
  \end{enumerate}
\end{exercise}
\begin{proof}[Hint]
  Use that $B(\cdot,\cdot)=\scalar{A\cdot}{\cdot}$ for some
  $A$~\cite[\S~III.5.1]{KirGvi82}. 
\end{proof}

\chapter{Miscellanea}
\label{sec:miscellanea}

\section{Functions of even 
Clifford numbers}
\label{pt:bivec-fun}
Let 
\begin{equation} \label{eq:p-def}
\n{a}= a_1 \n{p}_1 + a_2 \n{p}_2, \qquad \n{p}_1= \frac{1+e_1e_2}{2},\quad 
\n{p}_2= \frac{1-e_1e_2}{2},\quad a_1,a_2\in \Space{R}{}
\end{equation}
be an even Clifford number in $ \Cliff[1]{1} $. It follows from the 
identities
\begin{equation} \label{eq:p-prop}
\n{p}_1 \n{p}_2 = \n{p}_2 \n{p}_1  =0, \qquad \n{p}_1^2=\n{p}_1, \qquad
\n{p}_2^2=\n{p}_2, \qquad \n{p}_1 + \n{p}_2=1
\end{equation}
that $p( \n{a})= p(a_1) \n{p}_1 + p(a_2) \n{p}_2$ for any polynomial
$p(x)$. Let $\FSpace{P}{}$ be a topological space of functions 
$ \Space{R}{} \rightarrow \Space{R}{} $ such that polynomials are dense in 
it. Then for any $ f\in \FSpace{P}{} $ we can define $f(\n{a})$ by the 
formula
\begin{equation} \label{eq:bivec-fun}
f(\n{a})= f(a_1) \n{p}_1 + f(a_2) \n{p}_2.
\end{equation}
This definition gives continuous algebraic homomorphism. 

\section{Principal series representations of $\SL$}
\label{pt:principal}
We describe a realization of the principal series representations of 
$\SL$. The realization is deduced from the realization by left regular 
representation on the a space of homogeneous function of power $-is-1$ on $ 
\Space{R}{2} $ described in~\cite[\S~8.3]{MTaylor86}. We consider now the 
restriction of homogeneous function not to the unit circle as 
in~\cite[Chap.~8, (3.23)]{MTaylor86} but to the line $x_2=1$ in $ 
\Space{R}{2} $. Then an equivalent unitary representation of $\SL$ acts on 
the Hilbert space $\FSpace{L}{2}( \Space{R}{} )$ with the standard Lebesgue
measure by the transformations:
\begin{equation} \label{eq:principal}
[\pi_{is} (g) f](x)= \frac{1}{ \modulus{cx+d}^{1+is} } f \left( 
\frac{ax+b}{cx+d} \right), \qquad g^{-1}= \matr{a}{b}{c}{d}.
\end{equation}

\section{Boundedness of the Singular Integral Operator 
$\oper{W}_\sigma$} \label{pt:sio}
The kernel of integral operator $\oper{W}_\sigma$~\eqref{eq:sio1} is
singular in four points, which are the intersection of $\TSpace{T}{}$ and
the light cone with the origin in $\n{u}$. One can easily see
\begin{displaymath}
\modulus{\frac{(-\n{u} e_1 e^{e_{12} t} + \n{1})^\sigma}
{(-e^{-e_{12} t} e_1 \n{u} + \n{1})^{1+\sigma}}} =
\modulus{1+\n{u}^2}^{1/2} \frac{1}{ \modulus{t-t_0} }  + O(\frac{1}{ 
\modulus{t-t_0}^2 }).
\end{displaymath}
where $t_0$ is one of four singular points mentioned before for a fixed
$\n{u}$ and $t$ is a point in its neighborhood.  More over the kernel of
integral operator $\oper{W}_\sigma$ is changing the sign while $t$ crossing
the $t_0$.  Thus we can define $\oper{W}_\sigma$ in the sense of the
principal value as the standard singular integral operator.

Such defined integral operator $\oper{W}_\sigma$ becomes a bounded linear
operator $ \FSpace{L}{2} ( \TSpace{T}{} ) \rightarrow \FSpace{L}{2} (
\TSpace{T}{\lambda}) $, where $\TSpace{T}{\lambda}$ is the
circle~\eqref{eq:circle-l} in $\TSpace{R}{1,1}$ with center in the origin
and the ``radius'' $\lambda$.  Moreover the norm of the operator
$\lambda^{-2}\oper{W}_\sigma$ is uniformly bounded for all $\lambda$ and
thus we can consider it as bounded operator
\begin{displaymath}
\FSpace{L}{2} ( \TSpace{T}{} ) \rightarrow \FSpace{H}{\sigma} ( 
\TSpace{D}{} ) ,
\end{displaymath}
where
\begin{equation}
\FSpace{H}{\sigma} ( \TSpace{D}{} ) =\{ f(\n{u}) \such D_{\TSpace{D}{}} f(\n{u} 
)=0,\  \n{u}\in \TSpace{D}{},\
\modulus{\lambda}^{-2}\int_{\TSpace{T}{\lambda}} \modulus{f(\n{u})}^2\, du 
< \infty,\ \forall \lambda<0\}. \label{eq:hardy-b}
\end{equation}
is an analog of the classic Hardy space. Note that 
$\modulus{\lambda}^{-2}d\n{u}$ is exactly the invariant
measure~\eqref{eq:def-m-b} on $\TSpace{D}{}$. 

One can note the similarity of arising divergency and singularities with
the ones arising in \introind{quantum field theory}{field
  theory!quantum}. The similarity generated by the  
same mathematical object in basement: a pseudoeuclidean space with 
an indefinite metric.

\renewcommand{\thesection}{}

\bibliographystyle{plain}
\bibliography{abbrevmr,akisil,analyse,aphysics,acombin,arare}

\pdfbookmark[1]{Index}{index}
\printindex

\end{document}